# Split buildings of type $\mathsf{F}_4$ in buildings of type $\mathsf{E}_6$

Anneleen De Schepper[1] · N. S. Narasimha Sastry[2] · Hendrik Van Maldeghem[1]

**Abstract** A symplectic polarity of a building $\Delta$ of type $\mathsf{E}_6$ is a polarity whose fixed point structure is a building of type $\mathsf{F}_4$ containing residues isomorphic to symplectic polar spaces (i.e., so-called *split buildings* of type $\mathsf{F}_4$). In this paper, we show in a geometric way that every building of type $\mathsf{E}_6$ contains, up to conjugacy, a unique class of symplectic polarities. We also show that the natural point-line geometry of each split building of type $\mathsf{F}_4$ fully embedded in the natural point-line geometry of $\Delta$ arises from a symplectic polarity.

## Contents

Communicated by Ingo Runkel.

The first author is supported by the Fund for Scientific Research—Flanders (FWO—Vlaanderen).

The original version of this article was revised.

✉ Hendrik Van Maldeghem
  Hendrik.VanMaldeghem@UGent.be

  Anneleen De Schepper
  Anneleen.DeSchepper@UGent.be

  N. S. Narasimha Sastry
  nnsastry@gmail.com

[1] Department of Mathematics, Ghent University, Krijgslaan 281-S22, 9000 Ghent, Belgium

[2] Indian Institute of Technology, Dharwad, Walmi Campus, Near High Court, Dharwad 580011, India

98 A. De Schepper et al.

## 1 Introduction

Buildings are the natural geometries of algebraic groups and certain variants of them such as groups of mixed type and classical groups over division rings which are infinite dimensional over their centre. Galois descent and forms translate to fixed point buildings of an automorphism group. For the fundamental, so called *split* buildings over an algebraically closed field, this automorphism group is trivial. However, also these split buildings can arise as fixed point buildings under a suitable nontrivial automorphism group (but not a Galois group, of course) of other split buildings. There are four special examples of this phenomenon, and they beautifully fit together (other examples exist in abundance, e.g., a central collineation in a polar space fixes a polar space of one rank less). Each of them occurs over an arbitrary field. In each of those examples, the ambient building has a simply laced Dynkin diagram, and the fixed point building has a double or triple bond in its Dynkin diagram. The four examples can be listed in increasing complexity:

(1) A symplectic polarity (i.e., the polarity related to a nondegenerate alternating form on the underlying vector space) of a projective space of dimension $2n-1$ (i.e., a split building of type $\mathsf{A}_{2n-1}$) fixes a symplectic polar space (i.e., a split building of type $\mathsf{C}_n$);
(2) A non-type preserving involutory automorphism of the oriflamme complex of a hyperbolic quadric (i.e., a building of type $\mathsf{D}_{n+1}$) which pointwise fixes a hyperplane of the ambient projective space fixes a parabolic polar space (i.e., a split building of type $\mathsf{B}_n$);
(3) A triality of type $I_{\mathrm{id}}$ (in the terminology of [21]) of the oriflamme complex of a hyperbolic quadric with Witt index 4 (i.e., a split building of type $\mathsf{D}_4$) fixes a split Cayley generalised hexagon (i.e., a building of exceptional type $\mathsf{G}_2$);
(4) A symplectic polarity of a building of exceptional type $\mathsf{E}_6$ fixes a split building of exceptional type $\mathsf{F}_4$.

Remarkably, all types of Dynkin diagrams are involved (though not every single one, the exceptions being A$_{2n}$, E$_7$ and E$_8$). A common feature of all these examples is that, for some point-line approach to these buildings, the point-line geometry of fixed building is a *geometric hyperplane* of the point-line geometry of the ambient building:

(1) The lines of a symplectic polar space form a linear system of the line Grassmannian of the ambient projective space;
(2) The parabolic quadric is a hyperplane section of the hyperbolic quadric;
(3) The points of the hexagon are the points of a parabolic quadric arising as a hyperplane section of the hyperbolic quadric of Witt index 4;
(4) The points of the corresponding metasymplectic space (where the points are the elements of type 4, with Bourbaki labeling [4] as in Fig. 2) are in a hyperplane section of the 16-dimensional variety corresponding to the building of type E$_6$ (the points are the elements of type 1) in 26-dimensional projective space, cf. [8].

In fact, through the symplectic polarity of a building of type E$_6$, we can witness all the above features: A symplectic polarity is induced in every fixed 5-space (feature (1)); in a residue of type D$_5$ that is mapped onto a non-incident element, the map defined by the polarity and the projection in the sense of Corollary 4.13 below, is an involution pointwise fixing a parabolic polar space of type B$_4$ (feature (2)); finally, the principal of triality (feature (3)) will be crucial in many of our arguments (see, for instance, Lemma 6.12). One could even say that buildings of type E$_6$ exist thanks to triality!

In this paper, we study the symplectic polarities of buildings of type E$_6$. Our main goal is to explicitly construct a building of type E$_6$ from a given split building of type F$_4$. This will be accomplished in Theorem 6.36. There are a number of reasons why this is a worthwhile thing to do. Firstly, it provides an explicit geometric link between these buildings much deeper than just knowing that the point set of the building of type F$_4$ is a geometric hyperplane of the point-line geometry naturally associated with a building of type E$_6$ (Theorem 1). Secondly, it provides a geometric proof of the fact that, up to conjugacy, every building of type E$_6$ admits a unique symplectic polarity (Theorem 2). Thirdly, we will use the gained geometric insight to show uniqueness of the inclusion in question of the split building of type F$_4$ in a building of type E$_6$ (Theorem 3). And last but not least, it provides a wealth of properties of the metasymplectic spaces related to split buildings of type F$_4$, which can be used in other geometric problems (for instance extending the results in [11] is a good candidate).

Other constructions of buildings from smaller ones comprise the construction of the ambient projective 3-space from the embedded symplectic generalized quadrangle (see [16]), the construction of the ambient oriflamme geometry of the polar space of type D$_4$ from the embedded triality generalized hexagon (see [24]), and the construction of the ambient metasymplectic parapolar space of type F$_4$ from the embedded Ree–Tits generalised octagon (see [15,25]).

## 2 Preliminaries and main results

For undefined but basic notions of the theory of buildings (such as opposition relation, chambers, (Dynkin) diagram, etc.), we refer the reader to the excellent textbook [1].

Let $\Delta$ be a building of type E$_6$ over the field $\mathbb{K}$. The latter means that each rank 2 residue is either a generalised digon or a projective plane over $\mathbb{K}$. We label the types according to the Bourbaki conventions [4] of labeling Dynkin diagrams and call elements of type 1, 2, 3, 4, 5, 6 *points,* 5-*spaces, lines, planes,* 4-*spaces and quads*, respectively (see Fig. 1).

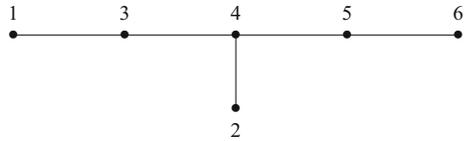

**Fig. 1** The Dynkin diagram of type $E_6$ with Bourbaki labeling

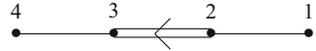

**Fig. 2** The Dynkin diagram of type $F_4$ with Bourbaki labeling

This way, we in fact identify $\Delta$ with its shadow space corresponding to the elements of type 1 (see e.g. [23]). The opposition relation on the types preserves the types 2 and 4 and switches type 1 with type 6 and type 3 with type 5. The corresponding point-line geometry is said to be *naturally associated with $\Delta$*.

Buildings of type $E_6$ are naturally associated with Chevalley groups of type $E_6$. It is well-known that each such group contains a maximal subgroup of type $F_4$, which is moreover pointwise fixed by an outer involutory automorphism. This involution induces a nontrivial involution (in fact, the opposition relation) on the diagram of $\Delta$ and can hence be seen as a polarity $\theta$ of $\Delta$. Geometrically, this maximal subgroup of type $F_4$ stabilizes a subbuilding $\Gamma$ of type $F_4$, consisting of some elements of types 1, 2, 3 and 4 of $\Delta$. This defines an embedding of $\Gamma$ in $\Delta$ with the following property: Every point, line and 5-space of $\Delta$ incident with a plane of $\Gamma$ belongs to $\Gamma$. We can choose types in $\Gamma$ such that the points, lines, planes and symplecta of $\Gamma$ are points, lines, planes and 5-spaces, respectively, of $\Delta$. The symplecta of $\Gamma$ are then symplectic polar spaces of rank 3, induced by $\theta$ on the $\theta$-fixed 5-spaces of $\Delta$. This motivates to call $\Gamma$ a *symplectic metasymplectic (parapolar) space* (see [7]; in general, a *metasymplectic space* is the point-line geometry obtained from *any* building of type $F_4$ by taking as points the objects either of type 1, or of type 4) and the polarity $\theta$ a *symplectic polarity* (see Proposition 4.14). It is unique up to conjugacy. Geometrically, this follows from the above property of the embedding of $\Gamma$ in $\Delta$, as we shall show in Theorem 6.37.

Now, the claims made in the previous paragraph are not easy to find in the literature. In the finite case, they follow from the classification of large almost simple maximal subgroups of groups of type $E_6$, see [12,13]. However, in the general case, hardly any literature exists about "non-Galois" automorphisms of buildings, or, more or less equivalently, of simple algebraic groups. In the present paper, we prove the above claims taking a geometric approach, and prove a new uniqueness result concerning the inclusion of split buildings of type $F_4$ into buildings of type $E_6$ (Theorem 3).

Hence we start with a split building $\Gamma$ of type $F_4$ over $\mathbb{K}$. This means that, again with Bourbaki labeling (see Fig. 2), the residues of type 1 are the buildings corresponding to the polar spaces of rank 3 defined by a symplectic polarity in $\mathsf{PG}(5, \mathbb{K})$. We will define additional elements using special substructures of $\Gamma$ to obtain $\Delta$, recover the polarity $\theta$ (Theorem 6.36) and prove its uniqueness up to conjugacy (Theorem 6.37). More exactly, we will prove the following theorems.

**Theorem 1** *Let $\Gamma$ be a symplectic metasymplectic parapolar space. Then there exist an explicitly defined geometry $(\mathfrak{E}, *)$, which is the geometry naturally associated with a building of type $E_6$, and a symplectic polarity $\theta$ of $(\mathfrak{E}, *)$ with fixed point structure $\Gamma$. More exactly, the sets of absolute points and absolute lines of $\theta$ are precisely the sets of points and lines of $\Gamma$, and the fixed planes and fixed 5-spaces of $\theta$ are the planes and symplecta, respectively, of $\Gamma$.*

**Theorem 2** *A building $\Delta$ of type* E$_6$ *admits, up to conjugacy, a unique symplectic polarity.*

In Sect. 7, we let $\Delta$ be a building of type E$_6$ over the field $\mathbb{K}$ and $\theta$ a symplectic polarity. An element of $\Delta$ of type 1,3,5 or 6 is called *absolute* if it is incident with its image. The absolute points, absolute lines, fixed planes and fixed 5-spaces with inherited incidence relation from E$_6$ are the points, lines, planes and symplecta of a metasymplectic parapolar space all of whose symplecta are symplectic polar spaces of rank 3. In other words, these four types of objects form a split building $\Gamma$ of type F$_4$. Now we view $\Delta$ and $\Gamma$ as independent point-line geometries (hence neglecting all objects other than the points and the lines and their mutual incidence) and say that $\Gamma$ *is point-line-embedded in* $\Delta$ if the point set of $\Gamma$ is a subset of the point set of $\Delta$, likewise for the line sets, and if incidence in $\Gamma$ is inherited from incidence in $\Delta$. The point-line-embedding is called *full* if all points of $\Delta$ on a line of $\Gamma$ are also points of $\Gamma$. In Sect. 7, we prove the following theorem.

**Theorem 3** *Let $\Gamma$ be a symplectic metasymplectic parapolar space. Let $\Delta$ be the natural point-line geometry associated with a building of type* E$_6$. *If $\Gamma$ is fully point-line-embedded in $\Delta$, then $\Gamma$ and $\Delta$ are defined over the same field and $\Gamma$ arises from a symplectic polarity of $\Delta$.*

Note that the condition of the embedding being full cannot be dispensed with since one can consider a symplectic metasymplectic parapolar space over a subfield of $\mathbb{K}$. Also, there exist metasymplectic parapolar spaces fully embedded in $\Delta$ which are not symplectic, i.e., which arise from a non-split building $\Gamma'$ of type F$_4$. This happens when $\mathbb{K}$ admits a Galois involution and $\Gamma'$ is the building associated with the corresponding twisted Chevalley group.

**Notation and Terminology**

- We will use the following convention to refer to certain buildings. With "type" of a building, we will always mean the relative and absolute Dynkin type (hence all buildings we consider are split) in the sense of algebraic groups, except when we explicitly mention that we only mean relative type. Hence,

  - A building of type B$_n$ is the simplical complex defined by a parabolic quadric (which is a nondegenerate quadric of Witt index $n$ in a projective space of dimension $2n$),
  - A building of type C$_n$ is the simplicial complex defined by a symplectic polarity in a projective space of dimension $2n - 1$, and
  - A building of type D$_n$ is the oriflamme complex (see Sect. 7.12 of [22]) of a $(2n-2)$-dimensional hyperbolic quadric, i.e., a nondegenerate quadric of maximal Witt index $n$ in a projective space of dimension $2n - 1$.

- We will recognise buildings via their standard point-line representations, which we will identify with the buildings themselves. This comprises:

  - Parabolic quadrics or polar spaces of type B$_n$ to mean that the given point-line geometry conforms to a building of type B$_n$ where points and lines are the elements of type 1 and 2, respectively;
  - Symplectic polar spaces, or polar spaces of type C$_n$, which conform to buildings of type C$_n$;
  - Quadrics of type D$_n$ are hyperbolic quadrics whose oriflamme complexes are buildings of type D$_n$;
  - Finally, a metasymplectic space is a point-line geometry (a so-called *parapolar space*, see [18]) associated with a (thick) building of type F$_4$.

- A *subspace* of a point-line geometry is a subset of points with the property that every line intersecting the set in at least two points is entirely contained in it.
  - A subspace is called *singular* if every pair of points is collinear.
  - A subspace is called a *(geometric) hyperplane* if every line intersects it in at least one point.
  - A geometric hyperplane is called *singular* if it consists of the set of all points not at maximal distance (in the point graph) from a certain point $x$ (for polar spaces, this just means the set of points collinear with or equal to $x$).

Hence there is some ambiguity when talking about a hyperplane to be singular: it can be singular as a subspace, or it can be singular as a hyperplane. We make the convention that when we write "singular (geometric) hyperplane", we always mean singular as a hyperplane, and when we write "singular subspace", then we mean that every pair of points is collinear, even if the subspace happens to be a geometric hyperplane.

*Remark 2.1* Let $G$ be the linear algebraic group of type $\mathsf{E}_6$ over the field $\mathbb{K}$. Let $G(\theta)$ be the centralizer of the standard symplectic polarity $\theta$ (in the terminology of [17], a *graph automorphism*). Then $G(\theta)$ is a linear algebraic group of split type $\mathsf{F}_4$ and a maximal closed connected subgroup of $G$. Theorem 15.1 of [17] asserts that, if $\mathbb{K}$ is algebraically closed, all maximal closed connected subgroups of $G$ of exceptional type are conjugate. This could be the base for an algebraic proof of Theorem 2.

## 3 Outline of the paper

We now explain the structure of the paper, in particular of the proofs of Theorems 1, 2 and 3.

The proofs of Theorems 1 and 2 comprise the first part of the paper, the proof of Theorem 3 will be referred to as the second part. We note that,

(I) In the first part, we start with an arbitrary symplectic metasymplectic parapolar space $\Gamma$ and
   (a) Construct a canonical building $\Delta$ of type $\mathsf{E}_6$ in which $\Gamma$ is embedded, and
   (b) Show that $\Gamma$ is the fixed point structure of a symplectic polarity of $\Delta$;
(II) In the second part, we start with a pair $(\Delta, \Gamma)$ consisting of a building $\Delta$ of type $\mathsf{E}_6$ and a split building $\Gamma$ of type $\mathsf{F}_4$, whose associated symplectic metasymplectic parapolar space is point-line-embedded in $\Delta$ and we again show that $\Gamma$ is the fixed point structure of a symplectic polarity of $\Delta$.

For Parts (I)(b) and (II), we need some background on buildings of type $\mathsf{E}_6$, in particular how to recognise symplectic polarities. This is provided in Sect. 4. Hence, strictly speaking, we do not need that section for Part (I)(a). However, we prove slightly more than strictly needed in Sect. 4 for the benefit of the reader. Indeed, we also point out the relations between some subspaces of a building of type $\mathsf{E}_6$ and a symplectic polarity, or its fixed point set. These relations can be kept in the back of the reader's mind as a distant guide when going through Sect. 6. In particular, it hints at how certain subspaces of $\Delta$ will be recovered from $\Gamma$.

Of course, for Part (I)(a) we need (basic) properties of symplectic metasymplectic parapolar spaces, and these are gathered in Sect. 5.1. The rest of Sect. 5 is devoted to preparing the (re)construction of the elements of the building $\Delta$ of type $\mathsf{E}_6$. The bulk of that is to get the quads of $\Delta$. We will recognise these by their intersections with $\Gamma$. One type of intersection is a substructure of $\Gamma$ which we shall call "extended equator geometry", and we define it in Sect. 5.2. We study some properties which shall enable us later to identify these with the intended intersections.

Now, we also have to construct a symplectic polarity once we defined $\Delta$. Such a polarity maps points to quads and vice versa. So, if we consider a point $p$ of $\Delta$, then we can look at the intersection $\widehat{E}$ of its image with $\Gamma$, or we can also look at the set $\widehat{T}$ of points of $\Gamma$ collinear with $p$. It turns out that there is a very neat relation between $\widehat{E}$ and $\widehat{T}$, and the latter can be constructed only knowing $\widehat{E}$. We show how to do this, and we call $\widehat{T}$ the "tropical circle geometry", hinting at the fact that it lies "in the neighbourhood of" and "surrounds" the (extended) equator geometry. Together, $\widehat{E}$ and $\widehat{T}$ and the points on a line joining them, define a geometric hyperplane of $\Gamma$, and we also prove this (see Sect. 5.4). At that point we can already define the point-line geometry related to $\Delta$ (Sect. 6.1). The rest of Sect. 6 is devoted to the construction of all elements of $\Delta$ and the symplectic polarity: the quads (Sect. 6.3), the 4-spaces (Sect. 6.4), the 5-spaces and the planes (Sect. 6.5), and the symplectic polarity (Sect. 6.6). In Sect. 6.7 we put all pieces together to prove Theorems 1 and 2.

Finally, in Sect. 7, we prove Theorem 3. We are given a point-line-embedding of the symplectic metasymplectic parapolar space geometry related to a split building $\Gamma$ of type F$_4$ in the natural point-line geometry related to a building $\Delta$ of type E$_6$. The first problem to solve is how the symplecta of $\Gamma$ are embedded in $\Delta$. It turns out that there are just two possibilities: either as a fully embedded symplectic polar space in a 5-space of $\Delta$, or as a subquadric of a quad (only possible in characteristic 2). We rule out that second possibility in Sect. 7.1. In Sect. 7.2, we analyse the relation between $\Gamma$ and a quad of $\Delta$ via the intersection of the latter with the former, and we use this in Sect. 7.3 to finally construct the associated symplectic polarity. Much of the argumentation should look somehow predictable by the gained insight into the structural properties of the standard inclusion of $\Gamma$ in $\Delta$ thanks to our analysis in the previous sections.

## 4 Some basic properties of buildings of type E$_6$

### 4.1 Generalities, types, apartments

We gather some facts about buildings $\Delta$ of type E$_6$. Throughout, we number the diagram as in Figs. 1 and 2, and choose to name the elements of type 1 *points*. We identify all other elements with the set of points incident with them. The elements of type 3 will be called *lines*, those of type 4 *planes*, those of type 5 will be called 4-*spaces*, those of type 6 *quads* and the elements of type 2 will be called 5-*spaces*. Note that 4-spaces are projective spaces of dimension 4 over $\mathbb{K}$, likewise 5-spaces are projective 5-spaces over $\mathbb{K}$, and quads are subspaces isomorphic to hyperbolic quadrics of type D$_5$ defined in some projective 9-space over $\mathbb{K}$, i.e., quadrics of maximal Witt-index 5 in such a space. A 4′-*space* is a hyperplane of a 5-space, but it does not conform to a (single) type in $\Delta$ (it can be considered as a flag of type {2, 6}). Also, a 3-*space* is some 3-space in a 5-space, or, equivalently, in a 4-space (it conforms to a flag of type {2, 5, 6}). So we obtain a point-line geometry, which we will call the *natural point-line geometry associated with* $\Delta$. The elements of types 2, 4, 5, 6 conform to subspaces of this point-line geometry, and their names mentioned above are chosen such that they reveal the structure of the subspace in question (e.g., planes are really projective planes, etc.).

If we choose the elements of type 6 to be the points, then we obtain an isomorphic geometry where the elements of type 5 are the lines, those of type 4 are the planes, those of type 3 are the 4-spaces, those of type 1 are the quads and the elements of type 2 are the 5-spaces. This is the *principle of duality* in buildings of type $\mathsf{E}_6$. It follows from the uniqueness of the building of type $\mathsf{E}_6$ over the field $\mathbb{K}$, see [22].

Given the natural point-line geometry associated with $\Delta$, we can go back to the incidence geometry or, equivalently, the numbered simplicial complex defined by $\Delta$ as follows. The elements of the geometry or the vertices of the simplicial complex are the points, lines, planes, 4-space, 5-spaces and quads, and the incidence or adjacency is given by the following rules. An object is *incident* with another object if, and only if, one of them is contained in the other, except when one object is a 5-space $V$, and the other is either a 4-space $W$ or a quad $\Sigma$. Then $V$ is incident with $W$ if, and only if, $V \cap W$ is a 3-dimensional subspace of both $V$ and $W$; $V$ is incident with $\Sigma$ if, and only if, $V \cap \Sigma$ is a 4-dimensional (singular) subspace of both $V$ and $\Sigma$, in which case it is a $4'$-space. The two families of maximal singular subspaces of $\Sigma$, characterised by the property that subspaces from the same family meet each other in even-dimensional subspaces, and members of different families meet each other in odd-dimensional subspaces, are the families of 4-spaces and of $4'$-spaces contained in $\Sigma$.

### 4.2 Some geometric properties of the point-line geometry $\Delta$

In this subsection, we assume that $\Delta$ is the natural point-line geometry associated with a (thick) building of type $\mathsf{E}_6$ defined over some field $\mathbb{K}$. Collinearity, from now on, refers to distinct points incident with (or, with our convention, contained in) at least one line of $\Delta$. *Opposite* elements are elements which are opposite in some apartment. A *flag* is a set of pairwise incident elements of different types, so a *chamber* in $\Delta$ is a flag consisting of 6 elements of $\Delta$.

Everything below is well-known, and we give precise references for most facts. Many things are contained in [20], but we also include references to [9], as the latter is easily accessible and provides an excellent source of information on buildings of type $\mathsf{E}_6$. Let us also remark that some of the properties are stated, without proof, in [14], where they are seen as results of "reading" the diagram. We were unable to find Facts 4.7, 4.10 and 4.11 explicitly in the literature, but these (and also the others) can be verified by the reader himself by including two appropriate flags (mostly just two elements) in an apartment. The assertion then becomes an assertion in a thin building of type $\mathsf{E}_6$. Such a thin building $\mathscr{A}$ is provided by the following easy construction (see Paragraph 10.3.4 in [3]): the 27 points of $\mathscr{A}$ are the 27 points of the generalized quadrangle $Q$ of order (2, 4) (arising from a nondegenerate bilinear form of Witt index 2 in a 5-dimensional projective space over the field of 2 elements). The lines of $\mathscr{A}$ are the non-collinear pairs of points of $Q$. The planes of $\mathscr{A}$ are the *triads* of $Q$ (i.e., the triples of non-collinear points). The 4-spaces are the intersections $p^\perp \cap q^\perp$, where $p$ and $q$ are two non-collinear points and $x^\perp$ denotes the set of elements collinear or equal to the point $x$ in $Q$. The 5-spaces through a point $p$ are obtained by taking some point $q$ not collinear with $p$ in $Q$, and then the points in $q^\perp \setminus (p^\perp \cup \{q\})$ together with $p$ form a 5-space. A quad simply is $p^\perp \setminus \{p\}$ for some point $p$ of $Q$. Opposition is also easily defined in $\mathscr{A}$. Indeed, a point $p$ is opposite the quad $p^\perp \setminus \{p\}$; a line $\{p, q\}$ is opposite the 4-space $p^\perp \cap q^\perp$; the plane $\{x, y, z\}$ is opposite the plane $x^\perp \cap y^\perp \cap z^\perp$ and the 5-space $\{p\} \cup (q^\perp \setminus (p^\perp \cup \{q\}))$ is opposite the 5-space $\{q\} \cup (p^\perp \setminus (q^\perp \cup \{p\}))$.

**Fact 4.1** (Lemma 18.7.1 of [9], Statement 3.7 of [20]) *Any pair of collinear points of $\Delta$ is contained in a unique line. Any pair of non-collinear points of $\Delta$ is contained in a unique quad.*

By duality, we have the following.

**Fact 4.2** *Two distinct quads of $\Delta$ either intersect in a unique 4-space, or in a unique point.*

**Fact 4.3** (Proposition 18.7.2(vii) of [9], Statements 3.5.4 and 3.9 of [20]) *Given a point $x$ and a quad $\Sigma$, then either $x \in \Sigma$, or $x$ is opposite $\Sigma$, which is equivalent to "no point of $\Sigma$ is collinear with $x$", or there is a unique 5-space $V$ incident with both $x$ and $\Sigma$. In the latter case, the intersection of $V$ with $\Sigma$ is precisely the set of points of $\Sigma$ collinear with $x$.*

In the last case, namely when there is a unique 5-space incident with both a point $x$ and a quad $\Sigma \not\ni x$, we say that $\Sigma$ *neighbors* $x$ (and $x$ *neighbors* $\Sigma$). This notion is standard in the theory of Hjelmslev planes and is inspired by the fact that $\Delta$ can be described as a Hjelmslev–Moufang plane over split octonions, see [19].

The previous fact has the following consequence.

**Corollary 4.4** *A quad $\Sigma$ is convex.*

*Proof* Let $x, y$ be two points of $\Sigma$. If they are collinear, then clearly, the line joining them is contained in the subspace $\Sigma$. If they are not collinear then let $z$ be a point collinear with $x$ and $y$ and suppose $z \notin \Sigma$. By Fact 4.3, $z$, $x$ and $y$ are contained in a unique 5-space, which is a contradiction. Hence $z \in \Sigma$. □

Fact 4.5 below is an immediate consequence of the fact that one can put any 4-space and any point in a common apartment, in which there are exactly two quads incident with the 4-space, and that in the apartment, each point is opposite only one quad.

**Fact 4.5** *At least one quad through a given 4-space is not opposite a given point.*

**Fact 4.6** (Proposition 18.7.2(v) of [9], Statement 3.5.3 of [20]) *Two 5-spaces are either disjoint, intersect in a point, or intersect in a plane. The latter case is equivalent to the 5-spaces being incident with a common plane (namely, their intersection). In particular, every 3-space is contained in a unique 5-space.*

**Fact 4.7** *Two disjoint 5-spaces are either opposite or there exists a 5-space intersecting them in disjoint planes.*

**Fact 4.8** (Proposition 18.7.2(v) of [9], Statement 3.2 of [20]) *Every 3-space is contained in a unique 4-space.*

**Fact 4.9** (Statement 3.2 of [20]) *Every $4'$-space is contained in a unique quad and in a unique 5-space.*

**Fact 4.10** *Given a point $x$ and a 5-space $V$, either $x$ and $V$ are incident, or $x$ is collinear with exactly one point of $V$, or $x$ is collinear with all points of a unique 3-space of $V$. In the latter case, the space spanned by $x$ and $x^\perp \cap V$ (i.e., the union of all lines through $x$ meeting $V$) is a 4-space.*

**Fact 4.11** *A point, line or plane is opposite a quad, 4-space, or plane, respectively, if, and only if, the collinearity relation between the two elements is empty. A 5-space is opposite another 5-space if, and only if, each point of the first is collinear with a unique point of the second 5-space if, and only if, each point of either of them is collinear with a unique point of the other.*

Now let $F$ and $F'$ be *opposite flags* in $\Delta$, i.e., each element of $F$ is opposite a unique element of $F'$ and vice versa. For every chamber $C$ containing $F$, there is a unique chamber $C'$ containing $F'$ at minimal distance from $C$ (where the distance of chambers is measured in the chamber graph, i.e., the graph with vertices the chambers, and two chambers are adjacent if they share 5 elements). We denote the map $C \mapsto C'$ by $\rho_{F,F'}$. The *residue of $F$* consists of all chambers containing $F$ and carries the structure of a spherical building. It is well-known that $\rho_{F,F'}$ can be naturally extended to *all* flags containing $F$, see Theorem 3.28 in [22]. In the following proposition and corollary, $\Delta$ denotes any spherical building. We will apply these statements to buildings of types $\mathsf{D}_4$, $\mathsf{F}_4$ and $\mathsf{E}_6$.

**Proposition 4.12** (Theorem 3.28 and Proposition 3.29 of [22]) *Let $F$ and $F'$ be opposite flags of $\Delta$. Then $\rho_{F,F'}$ is an isomorphism from the residue of $F$ to the residue of $F'$ and the type of the image of an element of type $i$ is the opposite in the residue of $F'$ of the opposite type of $i$ in $\Delta$. Also, chambers $C \supseteq F$ and $C' \supseteq F'$ are opposite in $\Delta$ if, and only if, $C'$ and $\rho_{F,F'}(C)$ are opposite in the residue of $F'$.*

There is a useful corollary.

**Corollary 4.13** *Let $\varphi$ be an automorphism of $\Delta$. Let $F$ and $F^\varphi$ be opposite flags of $\Delta$, and let $\sigma_{F,\varphi}$ be the automorphism of the residue of $F$ mapping a chamber $C \supseteq F$ onto $\rho_{F^\varphi,F}(C^\varphi)$. If $\varphi$ induces the natural opposition relation on the types of $\Delta$, then so does $\sigma_{F,\varphi}$ for the residue of $F$.* □

We end this section with three results proved in [27].

**Proposition 4.14** *Let $\theta$ be a duality of a building of type $\mathsf{E}_6$. The following are equivalent.*

- *$\theta$ is a symplectic polarity;*
- *$\theta$ maps no point to a neighboring quad and at least one point is absolute;*
- *$\theta$ maps no chamber to an opposite chamber;*
- *$\theta$ maps no line (or plane or 4-space or 5-space, respectively) to an opposite one.*

**Proposition 4.15** *Let $\Delta$ be a building of type $\mathsf{E}_6$ and let $\theta$ be a symplectic polarity of $\Delta$. Let $\Gamma$ be the building of type $\mathsf{F}_4$ consisting of the absolute points, absolute lines, fixed planes and fixed 5-spaces for $\theta$. Then a line $L$ of $\Delta$ containing at least two points of $\Gamma$ is entirely contained in $\Gamma$ and either $L$ is an absolute line, or $L$ is a hyperbolic line in some fixed 5-space $V$ (hyperbolic with respect to the symplectic polarity in $V$ induced by $\theta$ by relating a point $x \in V$ to the intersection $V \cap x^\theta$) in the sense of Sect. 5.1.*

**Proposition 4.16** *Let $\Delta$ be a building of type $\mathsf{E}_6$ and let $\theta$ be a symplectic polarity of $\Delta$. Let $V$ be a 5-space of $\Delta$. Then every point of $V$ is absolute for $\theta$ if, and only if, $V$ is fixed by $\theta$.*

### 4.3 Properties of symplectic polarities of buildings of type E$_6$

We now prove some additional properties of symplectic polarities. The goal is to gain some insight to enable us to recognise and define the elements of the building of type E$_6$ out of the elements of a building of type F$_4$ when knowing that the latter is the fixed point set in the former of a symplectic polarity.

Throughout, let $\Delta$ be a building of type E$_6$, with labeling and names of objects as above. In particular, we identify $\Delta$ with its natural point-line geometry, which is a partial linear space consisting of a set of points, endowed with certain subsets called lines, planes, 3-spaces, 4-spaces, 4′-spaces, 5-spaces and quads. Let $\theta$ be a symplectic polarity of $\Delta$. We denote by $H$ the set of its absolute points.

The first lemma holds for arbitrary polarities in E$_6$.

**Lemma 4.17** *Every polarity $\rho$ of $\Delta$ admits at least one non-absolute point.*

*Proof* Suppose for a contradiction that all points of $\Delta$ are absolute. Let $x$ be a point of $\Delta$. By our assumption, $x \in x^\rho$. Let $y$ be a point in $x^\rho$ not collinear with $x$. As $y \in x^\rho$, the quad $y^\rho$ contains $x$. By our assumption it also contains $y$ and hence, by Fact 4.1, $x^\rho = y^\rho$, a contradiction. □

The following is well-known, but we provide a proof for completeness. Note that a geometric hyperplane of a point-line geometry is called *proper* if it does not coincide with the full point set.

**Lemma 4.18** *The set $H$ is a proper geometric hyperplane of $\Delta$.*

*Proof* If a line $L$ has at least two points in common with $H$, then all points of $L$ belong to $H$ by Proposition 4.15. So we are left to show that no line is disjoint from $H$.

Let $L$ be any line of $\Delta$. By Proposition 4.14, $L^\theta$ is not opposite $L$. This means, in view of Fact 4.11, that some point $x \in L$ is collinear with some point $y \in L^\theta$. Since $x \in L$, we have $L^\theta \subseteq x^\theta$, and so $x$ is collinear with a point of its image. Proposition 4.14 together with Fact 4.3 imply that $x$ is absolute.

The properness of $H$ follows from Lemma 4.17. □

The next lemma tells us which lines are absolute and is a more detailed version of Proposition 4.15.

**Lemma 4.19** *Let $L$ be a line of $\Delta$. Then, $L$ is absolute if, and only if, it contains a point $x$ such that $L \subseteq x^\theta$. Moreover, every absolute line is contained in $H$.*

*Proof* Suppose first that $L$ contains a point $x$ with $L \subseteq x^\theta$. Consider a second point $y \in L \setminus \{x\}$. Then, since $y \in x^\theta$, we have $x \in y^\theta$ and Fact 4.3 and Proposition 4.14 imply that $y \in y^\theta$. It follows that $L \subseteq x^\theta \cap y^\theta = L^\theta$. Hence $L$ is absolute.

If $L$ is absolute and $x \in L$, then $x \in L \subseteq L^\theta \subseteq x^\theta$ and hence $L \subseteq H$. □

We now take a look at the relation of $H$ with a quad. Henceforth, the symbol "⊥" refers to collinearity in $\Delta$ and $x^\perp$ is the set of points equal to or collinear with the point $x$.

**Lemma 4.20** *Let $\Sigma$ be a quad of $\Delta$. Then, $\Sigma$ is absolute if, and only if, $\Sigma \cap H$ is a singular hyperplane of $\Sigma$, i.e., there exists a point $x \in \Sigma$ such that $\Sigma \cap H = x^\perp \cap \Sigma$. In this case, $x = \Sigma^\theta$. If $\Sigma$ is not absolute, then $\Sigma \cap H$ is a parabolic quadric of type B$_4$ (a geometric hyperplane of $\Sigma$ which is not singular).*

*Proof* If $\Sigma$ is absolute, then $\Sigma^\theta = x$ is an absolute point in $\Sigma \cap H$, and every line in $\Sigma$ through $x$ is absolute by Lemma 4.19. The same lemma also implies that every such line is contained in $H$, and hence $x^\perp \cap \Sigma \subseteq H \cap \Sigma$. If a point $z \in \Sigma$ not collinear with $x$ were absolute, then $z^\theta$ would contain $z$ (as $z$ is absolute) and $x$ (as $z$ is contained in $x^\theta$), and so $z^\theta = \Sigma$ by Fact 4.1, a contradiction. Consequently $\Sigma \cap H = x^\perp \cap \Sigma$.

Now suppose that $\Sigma \cap H = x^\perp \cap \Sigma$ for some $x \in \Sigma$. It suffices to show that $\Sigma$ is absolute. Assume, by way of contradiction, that it is not. By Fact 4.3 and Proposition 4.14, $\Sigma^\theta$ is not collinear with any point of $\Sigma$. So, by Fact 4.2, $\Sigma \cap p^\theta$ is a point, for every $p \in \Sigma$. Further, by Corollary 4.13, the map sending a point $p \in \Sigma$ to the intersection $p^\theta \cap \Sigma$ is an involution $\sigma$ of $\Sigma$ with $x^\perp \cap \Sigma$ as the set of its fixed points. Let $z$ be a point in $\Sigma$ not collinear with $x$. Since $x^\perp \cap z^\perp$ is a polar space of type $\mathsf{D}_4$, we can find two disjoint 3-spaces $U, U'$ in $x^\perp \cap z^\perp$, and each of them is contained in exactly two 4-spaces of $\Sigma$ (here we do not distinguish between the 4-spaces and 4'-spaces of $\Delta$), namely $\langle z, U \rangle$ and $\langle x, U \rangle$, and $\langle z, U' \rangle$ and $\langle x, U' \rangle$, respectively. But the 4-spaces $\langle x, U \rangle$ and $\langle x, U' \rangle$ and the 3-spaces $U$ and $U'$ are fixed by $\sigma$, hence the 4-spaces $\langle z, U \rangle$ and $\langle z, U' \rangle$ are also fixed, as is their intersection $\langle z, U \rangle \cap \langle z, U' \rangle = \{z\}$. So, $z \in x^\perp \cap \Sigma$, a contradiction. Hence, $\Sigma$ is absolute.

Hence, the first assertion is proved. The second assertion follows from the fact that $\Sigma \cap H$ is a proper geometric hyperplane of $\Sigma$ (the properness follows from the proof of Lemma 4.17) and $\Sigma$ has only two types of those, singular ones (and then the first assertion applies) and polar subspaces of type $\mathsf{B}_4$, as can easily be checked. □

In particular, we have the following statement.

**Corollary 4.21** *If $x \in H$, then $x^\theta \cap H = x^\theta \cap x^\perp$.* □

Next, we study the relation of $H$ with 4-spaces.

**Lemma 4.22** *Let $V$ be a 4-space of $\Delta$. Then, $V$ is absolute if, and only if, $V \subseteq H$.*

*Proof* First suppose that $V$ is contained in $H$, and put $L = V^\theta$. By Lemma 4.18, there exists a point $x \in L \cap H$. It follows that $x^\theta$ contains $V$, and hence $V \subseteq H \cap x^\theta$. By Corollary 4.21, $V \subseteq x^\perp \cap x^\theta$. All 4-spaces of $x^\theta$ in $x^\perp$ contain $x$. So, $x \in V$. Consequently, $L = V^\theta$ is contained in $x^\theta$ and contains $x$. By Lemma 4.19, $L$ is absolute, and hence, so is $V$.

Now suppose that $V$ is absolute. Then $L = V^\theta \subset V$. For any $x \in V$ we have that $L \subset x^\theta$, and as $x$ is collinear with $L$, Proposition 4.14 and Fact 4.3 imply $x \in x^\theta$. Hence, $V \subseteq H$. □

The relation of $H$ with a 5-space was the content of Proposition 4.16. Although we will not need a similar result for planes, we mention it for completeness.

**Lemma 4.23** *Let $\pi$ be a plane of $\Delta$. Then, $\pi$ is fixed by $\theta$ if, and only if, it contains two distinct points $x, y$ such that $\pi \subseteq x^\theta \cap y^\theta$. Moreover, every fixed plane is contained in $H$.*

*Proof* Suppose $x, y$ are distinct points of $\pi$ with $\pi \subseteq x^\theta \cap y^\theta$. Let $z$ be a point of $\pi \setminus xy$. As $z \in x^\theta \cap y^\theta$, we have that $x$ and $y$ belong to $z^\theta$, and therefore also $z \in z^\theta$ by Fact 4.3 and Proposition 4.14. This implies that $\pi$, which is generated by $x, y$ and $z$, is contained in $z^\theta$. It follows that $\pi = x^\theta \cap y^\theta \cap z^\theta = \pi^\theta$ and hence $\pi$ is fixed.

If $\pi$ is fixed and $x \in \pi$, then $x \in \pi = \pi^\theta \subseteq x^\theta$ and hence $\pi \subseteq H$. □

**Lemma 4.24** *Let $x$ be any point of $\Delta \setminus H$ and let $p$ be a point in $x^\perp \cap H$. Then, the set of points of $x^\theta \cap H$ collinear with $p$ forms a 3-space $U$ which is entirely contained in $p^\theta$. In particular, every line through $p$ intersecting $U$ is absolute.*

*Proof* Since $p$ is collinear with $x$, and the line $px$ is not absolute by Lemma 4.19, the quads $p^\theta$ and $x^\theta$ intersect in a non-absolute 4-space $V = (xp)^\theta$. Put $U = V \cap H$. Then, by Lemmas 4.18 and 4.22, $U$ is 3-dimensional. By Lemma 4.19, every line in $p^\theta$ through $p$ is absolute, and hence is contained in $H$. Since $p$ is collinear with every point of a 3-space of $V$ inside the polar space $p^\theta$, we see that $p$ is collinear with $U$ and with no point of $V \setminus U$. ⊓⊔

*Remark 4.25* The previous results imply that every absolute singular subspace of $\Delta$ (including the fixed ones) related to some node in the diagram (namely, nodes 2,3,4 and 5) is contained in $H$. This is also true for the $4'$-spaces, which are in fact flags of type $\{2, 6\}$, as the unique 5-space containing an absolute $4'$-space must be absolute (by definition of an absolute flag, which is a flag $F$ such that $F \cup F^\theta$ is also a flag) and hence belongs to $H$.

We can now start preparing for the construction of a building of type E$_6$ starting from a split building of type F$_4$.

## 5 Some special substructures of split buildings of type F$_4$

### 5.1 Properties of symplectic metasymplectic parapolar spaces

In this subsection we list basic properties of symplectic metasymplectic parapolar spaces. Most of them are just all possibilities of mutual position between two elements.

Let $\Gamma$ be a split building of type F$_4$ over $\mathbb{K}$. We view $\Gamma$ as a symplectic metasymplectic parapolar space. This means that we have a set of points (and this set is precisely the set of elements of type 4 of the building $\Gamma$, see Fig. 2), a set of lines (elements of type 3), a set of planes (elements of type 2) and a set of *symplecta* (elements of type 1) and these are such that each line, each plane and each symplecton is a proper convex subset of the set of points. In particular, $\Gamma$ is a partial linear space. The planes are projective planes when endowed with the lines of $\Gamma$ they contain; the lines and planes contained in a symplecton render it a symplectic polar space of rank 3. The opposition relation in $\Gamma$ ([22], Chapter 7) acts on the types as the identity. The basic properties of $\Gamma$ are stated below, as facts. As noted on page 80 of [26], these can be proved using the diagram of type F$_4$; they also follow from [7]. Facts 5.3, 5.4 and 5.5 are valid in any metasymplectic space, in particular in the dual of $\Gamma$ (the set of points of the dual is the set of symplecta of the original, lines of the dual correspond to planes of the original).

**Fact 5.1** *The symplecta, planes and lines of $\Gamma$ through a given point $p$, endowed with the natural incidence relation, form a polar space $R(p)$ of type B$_3$ over $\mathbb{K}$, where the points of that polar space are the symplecta through $p$, the lines are the planes through $p$, and the planes are the lines through $p$.*

In particular, it follows that the isomorphism class of the geometry $R(p)$ does not depend on $p$. It is usually called the *point residue geometry* of $\Gamma$. Another consequence is the following.

**Corollary 5.2** *Every singular subspace of $\Gamma$ is contained in some symplecton, and hence is either a point, a line or a projective plane.* ⊓⊔

**Fact 5.3** *Let $x$ and $y$ be two points of $\Gamma$. Then, precisely one of the following situations occurs.*

(0) $x = y$;
(1) There is a unique line incident with both x and y. In this case, we call x and y collinear. We denote the unique line joining them by $xy$ and write $x \perp y$;
(2) There is a unique symplecton incident with both x and y. In this case, there is no line incident with both x and y, and we call x and y symplectic, or say that $\{x, y\}$ is a symplectic pair, or say that x is symplectic to y. We denote the unique symplecton by $x \Diamond y$ and write $x \perp\!\!\!\perp y$;
(3) There is a unique point z collinear with both x and y. In this case, we call x and y special, or say that $\{x, y\}$ is a special pair, or say that x is special to y. We denote z by $x \bowtie y$. For every pair $\{x, z\}$ of collinear points, there is a point y such that $x \bowtie y = z$;
(4) There is no point collinear with both x and y. In this case, we call x and y opposite. For every point x there is at least one point y opposite x.

Moreover, each of these possibilities occurs.

**Fact 5.4** *The intersection of two symplecta is either empty, or a point, or a plane and each of these occurs. Also, the graph with vertices the symplecta, where two symplecta are adjacent if they meet in a plane, is connected.*

**Fact 5.5** *Let x be a point and S a symplecton of $\Gamma$. Then precisely one of the following situations occurs.*

(0) $x \in S$;
(1) The set of points of S collinear with x is a line L. Every point y of $S \setminus L$ which is collinear with each point of L is symplectic to x and $x \Diamond y$ contains L. Every other point z of S (i.e., every point z of S collinear with a unique point z' of L) is special to x and $x \bowtie z = z' \in L$. We say that x and S are close;
(2) There is a unique point u of S symplectic to x and $S \cap (x \Diamond u) = \{u\}$. All points v of S collinear with u are special to x and $x \bowtie v \notin S$. All points of S not collinear with u are opposite x. We say that x and S are far.

Moreover, each of these possibilities occurs.

The previous facts are fundamental and will sometimes be used without referring to them.

For a point x, we denote by $x^\perp$ and $x^{\perp\!\!\!\perp}$ the sets of points collinear or equal to x and symplectic or equal to x, respectively; likewise, for a set A of points we denote by $A^\perp$ the set $\bigcap_{a \in A} a^\perp$. We say that a point x is *collinear* with a set A if $A \subseteq x^\perp$. For a symplectic pair $\{x, y\}$ of points of $\Gamma$, the set of points $h(x, y) = (x^\perp \cap y^\perp)^\perp \subseteq x \Diamond y$ is called a *hyperbolic line*. It has the property that for each point z in $x \Diamond y$, the set $z^\perp \cap h(x, y)$ is either a singleton, or the whole set $h(x, y)$. Note that this implies that there are no planes in $x \Diamond y$ that consist of hyperbolic lines only, see Definition 5.20. Also, since $\Gamma$ is defined over $\mathbb{K}$, it is *thick*, meaning that all lines, and hence also all hyperbolic lines, have at least three points. Finally, for arbitrary distinct $x', y' \in h(x, y)$, we have $x \Diamond y = x' \Diamond y'$. Putting $h = h(x, y)$, we set $S(h) = x \Diamond y$.

In the next two lemmas, we establish the possible mutual positions of a point and a hyperbolic line, and of a point and a line.

**Lemma 5.6** *Let h be a hyperbolic line in $\Gamma$ and x a point. Then, exactly one of the following holds.*

(i) $x \in h$;
(ii) x is collinear with every point of h;

  (iii) $x$ is collinear with exactly one point of $h$ and symplectic to the other points of $h$;
  (iv) $x$ is collinear with exactly one point of $h$ and special to the other points of $h$;
   (v) $x$ is symplectic to every point of $h$;
  (vi) $x$ is symplectic to exactly one point of $h$ and special to the other points of $h$;
 (vii) $x$ is special to all points of $h$;
(viii) $x$ is symplectic to exactly one point of $h$ and opposite all other points of $h$;
  (ix) $x$ is special to exactly one point of $h$ and opposite all other points of $h$.

Also, $x \in S(h)$ in cases $(i), (ii), (iii)$; $x$ is close to $S(h)$ in cases $(iv), (v), (vi), (vii)$, and $x$ is far from $S(h)$ in cases $(viii)$ and $(ix)$.

*Proof* This follows directly from Fact 5.5 and the fact that a hyperbolic line in a symplectic polar space is a so-called *geometric line* of it. We recall that a geometric line in a point-line geometry is a set $g$ of points such that for each point $y$, exactly one point of $g$ is not at maximum distance (in the corresponding collinearity graph) from $y$, or no point of $g$ is at maximum distance from $y$. So for a hyperbolic line $h$ of $\Gamma$, if $x$ is any point of $S(h)$, then either all points of $h$ are collinear with $x$, or exactly one point of $h$ is collinear or equal to $x$. □

**Lemma 5.7** *Let L be a line in $\Gamma$ and $x$ a point. Then, exactly one of the following holds.*

   (i) $x \in L$;
  (ii) $x$ is collinear with every point of $L$;
 (iii) $x$ is collinear with exactly one point of $L$ and symplectic to the other points of $L$;
  (iv) $x$ is collinear with exactly one point of $L$ and special to the other points of $L$;
   (v) $x$ is symplectic to exactly one point of $L$ and special to the other points of $L$;
  (vi) $x$ is special to all points of $L$;
 (vii) $x$ is special to exactly one point of $L$ and opposite all other points of $L$.

*In particular, if two points $x$ and $y$ are opposite, and $L$ is a line through $x$, then $y$ is special to a unique point of $L$ and opposite all other points of $L$.*

*Proof* The proof follows by including $L$ into a symplecton and then consider all possible point-symplecton relations given in Fact 5.5. This way, one can also verify that all listed possibilities indeed occur. As an example, the last statement can be seen as follows. Let $x$ and $y$ be two opposite points and let $L$ be a line through $x$. Take any symplecton $S$ containing $L$. By Fact 5.5, $y$ is far from $S$. Moreover, the unique point $u$ of $S$ symplectic to $y$ is also symplectic to $x$. Hence, there is a unique point $z$ of $L$ which is collinear with $u$, and this point is special to $y$. All other points of $L$ are opposite $y$. In particular, it follows that a point can never be opposite all points of a line. □

### 5.2 The equator and extended equator geometries

In this subsection, we define the equator geometry, see also [11], Proposition 6.26. The structure of an equator geometry is a polar space of type B$_3$ (when endowed with the hyperbolic lines contained in it). It will turn out that each equator geometry is contained in the intersection of $\Gamma$ with a non-absolute quad of the geometry $\Delta$ of type E$_6$ to be defined in Sect. 6. The complete intersection, however, will have the structure of a polar space of type B$_4$. For this reason, we extend the equator geometry to the "extended equator geometry", which is indeed a polar space of type B$_4$, as we will show in Proposition 5.18. We further show that projective subspaces of $\Gamma$ all of whose lines are hyperbolic lines and which are at most 3-dimensional, are contained in an extended equator geometry (Lemma 5.21). Only much later, when our construction of the intended geometry of type $\mathsf{E}_6$ is complete, we will see that there are no such 4-spaces. We will not need to see this earlier.

**Definition 5.8** (*Equator Geometry*) Let $p, q$ be two opposite points of $\Gamma$. Let $\mathscr{S}_p$ denote the family of symplecta containing $p$. Then, by Fact 5.5, each member of $\mathscr{S}_p$ contains a unique point which is symplectic to $q$. The set of all such points is called the *equator geometry of the pair* $\{p, q\}$ and is denoted by $E(p, q)$. Using Fact 5.5(2), it is easy to see that this definition is symmetric in $p, q$.

The following was proved in Proposition 6.26 of [11].

**Proposition 5.9** *Let $p, q$ be two opposite points of $\Gamma$. Then, for any symplectic pair $\{u, v\}$ of points of $E(p, q)$, the hyperbolic line $h(u, v)$ is contained in $E(p, q)$. The geometry of points and hyperbolic lines of $E(p, q)$ is the point-line geometry of a polar space, which we also denote by $E(p, q)$, of type $\mathsf{B}_3$ over the field $\mathbb{K}$, isomorphic to the point residue geometry of $\Gamma$. A natural isomorphism from $E(p, q)$ to $R(p)$ is induced by the map $\varphi_{p,q}$ that sends a point $x \in E(p, q)$ to the symplecton $x \Diamond p$.* □

We will need the following property of polar spaces of type $\mathsf{B}_3$ (which holds for any polar space of rank at least 3).

**Lemma 5.10** *Any geometric hyperplane $G'$ of any geometric hyperplane $G$ of a polar space $\Pi$ of rank at least 3 contains two non-collinear points.*

*Proof* A geometric hyperplane $G'$ of a geometric hyperplane $G$ clearly has the property that each singular plane of $\Pi$ intersects $G'$ nontrivially. If each pair of points of $G'$ is collinear, then it is contained in a maximal singular subspace $U$. Let $U'$ be a maximal singular subspace disjoint from $U$ ($U'$ exists by [10]). Any plane $\pi'$ in $U'$ is disjoint from $G'$, a contradiction. □

The next result determines the mutual relations between any two points of $E(p, q)$.

**Lemma 5.11** *Let $p, q$ be opposite points, and $x, y \in E(p, q)$. Then either $x = y$, or $\{x, y\}$ is a symplectic pair, or $x$ is opposite $y$. All cases occur.*

*Proof* First suppose that $x \perp y \neq x$. Then $x$ is collinear with the point $y$ of the symplecton $p \Diamond y$ and is symplectic to the point $p$, which is not collinear with $y$, in contradiction to Fact 5.5.

Now suppose that $\{x, y\}$ is a special pair. Then, $x$ and $y$ are not collinear in $E(p, q)$. By Fact 5.4 and since $\varphi_{p,q}$ is an isomorphism, the symplecta $x \Diamond p$ and $y \Diamond p$ intersect in just $p$. Put $z = x \bowtie y$. Then $z$ is close to both $x \Diamond p$ and $y \Diamond p$. By Fact 5.5(1), $z$ is collinear with a line $L_x$ through $x$ in $x \Diamond p$, and to a line $L_y$ through $y$ in $y \Diamond p$. Since $p$ is not collinear with $x \in L_x$, there is a unique point $p_x$ on $L_x \setminus \{x\}$ collinear with $p$; likewise there is a unique point $p_y$ on $L_y \setminus \{y\}$ collinear with $p$. Now, since $p$ is not collinear with $L_x$, it is special to $z$ by Fact 5.5. But both $p_x$ and $p_y$, which are distinct, are contained in $p^\perp \cap z^\perp$, a contradiction. □

We are now ready to define the extended equator geometry for opposite points $p, q$.

**Definition 5.12** (*Extended Equator Geometry*) Let $p, q$ be two opposite points of $\Gamma$. Then define the point set

$$\widehat{E}(p, q) = \bigcup \{E(x, y) : x, y \in E(p, q), x \text{ opposite } y\}.$$

Note that, by Proposition 5.9 and Lemma 5.11, $E(p, q)$ contains pairs of opposite points. So, $\widehat{E}(p, q)$ is nonempty. The set $\widehat{E}(p, q)$, endowed with all the hyperbolic lines in it, is called the *extended equator geometry* for $p, q$. Further, $p, q$ and $E(p, q)$ are contained in $\widehat{E}(p, q)$. The latter follows from Proposition 5.9 and the trivial fact that every point of a polar space of rank 3 is collinear to two non-collinear points.

**Standing hypothesis.** From now on until Sect. 6.5 (included), we fix a pair of opposite points $p, q$ and write $\widehat{E} := \widehat{E}(p, q)$.

We now prove that $\widehat{E}$ does not contain collinear or special pairs of points, that it is closed under taking hyperbolic lines through symplectic pairs of its points, and that the geometry of its points and hyperbolic lines is the point-line geometry of a polar space of type B$_4$.

**Lemma 5.13** *No point $x \in \widehat{E}$ is collinear with any point of $E(p, q)$.*

*Proof* Suppose, for a contradiction, that $x \perp y \in E(p, q)$. Then, $x \notin E(p, q)$ by Lemma 5.11 and $x \notin \{p, q\}$ by the definition of $E(p, q)$. Let $x \in E(a, b)$, with $a, b \in E(p, q)$ and $a$ opposite $b$. Since $y \neq x$ and $(a \lozenge x) \cap (b \lozenge x) = \{x\}$, we may assume that $y$ does not belong to $a \lozenge x$. As $y \perp x$, $y$ is close to $a \lozenge x$ by Fact 5.5. Then $a$ cannot be opposite $y$ by Fact 5.5(1), so $a$ is symplectic to $y$, which implies that $a \perp x$ by the same reference, a contradiction. □

**Lemma 5.14** *No point $x \in \widehat{E}$ is special to any point of $E(p, q)$.*

*Proof* Suppose that $x \in \widehat{E}$ is special to the point $y \in E(p, q)$. Then $x \notin E(p, q) \cup \{p, q\}$. Let $a, b \in E(p, q)$ be opposite points such that $x \in E(a, b)$. Then, $a \neq y \neq b$. We first show that $y$ is opposite both $a$ and $b$. Indeed, suppose that $y$ is symplectic to $a$. Since $x$ is symplectic to $a$ and special to $y$, it follows from Lemma 5.6 that all points of $h(a, y) \setminus \{a\}$ are special to $x$. Lemma 5.6 implies that not all points of $h(a, y)$ can be opposite $b$ and Lemma 5.11 then tells us that there is a point of $h(a, y) \setminus \{a\}$ symplectic to $b$, and we can rename this point $y$. So we may assume that $y$ is symplectic to both $a$ and $b$. Then, $x, y \in E(a, b)$ and the assertion follows from Lemma 5.11.

Hence we know that $y$ is opposite both $a$ and $b$. Moreover, the above argument implies that no point in $E(p, q)$ symplectic to $a$ or $b$ is special to $x$. Hence, by Lemmas 5.6 and 5.11, the set of points $H_{x,a} = a^\perp \cap x^\perp \cap E(p, q)$ is a geometric hyperplane of $a^\perp \cap E(p, q)$. The latter is a geometric hyperplane of $E(p, q)$ and hence, by Lemma 5.10, $H_{x,a}$ contains two opposite points $a', b'$. If $y$ is symplectic to both of these, then the first part of the proof applies with $a', b'$ in the roles of $a, b$, respectively. Hence we may assume that $y$ and $a'$ are opposite. By Proposition 5.9, there is a point $a'' \in h(a, a')$ symplectic to $y$. Since $b$ is opposite $a$ and $b \perp\!\!\!\perp a'$, we know from Lemma 5.6 that $a''$ is opposite $b$. Recalling that $a'' \in H_{x,a}$ and as such $x \in E(a'', b)$, the first part of the proof implies that $y$ is opposite $a''$, contradicting the choice of $a''$. □

There is an interesting corollary to the previous three lemmas.

**Corollary 5.15** *Let $x \in \widehat{E}$. Then the set of points of $E(p, q)$ symplectic to or equal to $x$ is a geometric hyperplane of $E(p, q)$, viewed as a polar space, or coincides with it.*

*Proof* Let $h$ be a line of the polar space $E(p, q)$ (so $h$ is a hyperbolic line contained in $E(p, q)$). Then, by Lemma 5.6, either exactly one, or all points of $h$ are not opposite $x$. By Lemmas 5.13 and 5.14, $x$ can neither be collinear nor special to a point of $h \subseteq E(p, q)$. So either exactly one, or all points of $h$ are equal or symplectic to $x$. This completes the proof of the corollary. □

*Remark 5.16* It will follow from Proposition 5.18 that, if the characteristic of $\mathbb{K}$ is not equal to 2, then the only points of $\widehat{E}(p, q)$ which are symplectic to all points of $E(p, q)$ are $p$ and $q$. Indeed, in this case, $\widehat{E}(p, q)$ is a polar space of type $\mathsf{B}_4$ arising from an orthogonal polarity $\rho$ in $\mathsf{PG}(8, \mathbb{K})$. The set of points of $\widehat{E}(p, q)$ symplectic to all points of $E(p, q)$ is the image under $\rho$ of the codimension 2 space $\mathscr{A}$ of $\mathsf{PG}(8, \mathbb{K})$ generated by the points of $E(p, q)$. This is the line of $\mathsf{PG}(8, \mathbb{K})$ through $p$ and $q$, which intersects $\widehat{E}(p, q)$ in just $\{p, q\}$. If the characteristic of $\mathbb{K}$ is 2, then the associated polarity $\rho$ is symplectic in $\mathsf{PG}(7, \mathbb{K})$ and the image of $\mathscr{A}$, similarly defined as above, is a line of $\mathsf{PG}(7, \mathbb{K})$ all of whose points are contained in $\widehat{E}(p, q)$. Hence it contains at least three points of $\widehat{E}(p, q)$.

**Lemma 5.17** *Let* $x, y \in \widehat{E}$, $x \neq y$. *Then either* $\{x, y\}$ *is a symplectic pair, or* $x$ *is opposite* $y$. *If, moreover,* $\{x, y\}$ *is a symplectic pair, then* $h(x, y)$ *is contained in* $\widehat{E}$.

*Proof* By Lemma 5.10 and Corollary 5.15, we can find two opposite points $a, b \in E(p, q)$ symplectic to both $x$ and $y$. Hence $x, y \in E(a, b)$ and so the first assertion follows from Lemma 5.11. If, moreover, $x$ and $y$ are symplectic, then $h(x, y) \subseteq E(a, b) \subseteq \widehat{E}$, by Proposition 5.9 and the definition of $\widehat{E}$. □

We now have the following interesting proposition.

**Proposition 5.18** *The extended equator geometry* $\widehat{E}(p, q)$ *is a polar space of type* $\mathsf{B}_4$ *over the field* $\mathbb{K}$.

*Proof* We check the Buekenhout–Shult axioms of a polar space as given in [6]. We repeat these axioms for the convenience of the reader.

(1) *Every (hyperbolic) line contains at least* 3 *points.* This holds by Lemma 5.17 and the fact that a hyperbolic line contains at least 3 points.
(2) *There is no point collinear with every other point.* By definition of $\widehat{E}(p, q)$, any point $x \in \widehat{E}(p, q)$ is contained in an equator geometry, which is, by Proposition 5.9, a polar space of type $\mathsf{B}_3$, in which $x$ has an opposite point.
(3) *One-or-all axiom, i.e., either exactly one or all points of a given line are collinear with a given point.* This follows from Lemmas 5.17 and 5.6.
(4) *Finite rank, i.e., every nested family of singular subspaces is finite.* Again by Proposition 5.9, the residue in the point $p$ is isomorphic to the polar space induced on the set of points of $\widehat{E}(p, q)$ symplectic to both $p$ and $q$. Since in the whole of $\Gamma$, this is $E(p, q)$, it is also $E(p, q)$ in $\widehat{E}$. Hence the residue at $p$ is a polar space of type $\mathsf{B}_3$ and as such has rank 3. We conclude that the rank of $\widehat{E}(p, q)$ is 4 and hence finite.

The argument above implies that $\widehat{E}(p, q)$ is a polar space of type $\mathsf{B}_4$, as the residue in at least one point has type $\mathsf{B}_3$. The proposition is proved. □

The following is a straightforward consequence.

**Corollary 5.19** *A maximal singular subspace of* $\widehat{E}(p, q)$ *is a projective* 3-*space.* □

**Definition 5.20** A set of points of $\Gamma$ which is a projective space of dimension $i$, $i = 1, 2, 3$, when endowed with the hyperbolic lines it contains shall be called a *hyperbolic $i$-space* or *hyperbolic line* (if $i = 1$), *hyperbolic plane* (if $i = 2$), or *hyperbolic solid* (if $i = 3$) of $\Gamma$.

In view of the previous definition, we shall speak of singular *hyperbolic* subspaces of an extended equator geometry, and the meaning is clear.

**Lemma 5.21** *(i) Every hyperbolic line of $\Gamma$ is contained in a hyperbolic plane of $\Gamma$.*
*(ii) Every hyperbolic plane of $\Gamma$ is contained in a hyperbolic solid of $\Gamma$.*
*(iii) Every hyperbolic solid of $\Gamma$ is contained in an extended equator geometry of $\Gamma$.*

*Proof* We first provide an outline of the proof. Let $h$ be a hyperbolic line.

*Step 1* We construct a point $x$ symplectic to all points of $h$.
*Step 2* Given a hyperbolic plane $\alpha$ containing $h$, we show that each point of $\alpha \setminus h$ is as in the construction in the proof of Step 1.
*Step 3* For $x$ constructed in Step 1, we construct a point $p$ symplectic to each point of $\pi(x, h) := \{y \in h(x, z) : z \in h\}$. Note that, if $x$ and $h$ are contained in a common hyperbolic plane $\pi$, then $\pi(x, h) = \pi$.
*Step 4* Given a hyperbolic solid $S$ containing $\pi(x, h)$, we show that each point of $S \setminus \pi(x, h)$ is as in the construction in the proof of Step 3.
*Step 5* For $p$ constructed in Step 3, we construct an extended equator geometry that contains $S(p, x, h) := \{y \in h(p, z) : z \in \pi(h, x)\}$. Note that, if $x$, $p$ and $h$ are contained in a common solid $S$, then $S(p, x, h) = S$.

It then follows from Proposition 5.18 that $\pi(x, h)$ is a hyperbolic plane (showing $(i)$), and that $S(p, x, h)$ is a hyperbolic solid. In view of Steps 2 and 4, the Steps 3 and 5 apply to *each* hyperbolic plane (showing $(ii)$) and to *each* hyperbolic solid (showing $(iii)$).

Now we embark on the proof.

*Step 1* Let $h$ be a hyperbolic line of $\Gamma$. Since $S(h)$ is convex by the definition of parapolar space, $h^\perp \subseteq S(h)$. Since $S(h)$ is a polar space of rank 3, we can select a line $L$ in $S(h)$ contained in $h^\perp$. By Fact 5.1, we can select a plane $\pi$ of $\Gamma$ containing $L$ but not contained in $S(h)$. Let $x$ be any point of $\pi \setminus L$. Then Fact 5.5 implies that $x$ is symplectic to every point $u$ of $h$, and $L \subseteq x \lozenge u$.
*Step 2* Given a hyperbolic plane $\alpha$ and a hyperbolic line $h$ in it, each point $x$ of $\alpha$ not in $h$ is also not in $S(h)$ and, by Fact 5.5, is close to $S(h)$. Hence, $x^\perp \cap S(h)$ is a line $L$ contained in $h^\perp$.
*Step 3* Next, given $h$ and $x$ as in Step 1 above, we construct a point $p$ symplectic to each element of $\pi(x, h) = \{y \in h(x, z) : z \in h\}$. Let $u_1, u_2 \in h$, $u_1 \ne u_2$. With $L$ and $\pi$ as in the proof of Step 1 ($L$ and $\pi$ are uniquely determined by $h$ and $x$) and $i \in \{1, 2\}$, we have that $x \perp L \perp u_i$, so $x \lozenge u_1$ and $x \lozenge u_2$ intersect in the plane $\pi$. By Fact 5.1, we can select a symplecton $S_1$ containing $x$ and intersecting $\pi$ in a line $M_1$. Since $x \notin L$, the lines $L$ and $M_1$ intersect in a point $y_1$ distinct from $x$. Since $u_i \notin S_1$ (otherwise $S_1 = x \lozenge u_i$ contains $\pi$, a contradiction) and $u_i$ is collinear with $y_1 \in L$, $u_i$ is close to $S_1$. Let $L_i$ be the line in $S_1$ collinear with $u_i$. Then, all points of $L_i$ are collinear with $x$, as $x \perp\!\!\!\perp u_i$. Let $\pi_i$ be the plane of $\Gamma$ spanned by $x$ and $L_i$. Then, $\pi_1 \ne \pi_2$, as otherwise $\pi_i \subseteq x \lozenge u_i$ implies that $\pi_1 = \pi_2 \subseteq x \lozenge u_1 \cap x \lozenge u_2 = \pi$, contradicting the fact that $S_1$ intersects $\pi$ in a line. Let $p$ be any point of $S_1$ not collinear with $x$ but collinear with all points of both $L_1$ and $L_2$ (such a point exists as $L_1^\perp \cap L_2^\perp$ is a plane $\alpha$ in the projective 5-space underlying $S_1$ containing $x$, consisting of hyperbolic lines not containing $y_1$ and ordinary lines containing $y_1$. Hence $L_1^\perp \cap L_2^\perp \cap x^\perp = xy_1$ and $p$ can be chosen in $\alpha \setminus xy_1$ arbitrarily). Then, $p$ is symplectic to $x, u_1, u_2$. Note that $x \lozenge p$ does not intersect $h$ because otherwise $x \lozenge p$ would have to contain $L$. Consequently, by Lemma 5.6, $p$ is symplectic to each point of $\pi(x, h)$.
*Step 4* If $p$ is a point of a hyperbolic solid containing $\pi(x, h)$, with $p \notin \pi(x, h)$, then $x \lozenge p$ is disjoint from $h$. Indeed, otherwise $(x \lozenge p) \cap \pi(x, h)$ would contain a point collinear with $p$, a contradiction. Let $\pi$ and $L$ be defined as in the proof of Step 2 (with $\pi(x, h)$ playing the role of $\alpha$). Assume, for a contradiction, that $x \Diamond p$ contains $\pi$. In particular, $x \Diamond p$ contains $L$. However, every symplecton $S \neq S(h)$ through $L$ intersects $S(h)$ in a plane by Fact 5.4 and this plane intersects $h$, since $L^\perp$ is 3-dimensional in the 5-space underlying $S(h)$, and $h \subseteq L^\perp$. Hence $x \Diamond p$ would intersect $h$, a contradiction. Thus, $\pi$, and so $L$, is not contained in $x \Diamond p$. Consequently, as $x \Diamond p$ is convex, the point $p$ is not collinear with $L$. Now suppose that no point of $L$ is collinear with $p$. Then the line $N$ in $S(h)$ collinear with $p$ does not intersect $L$ ($p$ is indeed close to $S(h)$ as it is symplectic to each point of $h$). Also, every point $w \in L$ is special to $p$, as $w$ cannot be collinear with $N$ because $h^\perp$ does not contain planes (and it would contain $w$ and $N$). The point $w \bowtie p$ then belongs to $N$ and hence is collinear with all points of $h$. Let $u_1$ and $u_2$ be two distinct points of $h$. Then, likewise, $w \bowtie p$ belongs to the line $N_i$ in $x \Diamond u_i$ collinear with $p$, $i = 1, 2$. It follows that $w \bowtie p$ belongs to $x \Diamond u_1 \cap x \Diamond u_2 \cap S(h) = \pi \cap S(h) = L$. We conclude that $p$ is collinear with a unique point of $L$ after all, and we denote this point as $y_1$, emphasising the similarity with the construction in Step 3. It follows that $x \Diamond p$ intersects $\pi = (x \Diamond u_1) \cap (x \Diamond u_2)$ in the line $xy_1 = M_1$. With $x \Diamond p$ playing the role of $S_1$ above, we have shown that $p$ is constructed as in Step 3.

*Step 5* To find an extended equator geometry containing $S(p, x, h)$, we select a symplecton $S_2$ containing $x$ and intersecting $\pi$ in a line $M_2 \neq M_1$. Set $\{y_2\} = L \cap M_2$. The same argument as above provides a point $q \in S_2$ symplectic to each of $x, u_1, u_2$. We claim that $p$ is far from $S_2$. Indeed, if not, then $p$ is collinear with all points of a line $L'$ of $S_2$. Since $p \perp\!\!\!\perp x$, the line $L'$ is contained in $x^\perp \cap p^\perp$ and hence in $S_1$. Hence $S_1 \cap S_2 = \langle x, L' \rangle$. Since $y_1$ is collinear with a line of $S_2$ (containing $y_2$), and also, in the polar space $S_1$, with a line of $\langle x, L' \rangle$, these two lines coincide and so $y_2 \in S_1$, a contradiction. Hence the claim follows. But then, since $p \perp\!\!\!\perp x \perp\!\!\!\perp q$, we see that $p$ is opposite $q$. By construction, $x, u_1, u_2 \in E(p, q)$ and so $x, u_1, u_2, p \in \widehat{E}(p, q)$, and the latter is the promised extended equator geometry. □

Next we show that $\widehat{E}(p, q)$ is independent of the choice of the pair $\{p, q\}$ of opposite points in it. This and Lemma 5.17 imply that $\widehat{E}(p, q)$ is a convex subspace of $\Gamma$ relative to the hyperbolic lines of $\Gamma$ (see also Proposition 5.18).

**Proposition 5.22** *Let $a, b \in \widehat{E}(p, q)$ be two opposite points. Then, $\widehat{E}(p, q) = \widehat{E}(a, b)$.*

*Proof* We start by showing that $E(a, b) \subseteq \widehat{E}(p, q)$. Let $Y$ be the set of points of $\widehat{E}(p, q)$ symplectic to both $a$ and $b$. By Proposition 5.18, $Y$, endowed with all hyperbolic lines it contains, is a polar space of type $\mathsf{B}_3$, naturally contained in $E(a, b)$, which is also a polar space of type $\mathsf{B}_3$. Take two opposite points $x, y \in Y$. Then $x^\perp \cap y^\perp \cap Y$ is a polar space of type $\mathsf{B}_2$ naturally contained in the polar space $x^\perp \cap y^\perp \cap E(a, b)$; in fact this containment is "full", in the sense that common "lines" of both polar spaces (which are hyperbolic lines in $\Gamma$) have the same point sets. Consequently it follows from Proposition 5.9.4 of [26] that $x^\perp \cap y^\perp \cap Y = x^\perp \cap y^\perp \cap E(a, b)$, which readily implies $Y = E(a, b)$ by varying $x, y$. Hence $E(a, b) \subseteq \widehat{E}(p, q)$.

Now an arbitrary point of $\widehat{E}(a, b)$ is contained in $E(x, y)$ for some opposite points $x, y \in E(a, b)$. By the previous paragraph applied to $x, y$ instead of $a, b$, we know that $E(x, y) \subseteq \widehat{E}(p, q)$. Hence we have shown that $\widehat{E}(a, b) \subseteq \widehat{E}(p, q)$.

Now note that $E(a, b) \cap E(p, q)$ is the geometric hyperplane $H := a^\perp \cap b^\perp \cap E(p, q)$ of the geometric hyperplane $b^\perp \cap E(p, q)$ of $E(p, q)$, which is a polar space of type $\mathsf{B}_3$. By Lemma 5.10, $H$ contains two opposite points $x, y$. But then $p, q \in E(x, y) \subseteq \widehat{E}(a, b)$.

Hence, by the previous paragraph, switching the roles of $a, b$ and $p, q$, we obtain $\widehat{E}(p, q) \subseteq \widehat{E}(a, b)$. □

**Corollary 5.23** *Let $S$ be an arbitrary symplecton. Then, either $S$ is disjoint from $\widehat{E}(p, q)$ or $S \cap \widehat{E}(p, q)$ is a hyperbolic line. Hence, every symplecton that has a point $x$ in common with $\widehat{E}(p, q)$ intersects it in a hyperbolic line through $x$. In particular, any hyperbolic line in $\widehat{E}$ appears as the intersection of $\widehat{E}$ and a unique symplecton.*

*Proof* By Propositions 5.18 and 5.22, we may assume that $p \in S \cap \widehat{E}(p, q)$. By the definition of equator geometry, there is a unique point $a \in S$ in $E(p, q)$. But then, by Lemma 5.17, $h(p, a) \subseteq \widehat{E}(p, q) \cap S$. No other point of $\widehat{E}(p, q)$ is contained in $S$ as that point would then be collinear with at least one point of $h(p, a)$, contradicting Lemma 5.17. As a hyperbolic line defines a unique symplecton containing it, the rest of the corollary follows. □

### 5.3 The tropic circle geometries

A tropic circle geometry is related to an extended equator geometry. In the building of type E$_6$ we aim to construct, the first one is the set of points of $\Gamma$ collinear to a "new point", while the latter is the intersection of $\Gamma$ with the quad which is the image under the symplectic polarity of the said "new point". The intrinsic geometric connection between these two geometries is the fact that they are dual to each other, see Theorem 5.33.

The notion of tropic circle geometry is based on the following property of extended equator geometries. We continue with our notation $\widehat{E} = \widehat{E}(p, q)$ for two fixed opposite points $p, q$ in $\Gamma$.

**Proposition 5.24** *Let $x$ be a point of $\Gamma$ which is collinear with at least two points of $\widehat{E}$. Then $x^\perp \cap \widehat{E}$ is a hyperbolic solid.*

*Proof* By Propositions 5.18 and 5.22, we may assume that $p \perp x$. Let $a$ be a second point of $\widehat{E}$ collinear with $x$. By Lemma 5.17, $p \perp\!\!\!\perp a$. Hence, by Propositions 5.18 and 5.22, we can choose $q$ opposite $p$ and symplectic to $a$. Hence $a \in E(p, q)$. Then, by Fact 5.5, $x \in p \Diamond a$. By the canonical isomorphism $\varphi_{p,q}$, the set of intersections with $E(p, q)$ of the symplecta through $p$ and $x$ is a hyperbolic plane $\pi$ of $E(p, q)$. Let $b \in \pi$ and suppose $b \neq a$. Since $a$ is collinear with $x$ and $x \in p \diamond b$, the point $a$ is close to $p \Diamond b$. Since $b \in a^{\perp\!\!\!\perp}$, Fact 5.5 implies that $x \perp b$.

Hence all points $u$ of $\pi$ are collinear with $x$. But $x$ belongs to $u \Diamond p$, and in the latter symplectic polar space, $u$ and $p$ belong to $x^\perp$; hence, by the definition of the hyperbolic line $h(u, p)$, all points of $h(u, p)$ are collinear with $x$, implying that all points of the maximal singular hyperbolic subspace of $\widehat{E}(p, q)$ spanned by $\pi$ and $p$ are collinear with $x$. Since Lemma 5.17 implies that the set of points of $\widehat{E}$ collinear with $x$ is contained in a maximal singular subspace of $\widehat{E}$, there are no other points of $\widehat{E}$ collinear with $x$. The assertion is proved. □

**Definition 5.25** (*Tropic Circle Geometry*) The point set

$$\widehat{T}(p, q) = \left\{ x \in \Gamma : |x^\perp \cap \widehat{E}(p, q)| \geq 2 \right\},$$

endowed with all lines inside it, is called the *tropic circle geometry for* $\{p, q\}$.

Note that $\widehat{E}(p, q) \cap \widehat{T}(p, q) = \emptyset$. In particular, $p, q \notin \widehat{T}(p, q)$. We write $\widehat{T}$ instead of $\widehat{T}(p, q)$ if $\{p, q\}$ is understood.

**Definition 5.26** For $x \in \widehat{T}$, we denote by $\beta(x)$ the hyperbolic solid $x^\perp \cap \widehat{E}$.

**Corollary 5.27** *No point of $\widehat{T}(p, q)$ is opposite any point of $\widehat{E}(p, q)$.*

*Proof* Let $x \in \widehat{T}(p, q)$, and let $y \in \widehat{E}(p, q)$ be arbitrary. Then $y$ is symplectic to at least one point $z \in \beta(x)$. The last assertion of Lemma 5.7 now implies that $x$ ad $y$ are not opposite. □

In Proposition 5.32 below, we show that there is a neat connection between the dimension of the intersection of two maximal singular subspaces of $\widehat{E}(p, q)$ and the mutual position in $\Gamma$ of the corresponding points of $\widehat{T}(p, q)$. This will imply that the map taking $x \in \widehat{T}(p, q)$ to the hyperbolic solid $\beta(x)$ of $\widehat{E}(p, q)$ is an isomorphism from $\widehat{T}(p, q)$ to the dual polar space structure associated with $\widehat{E}(p, q)$.

**Lemma 5.28** *Let $x, y \in \widehat{T}$. If $\beta(x) = \beta(y)$, then $x = y$.*

*Proof* We may again suppose that $p$ belongs to $\beta(x)$. Then, $\beta(x)$ intersects $E(p, q)$ in a hyperbolic plane $\pi$. The intersection of all symplecta $S$ such that $\varphi_{p,q}^{-1}(S) \in \pi$ is, by the fact that $\varphi_{p,q}$ is an isomorphism of geometries (see Proposition 5.9), a line $L$ through $p$. Now, both $x$ and $y$ must be contained in all these symplecta, hence both are on $L$. Let $z \in \pi$ be arbitrary. Then in $z \lozenge p$, the point $z$ is collinear with exactly one point of $L$. Thus $x = y$. □

**Lemma 5.29** *Let $U$ be a hyperbolic solid of $\widehat{E}$. Then there is a unique point $x \in \widehat{T}$ with $\beta(x) = U$. Moreover, this is the only point of $\Gamma$ collinear with $U$.*

*Proof* By Lemma 5.28, we only need to prove the existence of $x$. We may suppose that $p \in U$. Then $U \cap E(p, q)$ is a hyperbolic plane $\pi$. As in the previous proof, there is a unique line $L$ through $p$ contained in all symplecta defined by $p$ and a point of $\pi$. Let $a, b \in \pi$ be arbitrary but distinct. Then, $b$ is not contained in $a \lozenge p$ and hence is close to it. So $b$ is collinear with a line $M \subseteq a \lozenge p$ and Fact 5.5 implies that $a$ and $p$ are also collinear with $M$. Clearly, $L$ is contained in the plane generated by $p$ and $M$, and so $\{x\} = L \cap M$ is collinear with both $a$ and $b$. Since $x$ is the unique point of $L$ collinear with $a$, we see, by varying $b \in \pi$, that $x$ is collinear with all points of $\pi$. Since $x \perp p$, the first assertion follows. As any other point in $\Gamma$ that is collinear with $U$ would belong to $\widehat{T}$, the second assertion also follows. □

The previous lemma shows that $\beta$ is bijective. We denote its inverse again by $\beta$. There is also an interesting corollary.

**Corollary 5.30** *For any hyperbolic solid $U$ in $\Gamma$, there is a unique point of $\Gamma$ collinear with $U$.*

*Proof* By Lemma 5.21, there is an extended equator geometry containing $U$. By Lemma 5.29, there is a unique point in $\Gamma$ collinear with $U$. □

**Definition 5.31** For any hyperbolic solid $U$, the unique point collinear with $U$ is denoted by $\beta(U)$.

We now relate the mutual position of two hyperbolic solids of $\widehat{E}$ to the mutual position of their images under $\beta$.

**Proposition 5.32** *Let $U$ and $V$ be hyperbolic solids in $\widehat{E}$. Then*

(i) *$U$ and $V$ intersect in a hyperbolic plane $\pi$ if, and only if, $\beta(U)$ and $\beta(V)$ are collinear in $\Gamma$. In this case, every point of the line of $\Gamma$ joining $\beta(U)$ and $\beta(V)$ belongs to $\widehat{T}$ and is collinear with all points of $\pi$. Also, if some point is collinear with all points of $\pi$, then it belongs to the line joining $\beta(U)$ and $\beta(V)$. Consequently, $\widehat{T}$ is a subspace of $\Gamma$.*

(ii) $U \cap V$ is a hyperbolic line if, and only if, $\beta(U)$ and $\beta(V)$ are symplectic in $\Gamma$. In this case, every point of $h(\beta(U), \beta(V))$ belongs to $\widehat{T}$ and is collinear with all points of $U \cap V$.

(iii) $U \cap V$ is a singleton $\{z\}$ if, and only if, $a = \beta(U)$ and $b = \beta(V)$ are special in $\Gamma$. In this case, $z = a \bowtie b$.

(iv) $U \cap V = \emptyset$ if, and only if, $a = \beta(U)$ and $b = \beta(V)$ are opposite in $\Gamma$.

*Proof* (i) Put $a = \beta(U)$, $b = \beta(V)$, and suppose first that $a \perp b$ (so $a \neq b$). Since $U \neq V$ by Lemma 5.28, the union $U \cup V$ contains a pair of opposite points. Hence by Lemma 5.18 we may assume $p \in U \setminus V$ and $q \in V \setminus U$. Thus we have $p \perp a \perp b \perp q$, with $p$ opposite $q$. The last assertion of Lemma 5.7 implies that $\{p, b\}$ and $\{q, a\}$ are special pairs. Let $S$ be any symplecton through $q$ and $b$, and let $\{x\} = E(p, q) \cap S$. Then $\{p, x\}$ is symplectic and $p$ and $S$ are far. Since $p$ is special to $b$, Fact 5.5 implies $b \perp x$. Now consider $T = p \Diamond x$. Since $b \perp x$, $b$ is close to $T$. Hence there is a line $L$ in $T$ containing $x$ such that $L$ is collinear with $b$. Moreover, $b \bowtie p$ is contained in $L$ (see Fact 5.5). Since $b \bowtie p = a$, we see that $a \perp x$. Varying $S$ over all symplecta containing $q$ and $b$ and using the isomorphism $\varphi_{q,p}$, the point $x$ varies over a plane of $E(p, q)$, which must coincide with $U \cap V$ since $x \perp a$ and $x \perp b$.

By Lemma 5.7, no point of the line $ab$ is symplectic to or opposite $p$. Lemma 5.11 then implies that the line $ab$ has empty intersection with $\widehat{E}$.

Let $z$ be any point of $U \cap V$. Then $a \perp z \perp b$, and so every point of the line $ab$ is collinear with $z$ and hence with all points of $\pi$.

Now assume that $U$ and $V$ intersect in a plane $\pi$. Then we can assume that $p \in U \setminus \pi$ and $q \in V \setminus \pi$; $p$ and $q$ are opposite. Hence $\pi \subseteq E(p, q)$. Consider two points $x, y \in \pi$. Then both $a$ and $b$ are collinear with both $x, y$ and hence both are contained in $x \Diamond y$. It follows that $a, b$ are either symplectic or collinear. If they were symplectic, then $a \Diamond b$ would contain $\pi$, a contradiction.

From the first part of the proof, we already know that every point of the line $ab$ is collinear with all points of $\pi$. Now suppose some point $c$ is collinear with all points of $\pi$. Then $c \in \widehat{T}$ and we have just shown that $a \perp c \perp b$. Suppose $c$ does not belong to the line $ab$, then take two points $u, v \in \pi$. It follows that $a, b, c \in u \Diamond v$, contradicting the fact that $u \Diamond v$ is a polar space of rank 3 and hence no plane can be contained in the intersection of perps $u^{\perp} \cap v^{\perp}$.

(ii) Again put $a = \beta(U)$ and $b = \beta(V)$. Assume first that $U$ and $V$ intersect in a hyperbolic line $h$. We can again assume that $p \in U \setminus h$ and $q \in V \setminus h$. Then, $h \subseteq E(p, q)$. Consider two points $x, y \in h$. Then both $a$ and $b$ are collinear with both $x, y$ and hence contained in $x \Diamond y$. It follows that $a, b$ are either symplectic or collinear. But they are not collinear by (i).

Now assume that $\{a, b\}$ is a symplectic pair. Then by Lemma 5.28 and (i), we know that $U \cap V$ contains at most a hyperbolic line. Hence we may again assume that $p \in U \setminus V$ and $q \in V \setminus U$. Then both $p$ and $q$ are close to $a \Diamond b$. Hence $p$ is collinear with the points of a line $L \subseteq a \Diamond b$, and $q$ is collinear with the points of a line $M \subseteq a \Diamond b$. If some point $u$ of $L$ was collinear with all points of $M$, then $\{u, q\}$ would be a symplectic pair with $p$ close to $u \Diamond q$, contradicting Fact 5.5 and the fact that $p$ and $q$ are opposite. Hence $L$ and $M$ are, viewed as lines of the symplectic polar space $a \Diamond b$, opposite lines. This implies that $a \Diamond b$ contains a unique hyperbolic line $h$ all of whose points are collinear with $L$ and $M$, i.e., $h = L^{\perp} \cap M^{\perp}$. In particular, $h$ is contained in $a^{\perp} \cap b^{\perp}$. By Fact 5.5, all points of $h$ are symplectic to both $p$ and $q$, hence $h \subseteq E(p, q)$. So $h \subseteq U \cap V$, implying $h = U \cap V$.

(iii) Suppose first that $U$ and $V$ intersect in a point. Then $a$ and $b$ are collinear with a common point and hence cannot be opposite. Moreover, they are neither symplectic nor collinear by $(i)$ and $(ii)$. Consequently, they are special.

Now suppose that $a$ and $b$ are special. We show that $z = a \bowtie b$ belongs to $\widehat{E}$, which will complete the proof of $(iii)$.

Suppose $z \notin \widehat{E}$. Then $U \cap V = \emptyset$. Let $h$ be a hyperbolic line in $U \cup V$. Let $S$ be the unique symplecton containing $h$. By convexity, $S$ contains either $a$ (if $h \subseteq U$) or $b$ (if $h \subseteq V$). In any case, since $z \perp a, b$, we see that $z$ is either close to $S$ or contained in $S$. Suppose first that $z$ is close to $S$. Then $z$ is collinear with all points of a line $L \subseteq S$. Note that $a$ or $b$ is on $L$. Let $c$ be a point of $L$ distinct from $a$ and $b$. Since $h \subseteq a^\perp$ (if $a \in L$) or $h \subseteq b^\perp$ (if $b \in L$), we have that $L^\perp \cap h = c^\perp \cap h$ is either a singleton or $h$. It follows from Fact 5.5 that either one or all points of $h$ are symplectic to $z$.

Now suppose $z \in S$. If all points of $h$ are collinear with $z$, then $z \in \widehat{T}$, which is impossible since $\beta(z) \cap \beta(a)$ and $\beta(z) \cap \beta(b)$ are hyperbolic planes in $\beta(z)$ by $(i)$, contradicting $\beta(a) \cap \beta(b) = \emptyset$. Hence, by Lemma 5.6, $z$ is collinear with a unique point of $h$ and symplectic to the other points of $h$. Now notice that, if $z$ would be contained in a symplecton defined by any other hyperbolic line $h'$ in $U \cup V$ disjoint from $h$, then the same argument implies that $z$ is collinear with a point of $h'$, but then again $z \in \widehat{T}$, a contradiction.

Since $x$ is collinear to at most one point of $U \cup V$, we may hence assume that $x$ is collinear to no point of $U$. So, the previous arguments imply that $z^{\perp\!\!\!\perp} \cap U$ contains a hyperbolic plane $H$ and $z^{\perp\!\!\!\perp} \cap V$ contains at least two points $u, v$. It is easy to see that $H$ contains a point $x$ not symplectic to one of $u, v$, say $u$, as otherwise $H$ and $h(u, v)$ would be contained in a singular hyperbolic subspace of $\widehat{E}$ with dimension at least 4, a contradiction. But then, $x$ and $u$ are opposite by Lemma 5.17 and $z \in E(x, u) \subseteq \widehat{E}(p, q)$.

(iv) This follows by elimination and the previous cases. □

Recall that the dual polar space associated with a polar space $\Omega$ of type $\mathsf{B}_4$ is naturally embedded in the half spin geometry associated to a polar space $\Omega'$ of type $\mathsf{D}_5$ as each 3-space $M$ of $\Omega$ is contained in a unique member $M'$ of one of the families of maximal singular subspaces of $\Omega'$. Two maximal singular subspaces $M'_1$ and $M'_2$ of $\Omega'$ intersecting in a singular 3-space of $\Omega'$ correspond to distinct members $M_1$ and $M_2$ of $\Omega$ that intersect each other in at least a line. The previous lemmas now readily imply the following geometric connection between $\widehat{E}(p, q)$ and $\widehat{T}(p, q)$.

**Theorem 5.33** *(i) The set $\widehat{T}(p, q)$ endowed with the lines of $\Gamma$ contained in it is isomorphic to the dual polar space of type $\mathsf{B}_4$ over the field $\mathbb{K}$.*

*(ii) The set $\widehat{T}(p, q)$ endowed with the lines and hyperbolic lines of $\Gamma$ contained in it is isomorphic to the half spin geometry of type $\mathsf{D}_5$ over the field $\mathbb{K}$.* □

By definition, $\widehat{T}(p, q)$ is, as a set of points, uniquely determined by $\widehat{E}(p, q)$. By Proposition 5.32$(iii)$, the set of points $\{a \bowtie b \mid a, b \in \widehat{T}(p, q) \text{ and } a \text{ special to } b\}$ coincides with $\widehat{E}(p, q)$. Thus, $\widehat{E}(p, q)$ is also determined by $\widehat{T}(p, q)$ as a set of points. Note that the structure of $\widehat{E}(p, q)$ as a polar space of type $\mathsf{B}_4$ and the structure $\widehat{T}(p, q)$, both as a dual polar space of type $\mathsf{B}_4$ and as a half spin geometry of type $\mathsf{D}_5$, are inherited from $\Gamma$. For $\widehat{T}(p, q)$, this follows from Theorem 5.33.

**Definition 5.34** (*Imaginary completion of $\widehat{E}(p, q)$*) The geometry of type $\mathsf{D}_5$ corresponding to the half spin geometry $\widehat{T}(p, q)$ will be denoted by $\Theta(\widehat{T}(p, q))$. By Theorem 5.33 we can assume that it contains $\widehat{E}(p, q)$ as a geometric hyperplane. We call it the *imaginary* completion of $\widehat{E}(p,q)$, and the points of $\Theta(\widehat{T}(p,q))\setminus\widehat{E}(p,q)$ are called the *imaginary points* of $\widehat{E}(p,q)$. We will provide an interpretation of these imaginary points of $\widehat{E}(p,q)$ in Corollary 6.11.

**Corollary 5.35** *If $p', q'$ are two opposite points of $\Gamma$ and $\widehat{T}(p',q') = \widehat{T}(p,q)$, then $p', q' \in \widehat{E}(p,q)$. In other words, $\widehat{E}(p',q') = \widehat{E}(p,q)$.*

*Proof* Let $x \in \widehat{E}(p,q)$. We can find two hyperbolic solids $U_1, U_2$ of $\widehat{E}(p,q)$ such that $U_1 \cap U_2 = \{x\}$. Then by Proposition 5.32(*iii*), we know that $a_1 = \beta(U_1)$ and $a_2 = \beta(U_2)$ are special and $x = a_1 \bowtie a_2$. Again, by Proposition 5.32(*iii*), and the fact that $a_1, a_2 \in \widehat{T}(p',q')$, we also know that $a_1 \bowtie a_2 \in \widehat{E}(p',q')$. Hence $\widehat{E}(p,q) \subseteq \widehat{E}(p',q')$, implying equality. □

### 5.4 The hyperplane $\widehat{H}(p,q)$ of $\Gamma$

We denote by $\widehat{H}(p,q)$ the set of points of $\Gamma$ collinear or equal to at least one point of $\widehat{E}(p,q)$. The notation $\widehat{H}(p,q)$ comes from "hyperplane". We indeed intend to prove that $\widehat{H}(p,q)$ is a geometric hyperplane of $\Gamma$ and that each member of $\widehat{E}(p,q) \cup \widehat{T}(p,q)$ is a *deep point* of it (a *deep point* of a geometric hyperplane $H$ is a point $x$ for which $x^\perp$ is contained in $H$). But first we need another lemma.

Let $x$ be any point of $\Gamma$ and let $\mathcal{N}_x$ denote the set of lines of $\Gamma$ through $x$.

**Lemma 5.36** *Let $x$ be any point of $\Gamma$. Furnished with the planar point pencils, $\mathcal{N}_x$ has the structure of a dual polar space of type B$_3$ over $\mathbb{K}$. If we furnish it further with all subsets of $\mathcal{N}_x$ consisting of all lines intersecting some hyperbolic line contained in $x^\perp$, then $\mathcal{N}_x$ has the structure of a polar space of type D$_4$ over $\mathbb{K}$.*

*Proof* The first assertion follows immediately from the diagram of $\Gamma$ since the geometry under consideration is just a point-line truncation of the residue at $x$. The second assertion follows from the fact that we obtain a half spin D$_4$-geometry from a dual polar space of type B$_3$ if we add the "hyperbolic lines" of all "quads" (using the language of dual polar spaces and near polygons). These "quads" are symplectic quadrangles, and they correspond precisely to the residual geometries at $x$ of the symplecta through $x$. The "hyperbolic lines" are then the sets of lines through $x$ in the symplecta meeting a common hyperbolic line not containing $x$ but contained in a symplecton through $x$. Hence the lemma follows. □

We denote the point-line geometry of type D$_4$ on $\mathcal{N}_x$ by D$_4(\mathcal{N}_x)$.

**Lemma 5.37** (i) *$\widehat{H}(p,q)$ is the union of all lines containing a point of $\widehat{E}(p,q)$ and a point of $\widehat{T}(p,q)$.*
(ii) *$\widehat{H}(p,q)$ is the union of all symplecta containing a point of $\widehat{E}(p,q)$.*
(iii) *$\widehat{H}(p,q)$ is the union of all lines containing a point of $\widehat{T}(p,q)$.*
(iv) *$\widehat{H}(p,q)$ is a proper geometric hyperplane of $\Gamma$ and its set of deep points coincides with $\widehat{E}(p,q) \cup \widehat{T}(p,q)$.*

*Proof* Let $H = \widehat{H}(p,q)$, let $H^-$ be the union of all lines containing a point of both $\widehat{E}(p,q)$ and $\widehat{T}(p,q)$, and let $H^+$ be the union of all symplecta containing a point of $\widehat{E}(p,q)$. Then, obviously, $H^- \subseteq H$ and by Corollary 5.2, we also have $H \subseteq H^+$. In order to show (*i*) and (*ii*), it suffices to prove that $H^+ \subseteq H^-$. So let $x \in H^+$, and let $S$ be a symplecton containing $x$ and a point of $\widehat{E}(p,q)$, which we can choose to be $p$. By Corollary 5.23, $S$ contains a hyperbolic line $h$ contained in $\widehat{E}(p,q)$. We may assume that $x \notin h$. Then $x$ is collinear with at least one point $y \in h$. Consider a second point $y' \in h \setminus \{y\}$. Then $y'$ is collinear with exactly one point $z \in xy$ and so $z \in \widehat{T}(p,q)$, whereas $y \in \widehat{E}(p,q)$ and $x \in yz$, showing $x \in H^-$.

Now let $H^*$ be the set of all points collinear in $\Gamma$ with at least one point of $\widehat{T}(p,q)$. Clearly, by $(i)$, $H \subseteq H^*$. Hence, in order to show $(iii)$, it suffices to prove $H^* \subseteq H$. So let $x \in H^*$ and let $a \in \widehat{T}(p,q)$ be collinear with $x$. Considering $\mathcal{N}_a$, Lemma 5.36 yields a symplecton containing $ax$ and a point of $\beta(a)$. Consequently, $x$ is collinear or symplectic to at least one point of $\widehat{E}(p,q) \supseteq \beta(a)$. By $(ii)$, this suffices to conclude $x \in H$ and $(iii)$ follows.

We now show that $H$ is a proper geometric hyperplane of $\Gamma$. First we prove that it is a subspace. Let $x_1, x_2 \in H$ be collinear points. If one of $x_1, x_2$ belongs to $\widehat{E}(p,q) \cup \widehat{T}(p,q)$, the definition of $H$ and $(iii)$ imply that the line $x_1 x_2$ belongs to $H$. So we may assume that neither of $x_1, x_2$ belongs to $\widehat{E}(p,q) \cup \widehat{T}(p,q)$. But by $(i)$, $x_i$, $i = 1, 2$, belongs to a line $L_i$ intersecting $\widehat{E}(p,q)$ in some point $y_i$ and intersecting $\widehat{T}(p,q)$ in some point $z_i$. If $\beta(z_1) \cap \beta(z_2) = \emptyset$, then, by Proposition 5.32$(iv)$, $z_1$ and $z_2$ are opposite. Since $z_1$ is not opposite $x_2$, it is opposite $y_2$ by Lemma 5.7. This contradicts Corollary 5.27.

Hence we may assume that $\beta(z_1) \cap \beta(z_2)$ contains some point $b$. Let $i \in \{1, 2\}$. Since $b$ is collinear with $z_i$ and also either equal to or symplectic to $y_i$, the pair $\{b, x_i\}$ is either a collinear pair or a symplectic pair, respectively, and vice versa. Note that, if $b \neq y_i$, then both $z_i$ and $x_i$ belong to $b \Diamond y_i$. First suppose that $\{b, x_1\}$ is symplectic. If $x_2$ belongs to $b \Diamond x_1$, then by $(ii)$, all points of $x_1 x_2$ belong to $H$. So we may assume that $x_2$ is close to $b \Diamond x_1$. This implies that $x_2$ is collinear with all points of a line $M_2 \subseteq b \Diamond x_1$. Since $\{b, x_2\}$ is a collinear or symplectic pair, $b$ is either on $M_2$ or collinear with all points of $M_2$. In any case $b$ is collinear with $x_1$, a contradiction. Similarly, $b \perp x_2$ leads to a contradiction. So we may assume that both $x_1$ and $x_2$ are collinear with $b$. In this case, $L_1$ meets $L_2$ in $b = y_1 = y_2$ and every point of the line $x_1 x_2$ is collinear with $b$, which proves that $x_1 x_2 \subseteq H$. Thus, $H$ is a subspace of $\Gamma$.

By the foregoing, all members of $\widehat{E}(p,q) \cup \widehat{T}(p,q)$ are deep points. Now we show that no point of $\widehat{H}(p,q) \setminus (\widehat{E}(p,q) \cup \widehat{T}(p,q))$ is deep. This will also imply that $H$ is proper.

Let $x_1$ be any point of $\widehat{H}(p,q) \setminus (\widehat{E}(p,q) \cup \widehat{T}(p,q))$. By $(i)$, $x_1$ belongs to a line $L_1$ intersecting $\widehat{E}(p,q)$ in some point $y_1$ and intersecting $\widehat{T}(p,q)$ in some point $z_1$. Let $x_2$ be any point collinear with $x_1$ such that $\{x_2, y_1\}$ is a special pair (such a point exists by Fact 5.3). Then Lemma 5.7 implies that $\{x_2, z_1\}$ is also a special pair. By Lemma 5.17, $x_2 \notin \widehat{E}(p,q)$, and by Proposition 5.32$(iii)$, $x_2 \notin \widehat{T}(p,q)$ since $z_1 \bowtie x_2 = x_1 \notin \widehat{E}(p,q)$. Assume, for a contradiction, that $x_2 \in H$. Then $x_2$ belongs to a line $L_2$ intersecting $\widehat{E}(p,q)$ in some point $y_2$ and intersecting $\widehat{T}(p,q)$ in some point $z_2$. Note that $y_1$ and $y_2$ are not opposite since this would imply by Lemma 5.7 that $y_1$ and $z_2$ are opposite (taking into account that $y_1$ and $x_2$ are not opposite), contradicting Corollary 5.27. Hence $y_1 \perp y_2$. Since $x_2 \perp y_2$ and $x_2$ is special to $y_1$, we have $y_1 \neq y_2$ and the point $x_2$ is close to $y_1 \Diamond y_2$. By Fact 5.5, $x_2 \bowtie y_1 \in y_1 \Diamond y_2$, and so, since $x_2$ is collinear with a line of $y_1 \Diamond y_2$, the points $x_2 \bowtie y_1 = x_1$ and $y_2$ are collinear. By Lemma 5.7, also $x_1$ and $z_2$ are collinear. Also, since $y_1 \perp y_2$ and $y_1$ is special to $x_2$, the same lemma implies that $y_1$ is special to $z_2$. Again the same lemma then implies that $z_2$ and $z_1$ are special (since $z_2$ is special to $y_1$ and collinear with $x_1$). But then $z_1 \bowtie z_2 = x_1 \notin \widehat{E}(p,q)$, contradicting Proposition 5.32$(iii)$. Hence $x_2 \notin H$ and so $x_1$ is not a deep point of $H$.

We now show that $H$ is a geometric hyperplane of $\Gamma$. Since we have already proven that $H$ is a subspace, it suffices to show that every line $M$ of $\Gamma$ not inside $H$ intersects $H$ in a point. Let $r \in M$ be a point not belonging to $H$. Since $H = H^+$, any point of $\widehat{E}(p,q)$ is either special to or opposite $r$. Let $H_r$ be the set of points of $\widehat{E}(p,q)$ which are special to $r$. We prove a number of claims (and forget the notation for points already used in the current proof).

- $H_r$ *is a subspace of* $\widehat{E}(p,q)$, *viewed as a polar space.* Indeed, let $\{a,b\}$ be a symplectic pair in $H_r$. There are two possibilities. The first possibility is that $r$ is far from $a \lozenge b$. Let $s$ be the unique point of $a \lozenge b$ that is symplectic to $r$. Then both $a$ and $b$ are collinear with $s$, and, by Lemma 5.6, so is every point of $h(a,b)$. Hence, $h(a,b) \subseteq H_r$ by Fact 5.5 and Lemma 5.17. The second possibility is that $r$ is close to $a \lozenge b$. Then no point of $h(a,b)$ can be opposite $r$; so they are all special to $r$ and again $h(a,b) \subseteq H_r$.
- $H_r$ *is either a geometric hyperplane of* $\widehat{E}(p,q)$, *viewed as a polar space, or coincides with it.* Suppose $H_r$ does not coincide with $\widehat{E}(p,q)$, and let $h$ be a hyperbolic line in $\widehat{E}(p,q)$ containing at least one point opposite $r$. Then, Lemma 5.6 implies that a unique point $x$ of $h$ is not opposite $r$, and hence it is special to $r$ since $r \notin \widehat{H}(p,q)$.

To every point $x$ of $H_r$, we associate the line $L_x$ through $r$ containing $r \bowtie x$. As before, we denote the set of lines of $\Gamma$ through $r$ by $\mathcal{N}_r$, and this mapping by $\lambda : H_r \to \mathcal{N}_r : x \mapsto \lambda(x) = L_x$.

- *The map $\lambda$ just defined is injective.* Indeed, suppose $\lambda(x_1) = \lambda(x_2)$, for two distinct points $x_1, x_2 \in H_r$. If $r \bowtie x_1 \neq r \bowtie x_2$, then, noting both belong to $H$, the line joining $r \bowtie x_1$ and $r \bowtie x_2$ belongs to $H$ as $H$ is a subspace. Consequently $r \in H$, a contradiction. Hence we may assume that $y = r \bowtie x_1 = r \bowtie x_2$. But then $y \in \widehat{T}(p,q)$ and $(iii)$ implies $r \in H$, again the same contradiction.
- *The map $\lambda$ maps opposite pairs of points of $H_r$ onto pairs of lines not contained in a symplecton.* Let $x_1, x_2 \in H_r$ be opposite points and suppose $L_1 = \lambda(x_1)$ and $L_2 = \lambda(x_2)$ are contained in a common symplecton $S$. Both $x_1$ and $x_2$ are close to $S$ and so collinear with respective lines $M_1$ and $M_2$ belonging to $S$. Inside $S$, $M_1$ and $M_2$ are opposite because if a point of $u$ were collinear to all point of $M_2$, then $x_2 \perp\!\!\!\perp u$ and by the last assertion of Lemma 5.7 and the fact that $x_1$ and $x_2$ are opposite, this is a contradiction. Still inside $S$, every point of the hyperbolic line $M_1^\perp \cap M_2^\perp$ is symplectic to both $x_1$ and $x_2$ and hence belongs to $\widehat{E}(p,q)$; hence $r \in S$ belongs to $H$ by $(ii)$, a contradiction.
- *The map $\lambda$ sends the points of any hyperbolic line $h \subseteq H_r$ either to all lines of a planar line pencil through $r$, or to all lines through $r$ intersecting a certain hyperbolic line $g \subseteq r^\perp$.* Set $S = S(h)$. By $(ii)$, $r \notin S$. There are two possibilities. The first possibility is that $r$ is close to $S$. Let $L = r^\perp \cap S$ be the unique line of $S$ collinear with $r$. By the injectivity of $\lambda$, the projection of $h$ onto $L$ in $S$ is injective, and since $h$ is a hyperbolic line in $S$, it is also surjective. So $\lambda(h)$ consists of all lines of the plane through $r$ and $L$, and hence is a (full) planar line pencil.

  The second possibility is that $r$ and $S$ are far. Let $z \in S$ be the unique point symplectic to $r$. Since all points of $h$ are special to $r$, Fact 5.5 implies that $h \subseteq z^\perp$. Set $T = r \lozenge z$. Then, each point $x \in h$ is close to $T$ and so $x \bowtie r$ belongs to $T$ and is collinear with $z$. This implies $\lambda(h) \subseteq T$. Now consider any point $y \in h^\perp$, with $y \neq z$. Then $r$ and $y$ are opposite, and hence, projection (in the sense of Proposition 4.12) defines an isomorphism between $\mathcal{N}_y$ and $\mathcal{N}_r$ (and note, here, that projection means "meeting a common line"). But the line $xy$, for $x \in h$, corresponds to $\lambda(x)$. Since $x$ runs through a hyperbolic line, also $r \bowtie x$ runs through a hyperbolic line and the assertion is proved.

All this now implies that $\lambda(H_r)$ is a fully embedded subgeometry and a subspace of the geometry $\mathsf{D}_4(\mathcal{N}_r)$. However, no subspace of a polar space of type $\mathsf{D}_4$ is isomorphic to a geometric hyperplane of a polar space of type $\mathsf{B}_4$, except for the full polar space itself. Hence $\lambda$ is surjective. This implies that every line $R$ through $r$ contains a point collinear with some point of $H_r$; hence $R$ contains a point of $H$. □

The following are immediate corollaries of the proof of Lemma 5.37.

**Corollary 5.38** *Let $r$ be a point not belonging to $\widehat{H}(p, q)$. Then the set of points of $\widehat{E}(p, q)$ not opposite $r$ (hence special to $r$) induces a polar subspace of type $\mathsf{D}_4$ in $\widehat{E}(p, q)$, viewed as a full polar space of type $\mathsf{B}_4$ over $\mathbb{K}$.* □

**Corollary 5.39** *Let $L$ be any line of $\Gamma$ containing a point of $\widehat{E}(p, q)$. Then $L$ contains precisely one point of $\widehat{T}(p, q)$.* □

## 6 Constructing a building of type $\mathsf{E}_6$ from a split building of type $\mathsf{F}_4$

In this section we conclude our construction. The point set of the building of type $\mathsf{E}_6$ is the union of the point set of $\Gamma$ and the family of all extended equator geometries of $\Gamma$. In order to well-define the lines, we need another lemma, which we prove in Sect. 6.1. In the rest of the section, we identify suitable substructures of the split building $\Gamma$ of type $\mathsf{F}_4$ as elements of various types of the building $\Delta$ of type $\mathsf{E}_6$, and we define a suitable incidence relation among them to conclude the construction of $\Delta$. Further, we construct the symplectic polarity of $\Delta$ whose fixed point geometry is precisely $\Gamma$, and prove Theorems 1 and 2.

### 6.1 The point-line $\mathsf{E}_6$-geometry

The following proposition is the basis for the definition of new lines.

**Proposition 6.1** *Let $\widehat{E} = \widehat{E}(p, q)$ be an extended equator geometry, with $p, q$ two opposite points of $\Gamma$. Let $x$ be a point of $\widehat{T}(p, q)$ collinear with $p$ and put $U = \beta(x)$. Let $y$ be the point of $\Gamma$ collinear with all points of the hyperbolic solid $V$ of $\widehat{E}$ containing $q$ and the hyperbolic plane $U \cap E(p, q)$. Then an extended equator geometry $\widehat{E}'$ contains $U$ if, and only if, it can be written as $\widehat{E}(p, q')$, with $q' \in qy\setminus\{y\}$. Also, if $q' \neq q$, then $\widehat{E} \cap \widehat{E}' = U$. Consequently, if $q', q'' \in qy\setminus\{y\}$ with $q' \neq q''$, then $\widehat{E}(p, q') \cap \widehat{E}(p, q'') = U$.*

*Proof* Let $q'$ be a point of $qy\setminus\{y\}$. Since $p \perp x \perp y$ (see Proposition 5.32(i)), we have that $p$ is not opposite $y$. Since $p$ is opposite $q$, Lemma 5.7 implies that $p$ is opposite $q'$. Since each point $a$ of $U \cap E(p, q)$ is collinear with $y$ and symplectic to $q$, it is symplectic to $q'$, as follows from Lemma 5.7. It follows that $a \in E(p, q')$. Hence, $\widehat{E}(p, q')$ contains $U \cap E(p, q)$. Since it also contains $p$, we easily deduce, using Proposition 5.32(ii), that $U \subseteq \widehat{E}(p, q')$

Let $\widehat{E}'$ be an extended equator geometry containing $U$. Let $a, b, c$ be three points in $E(p, q) \cap U = V \cap U$ not on a common hyperbolic line. Then $a, b, c, p$ generate $U$ inside $\widehat{E}(p, q)$. Let $S_a, S_b, S_c$ be the symplecta containing $q$ and $a, b, c$, respectively. Since $a, b, c \in \widehat{E}'$, each of $S_a, S_b, S_c$ contains, by Corollary 5.23, a hyperbolic line $h_a, h_b, h_c$, respectively, entirely contained in $\widehat{E}'$. Also, since $y$ is collinear with $q, a, b, c$, it belongs to $S_a, S_b, S_c$. By the isomorphism $\varphi_{q,p}$ (Proposition 5.9), the intersection of $S_a, S_b, S_c$ is precisely the line $L := qy$. Since hyperbolic lines of symplecta are geometric lines, the point $q$ is collinear in $\Gamma$ with unique points $q_a, q_b, q_c$ of $h_a, h_b, h_c$, respectively. Suppose, by way of contradiction, that $q_a \neq q_b$. Then $q$ belongs to the tropic circle geometry $\widehat{T}'$ of $\widehat{E}'$. But also $y$ belongs to $\widehat{T}'$. Since $y \perp q$, Proposition 5.32(i) implies that $y^\perp \cap q^\perp \cap \widehat{E}'$ is a hyperbolic plane $\pi$. Clearly, $\pi$ intersects $U \cap V$ in a hyperbolic line, contradicting the fact that $q$ is not collinear with any point of $U$. Hence $q_a = q_b = q_c =: q' \in L$, with $q' \neq y$. It follows from Lemma 5.7(7) that $q'$ is opposite $p$ and so $\widehat{E}' = \widehat{E}(p, q')$ by Proposition 5.22. The same proposition also implies the last assertion, as $U$ is a maximal singular subspace of both $\widehat{E}$ and $\widehat{E}'$ and, by Proposition 5.22, $U \subsetneq \widehat{E} \cap \widehat{E}'$ would mean that $\widehat{E} = \widehat{E}'$. □

**Definition 6.2** (*Point-Line* E$_6$-*Geometry*) We define the *point-line* E$_6$-*geometry* as the pair $(\mathscr{P}, \mathscr{L})$, where $\mathscr{P}$ is the point set of $\Gamma$ union the family $\mathscr{E}$ of extended equator geometries, and $\mathscr{L}$ is the set of ordinary and hyperbolic lines of $\Gamma$ union the following family $\mathscr{F}$ of subsets of $\mathscr{P}$. Let $U$ be any hyperbolic solid of $\Gamma$. Then $\beta(U)$ together with all extended equator geometries containing $U$ is a general element of $\mathscr{F}$ (see Proposition 6.1). Inclusion between the elements of $\mathscr{P}$ and those of $\mathscr{L}$ defines incidence. Members of $\mathscr{E}$ and $\mathscr{F}$ will frequently be referred to as *the new points* and *the new lines*, respectively.

Note that $(\mathscr{P}, \mathscr{L})$ is a *partial linear space*, i.e., every pair of distinct points is contained in at most one line. When we write about collinear points of $\Gamma$, we will always mean the collinearity in $\Gamma$, and not in $(\mathscr{P}, \mathscr{L})$, unless explicitly mentioned otherwise.

Note also that, by Lemma 5.17 and Proposition 5.32($ii$), the sets $\widehat{T}(p, q)$ and $\widehat{E}(p, q)$ are subspaces of $(\mathscr{P}, \mathscr{L})$.

## 6.2 $(\mathscr{P}, \mathscr{L})$ corresponds to a building of type E$_6$: First observations

We are left to show that $(\mathscr{P}, \mathscr{L})$ is the point-line geometry of the building of type E$_6$ over the field $\mathbb{K}$. Towards that end, we continue our series of lemmas.

A *full pencil* in a tropic circle geometry $\widehat{T}$ is the intersection of $\widehat{T}$ with the union $x^{\perp} \cup x^{\perp\!\!\!\perp}$ for a certain point $x$ of $\widehat{T}$ (hence it is the set of points of $\widehat{T}$ collinear with or equal to $x$, with "collinearity" in the half spin geometry of type D$_5$ defined by $\widehat{T}$, see Theorem 5.33($ii$); this collinearity coincides precisely with collinearity in $(\mathscr{P}, \mathscr{L})$). The point $x$ is then called the *centre* of the full pencil. Obviously, a full pencil has a unique centre. Also, for a new point $\mathfrak{e} \in \mathscr{E}$, we denote by $T_{\mathfrak{e}}$ the corresponding tropic circle geometry. The element of $\mathscr{L}$ through distinct elements $\alpha$ and $\alpha'$ in $\mathscr{P}$ will be denoted by $\langle \alpha, \alpha' \rangle$. Note that the definition of $\mathscr{F}$ yields a natural bijective correspondence between $\mathscr{F}$ and the family of hyperbolic solids of $\Gamma$.

**Lemma 6.3** *Let $\mathfrak{e}, \mathfrak{e}' \in \mathscr{E}$. Then, $\mathfrak{e}$ and $\mathfrak{e}'$ are collinear in $(\mathscr{P}, \mathscr{L})$ if, and only if, $T_{\mathfrak{e}} \cap T_{\mathfrak{e}'}$ contains a full pencil in both $T_{\mathfrak{e}}$ and $T_{\mathfrak{e}'}$. In this case, the intersection is that full pencil and it has centre $\langle \mathfrak{e}, \mathfrak{e}' \rangle \cap \Gamma$.*

*Proof* First assume that $\mathfrak{e}$ and $\mathfrak{e}'$ are collinear in $(\mathscr{P}, \mathscr{L})$ and let $U = \mathfrak{e} \cap \mathfrak{e}'$. Let $x \in T_{\mathfrak{e}}$ be collinear with all points of $U$. Then clearly $x \in T_{\mathfrak{e}'}$. Now suppose $y \in T_{\mathfrak{e}}$ is collinear or symplectic to $x$. Then, by Proposition 5.32($i$) and ($ii$), $y$ is collinear with at least two points of $U$, and hence $y$ also belongs to $T_{\mathfrak{e}'}$. Now suppose that some point $z \in T_{\mathfrak{e}}$ also belongs to $T_{\mathfrak{e}'}$, with $z$ special to or opposite $x$. Let $V = z^{\perp} \cap \mathfrak{e}$. Then, by Lemma 5.28, $V \neq U$ and, by Proposition 5.18, there is a hyperbolic solid $W$ of $\mathfrak{e}$ intersecting $V$ in just a point $z' \notin U$ and intersecting $U$ in at least two points (and hence in a hyperbolic line or in a hyperbolic plane). The latter implies that the point $z''$ of $T_{\mathfrak{e}}$ corresponding to $W$ belongs to $T_{\mathfrak{e}} \cap T_{\mathfrak{e}'}$. Further, Proposition 5.32($iii$) implies that $\{z'', z\}$ is a special pair. Hence $z' = z \bowtie z''$ belongs to both $\mathfrak{e}$ and $\mathfrak{e}'$, contradicting our assumption $z \notin U$.

Now suppose that $T_{\mathfrak{e}} \cap T_{\mathfrak{e}'}$ contains a full pencil in both $T_{\mathfrak{e}}$ and $T_{\mathfrak{e}'}$. Let $x \in T_{\mathfrak{e}}$ be the centre of such a full pencil in $T_{\mathfrak{e}}$. Recall that, for $z \in T_{\mathfrak{e}}$, we denote $\beta(z) = z^{\perp} \cap \mathfrak{e}$. Given any point $u \in \beta(x)$, we can choose $y, z \in T_{\mathfrak{e}}$ such that $\beta(y) \cap \beta(x)$ and $\beta(z) \cap \beta(x)$ are hyperbolic lines, whereas $\beta(y) \cap \beta(z) = \{u\}$. Consequently, $y$ and $z$ are both symplectic to $x$ and special to each other. As $u = y \bowtie z$ and $y, z \in T_{\mathfrak{e}'}$, Proposition 5.32($iii$) implies that $u \in \mathfrak{e}'$. Since $u \in \beta(x)$ was arbitrary, $\beta(x) \subseteq \mathfrak{e}'$ and by definition, this yields that $\mathfrak{e}'$ is collinear with $\mathfrak{e}$ in $(\mathscr{P}, \mathscr{L})$. It follows from the first paragraph that $T_{\mathfrak{e}} \cap T_{\mathfrak{e}'}$ coincides with the full pencil in both $T_{\mathfrak{e}}$ and $T_{\mathfrak{e}'}$ with centre $\langle \mathfrak{e}, \mathfrak{e}' \rangle \cap \Gamma$. □

We can be a little more specific about the geometric structure of a full pencil in a tropic circle geometry, and its relation to the corresponding extended equator geometry.

**Lemma 6.4** *Let $U$ be a hyperbolic solid in $\widehat{E}(p,q)$ and let $P_U$ be the set of points of $\Gamma$ collinear with at least two points of $U$. Then,*

(i) *$P_U \subseteq \widehat{T}(p,q)$ and the structure on it induced by $(\mathscr{P}, \mathscr{L})$ is a cone with vertex $\beta(U)$ over a geometry isomorphic to the line Grassmannian of a projective 4-space $W$ over $\mathbb{K}$;*
(ii) *The subgeometry $P_U^\Gamma$ of $P_U$ restricted to the ordinary lines of $\Gamma$ through $\beta(U)$ is isomorphic to a cone over a 3-space over $\mathbb{K}$; in the line Grassmannian, this 3-space corresponds to all lines through a point of $W$. Hence, $P_U^\Gamma$ is a projective 4-space over $\mathbb{K}$.*

*Proof* (i) By the definition of $\widehat{T}(p,q)$, every point of $P_U$ belongs to $\widehat{T}(p,q)$. Since $P_U$ is the set of points $x$ of $\widehat{T}(p,q)$ such that $x^\perp \cap U$ contains a hyperbolic line, Lemmas 5.28 and 5.29, and Proposition 5.32 imply that $P_U$ is the full pencil in $\widehat{T}(p,q)$ with deep point $\beta(U)$. By Theorem 5.33(ii), $\widehat{T}(p,q)$, endowed with its ordinary and hyperbolic lines, is a half spin geometry of type $\mathsf{D}_5$, say with corresponding polar space $\Omega$ of type $\mathsf{D}_5$ and system $\Phi^+$ of maximal singular subspaces. So $\beta(U)$ corresponds to an element $U^+ \in \Phi^+$. Then the set $P_U \setminus \{\beta(U)\}$ corresponds to the set of elements of $\Phi^+$ intersecting $U^+ \in \Phi^+$ in planes. Clearly, elements of $\Phi^+$ intersecting $U^+$ in the same plane correspond to points in a(n ordinary or hyperbolic) line of $\Gamma$ through $\beta(U)$, and vice versa. Hence the lines (of $\widehat{T}(p,q)$ viewed as a half spin geometry of type $\mathsf{D}_5$) through $\beta(U)$ correspond bijectively and naturally to the points of the plane Grassmannian $G_2(U^+)$ of $U^+$. We show that this bijective correspondence is an isomorphism. Note that, by Lemmas 5.6 and 5.7, two lines $L$ and $L'$ through $\beta(U)$ in $\widehat{T}(p,q)$ are "collinear" (meaning that each point of $L \setminus \{\beta(U)\}$ is either collinear or symplectic to each point of $L' \setminus \{\beta(U)\}$) if, and only if, there is a point $y \in L \setminus \{\beta(U)\}$ collinear or symplectic to a point $y' \in L' \setminus \{\beta(U)\}$. Two planes $\pi_1$ and $\pi_2$ are collinear in $G_2(U^+)$ precisely if they intersect in a line. This happens if, and only if, some members $W_1, W_2 \in \Phi^+$, with $\pi_i \subseteq W_i$, $i = 1, 2$, intersect in a plane, hence, if, and only if, the corresponding lines are "collinear". Note that $G_2(U^+)$ is equal to the line Grassmannian of the dual space $W$ of $U^+$.

(ii) The ordinary lines through $\beta(U)$ in $P_U$ correspond, according to Proposition 5.32(i), to the members of $\Phi^+$ intersecting $U^+$ in planes contained in $U$. Hence we obtain a 3-space in the plane Grassmannian of $U^+$, implying that $P_U^\Gamma$ is a cone with vertex $\beta(x)$ over a 3-space. Dualising, we see that the ordinary lines of $\Gamma$ through $\beta(U)$ correspond to all lines through a point of $W$. □

Let $\widehat{T}$ be a tropic circle geometry and $\widehat{E}$ its associated extended equator geometry. For a point $x \in \widehat{E}$, the set of points $\widehat{T} \cap x^\perp$, endowed with the lines and hyperbolic lines contained in it, is a subspace of $\widehat{T}$ isomorphic to a polar space of type $\mathsf{D}_4$, as follows from Theorem 5.33(ii). Such a polar space will be referred to as a *standard $\mathsf{D}_4$ in $\widehat{T}$*. A geometric hyperplane of $\widehat{T} \cap x^\perp$ isomorphic to a polar space of type $\mathsf{B}_3$ will be referred to as a *standard $\mathsf{B}_3$ in $\widehat{T}$ with centre $x$*. Any singular geometric hyperplane of $\widehat{T} \cap x^\perp$ consists of the set $H_a$ of points of $\widehat{T} \cap x^\perp$ collinear, symplectic or equal to a given point $a$ of $\widehat{T} \cap x^\perp$, which is then a deep point of $H_a$. It is unique since for every point $b \in H_a \setminus \{a\}$, we can find a point $c \in H_a$ special to $b$. Hence $b$ cannot be a deep point as well.

For the next lemma, we recall that the lines of the tropic circle geometries are the ordinary and hyperbolic lines of $\Gamma$ contained in it. As such, a geometric hyperplane of a subspace $A$ of a tropic circle geometry is defined relative to the ordinary and hyperbolic lines contained in $A$.

**Lemma 6.5** *Let $p$ be any point of $\Gamma$ and let $p \in \mathfrak{e} \cap \mathfrak{e}'$ with $\mathfrak{e}, \mathfrak{e}' \in \mathscr{E}$. Then,*

(i) *$T_\mathfrak{e} \cap T_{\mathfrak{e}'} \cap p^\perp$ is a geometric hyperplane of both $T_\mathfrak{e} \cap p^\perp$ and $T_{\mathfrak{e}'} \cap p^\perp$;*
(ii) *$\mathfrak{e}$ and $\mathfrak{e}'$ are collinear in $(\mathscr{P}, \mathscr{L})$ if, and only if, $T_\mathfrak{e} \cap T_{\mathfrak{e}'} \cap p^\perp$ is a singular geometric hyperplane in $T_\mathfrak{e} \cap p^\perp$ and $T_{\mathfrak{e}'} \cap p^\perp$. If this is the case, then the (unique) deep point is the unique point of $T_\mathfrak{e} \cap T_{\mathfrak{e}'}$ collinear in $\Gamma$ with every point of the hyperbolic solid $\mathfrak{e} \cap \mathfrak{e}'$;*
(iii) *$\mathfrak{e}$ and $\mathfrak{e}'$ are not collinear in $(\mathscr{P}, \mathscr{L})$ if, and only if, $T_\mathfrak{e} \cap T_{\mathfrak{e}'} \cap p^\perp$ is a standard B$_3$. If this is the case, then $\mathfrak{e} \cap \mathfrak{e}' = \{p\}$.*

*Proof* (i) Let $L$ be an ordinary or hyperbolic line in $T_\mathfrak{e} \cap p^\perp$. Then the structure induced by the lines and hyperbolic lines on the set of points obtained by joining $p$ with the points of $L$ is a projective plane $\pi$ (noting that, if $L$ is a hyperbolic line, then $p \in S(L)$). By Corollary 5.39, each line through $p$ in $\pi$ contains a unique point of $T_{\mathfrak{e}'}$, and Proposition 5.32$(i)$ and $(ii)$ imply that the set of points thus obtained is closed under taking lines and hyperbolic lines. It follows that $T_{\mathfrak{e}'} \cap \pi$ is a line $L'$ of $\pi$. Hence either $L = L'$ or $L \cap L'$ is a unique point, as required.

(ii) If $\mathfrak{e}$ and $\mathfrak{e}'$ are collinear in $(\mathscr{P}, \mathscr{L})$, then Lemma 6.3 implies that $T_\mathfrak{e} \cap T_{\mathfrak{e}'} \cap p^\perp$ is a singular geometric hyperplane in both $T_\mathfrak{e} \cap p^\perp$ and $T_{\mathfrak{e}'} \cap p^\perp$.

Now suppose that $T_\mathfrak{e} \cap T_{\mathfrak{e}'} \cap p^\perp$ is a singular geometric hyperplane in both $T_\mathfrak{e} \cap p^\perp$ and $T_{\mathfrak{e}'} \cap p^\perp$, and let $x$ be the deep point. Let $p' \neq p$ belong to $x^\perp \cap \mathfrak{e}$. Let $L$ be an arbitrary line in $p \Diamond p'$ through $p$. Then the point $x_L := p'^\perp \cap L$ belongs to $T_\mathfrak{e}$. Since also $x$ belongs to $p \Diamond p'$,

- The assumption that $T_\mathfrak{e} \cap T_{\mathfrak{e}'} \cap p^\perp$ is a singular geometric hyperplane in $T_{\mathfrak{e}'} \cap p^\perp$,
- The fact that $x$ is a deep point, and
- $x \perp\!\!\!\perp x_L$,

imply that $x_L \in T_{\mathfrak{e}'}$. Hence $p^\perp \cap p'^\perp$ (which lies automatically in $p \Diamond p'$) is contained in $T_{\mathfrak{e}'}$. Since $p \in \mathfrak{e}'$, Corollary 5.23 implies that some hyperbolic line $h'$ in $p \Diamond p'$ through $p$ is contained in $\mathfrak{e}'$, and hence $h'^\perp$ belongs to $T_{\mathfrak{e}'}$. If $h' \neq h(p, p')$, then $h(p, p')^\perp$ and $h'^\perp$ generate in $(\mathscr{P}, \mathscr{L})$ the 4-space $p^\perp \cap (p \Diamond p')$. Since $T_{\mathfrak{e}'}$ is a subspace, this would imply that $p \in T_{\mathfrak{e}'} \cap \mathfrak{e}'$, a contradiction. Consequently $h(p, p') \subseteq \mathfrak{e}'$ and hence $\mathfrak{e}$ and $\mathfrak{e}'$ intersect in the 3-space $x^\perp \cap \mathfrak{e}$.

(iii) The only nonsingular geometric hyperplanes of a polar space of type D$_4$ are polar spaces of type B$_3$ obtained by slicing with a hyperplane in the standard embedding. Hence, by $(ii)$, $\mathfrak{e}$ and $\mathfrak{e}'$ are not collinear in $(\mathscr{P}, \mathscr{L})$ if, and only if, $T_\mathfrak{e} \cap T_{\mathfrak{e}'} \cap p^\perp$ is a standard B$_3$ in both $T_\mathfrak{e}$ and $T_{\mathfrak{e}'}$.

Finally, assume that $\mathfrak{e} \cap \mathfrak{e}'$ contains two points $p$ and $p'$ with $p \neq p'$. Then, $p^\perp \cap p'^\perp$ is contained in $T_\mathfrak{e} \cap T_{\mathfrak{e}'} \cap p^\perp$ and is a 3-dimensional subspace in $p \Diamond p'$. Hence the geometric hyperplane $T_\mathfrak{e} \cap T_{\mathfrak{e}'} \cap p^\perp$ of $T_\mathfrak{e} \cap p^\perp$ is singular as it cannot be a standard B$_3$. Hence, $\mathfrak{e}$ and $\mathfrak{e}'$ are collinear in $(\mathscr{P}, \mathscr{L})$ and therefore they meet in a hyperbolic solid. □

## 6.3 The quads of $(\mathscr{P}, \mathscr{L})$

**Proposition 6.6** *Let $p$ be any point of $\Gamma$. Put $X = p^\perp \cup \{\widehat{E}(p, q) : q \text{ is opposite } p\}$. Then $X$, endowed with the members of $\mathscr{L}$ it contains as its lines, is a polar space of type D$_5$ over $\mathbb{K}$, denoted by $\Sigma(p)$.*

*Proof* By Lemma 5.36, the structure on the set of lines of $\Gamma$ through $p$ induced by the ordinary and hyperbolic lines is a polar space of type D$_4$. Since $p$ is not collinear with any element of $X \cap \mathcal{E}$, we see that, if $X$ is a polar space, then it is of type $\mathsf{D}_5$. In fact, the verification of all the axioms of a polar space of finite rank is immediate, except for the one-or-all axiom, which is exactly what we proceed to do now. There are a few cases to consider.

**Case 1** *Let $x$ be a point of $\Gamma$ in $p^\perp$ and $L$ an ordinary or a hyperbolic line contained in $p^\perp$.* In this case, the one-or-all axiom follows straight from the fact that the lines of $\Gamma$ through $p$ form a polar space of type $\mathsf{D}_4$ when two lines are considered collinear if they are contained in a common symplecton.

**Case 2** *Let $\mathfrak{e}$ be an extended equator geometry containing $p$ and let $L$ be an ordinary or a hyperbolic line contained in $p^\perp$.* If $L$ contains $p$, then Corollary 5.39 implies that $\mathfrak{e}$ is collinear with the unique point of $L$ in $T_\mathfrak{e}$. If $L$ does not contain $p$, then the structure induced by the lines and hyperbolic lines on the set of points obtained by joining $p$ to each point of $L$ is a projective plane $\pi$. By Corollary 5.39, all lines of $\pi$ through $p$ contain a unique point of $T_\mathfrak{e}$. Since $T_\mathfrak{e}$ is closed under taking lines and hyperbolic lines by Proposition 5.32(*i*) and (*ii*), respectively, these points are on a line $L'$ of $\pi$. As $L'$ intersects $L$ in a unique point or coincides with it, $L$ contains at least one point collinear in $(\mathcal{P}, \mathcal{L})$ to $\mathfrak{e}$.

**Case 3** *Let $x \in p^\perp$ be a point of $\Gamma$ and let $\mathfrak{L} \in \mathscr{L}$ be a new line contained in $X$.* Then, $p$ belongs to the hyperbolic solid $U$ contained in each new point of $\mathfrak{L}$. Let $s$ be the unique point of $\mathfrak{L}$ in $\Gamma$. We may assume that $x \notin \mathfrak{L}$, i.e., $x \neq s$. So suppose first that $x$ and $s$ are collinear or symplectic in $\Gamma$ (note that both belong to $p^\perp$ and hence cannot be opposite each other). If $x$ belongs to $T_\mathfrak{e}$ for some $\mathfrak{e} \in \mathfrak{L}\setminus\{s\}$, then Proposition 5.32(*i*) and (*ii*) imply that $x$ is collinear with at least two points of $U$ and hence it belongs to $T_\mathfrak{e}$ for each $\mathfrak{e} \in \mathfrak{L}\setminus\{s\}$. Consequently, $x$ is collinear in $(\mathcal{P}, \mathcal{L})$ with either only the unique point $s$ of $\mathfrak{L}$ in $\Gamma$, or to all points of $\mathfrak{L}$ according to whether $x$ does not belong to $T_\mathfrak{e}$ for each $\mathfrak{e} \in \mathfrak{L}\setminus\{s\}$ or does belong to $T_\mathfrak{e}$ for some $\mathfrak{e} \in \mathfrak{L}\setminus\{s\}$.

Now suppose that $x$ and $s$ are special. Let $\mathfrak{e}$ be an arbitrary extended equator geometry of $\mathfrak{L}$. By Corollary 5.39, there is a unique point $x_\mathfrak{e}$ on the line $M = px$ which is contained in $T_\mathfrak{e}$. Then $x_\mathfrak{e}^\perp \cap \mathfrak{e}$ is a 3-space $U_{x,\mathfrak{e}}$ sharing only $p$ with $U$. Let $y$ be a point of $T_\mathfrak{e}$ such that $y^\perp \cap U$ is a plane $\pi$ not containing $p$. Put $L = qy$ for each point $q \in (y^\perp \cap \mathfrak{e})\setminus\pi$ (note that $q$ is automatically opposite $p$). By Proposition 6.1, there is a bijective correspondence between the points of $L\setminus\{y\}$ and the extended equator geometries in $\mathfrak{L}$, given by $q' \mapsto \widehat{E}(p, q')$. It also follows that $p$ is opposite all points of $L\setminus\{y\}$ and is special to $y$. Now, since $y^\perp \cap U_{x,\mathfrak{e}}$ is necessarily empty, $y$ and $x_\mathfrak{e}$ are opposite by Proposition 5.32(*iv*). So Lemma 5.7 implies that $y$ (which is special to $p$) is opposite all points of $M\setminus\{p\}$. In addition, this lemma also implies that being not opposite defines a bijection, say $\sigma$, between the points of $M$ and those of $L$. Now note that, for any $q' \in L\setminus\{y\}$, Corollary 5.39 implies that there is a unique point $z$ of $M$ that belongs to $\widehat{T}(p, q')$. Since $q'$ is symplectic to at least one point of $z^\perp \cap \widehat{E}(p, q')$, Lemma 5.7 says that $z$ cannot be opposite $q'$. Hence $\sigma(z) = q'$. Let $q^* = \sigma(x)$. Then clearly, $\widehat{E}(p, q^*)$ is the unique point of $\mathfrak{L}$ collinear with $x$. This completes the proof in this case.

**Case 4** *Let $\mathfrak{e}$ be an extended equator geometry containing $p$ and let $\mathfrak{L} \in \mathscr{L}$ be a new line contained in $X$.* Note that all extended equator geometries belonging to $\mathfrak{L}$ contain $p$. We may assume that $\mathfrak{e} \notin \mathfrak{L}$. Let $s$ be the unique point of $\mathfrak{L}$ in $\Gamma$ and let $x$ be the unique point of $T_\mathfrak{e}$ on the line $sp$. Let $H_s = T_\mathfrak{e} \cap (s^\perp \cup s^\perp\!\!\!\perp) \cap p^\perp$ and, for a point $\mathfrak{f}$ of $\mathfrak{L}\setminus\{s\}$, put $H_\mathfrak{f} = T_\mathfrak{e} \cap T_\mathfrak{f} \cap p^\perp$.

First note that $T_\mathfrak{e} \cap p^\perp$ only has points of $\Gamma$ and the lines of $(\mathcal{P}, \mathcal{L})$ in it are either ordinary or hyperbolic lines of $\Gamma$. It then follows from cases 1 and 2, respectively, that the sets $H_s$ and $H_\mathfrak{f}$, for every $\mathfrak{f}$ in $\mathfrak{L}\setminus\{s\}$, are geometric hyperplanes of $T_\mathfrak{e} \cap p^\perp$. Moreover, case 3 implies that every point of $T_\mathfrak{e} \cap p^\perp$ is contained in either one or each of those hyperplanes. Let $H$ denote the intersection of all these hyperplanes. By the previous cases, $H$ is a geometric hyperplane of each of them, in particular of $H_s$, which is equal to $T_\mathfrak{e} \cap (x^\perp \cup x^\perp\!\!\!\perp) \cap p^\perp$ by Lemma 5.7.

Note that $H_s$ is a cone over a polar space $D$ of type $\mathsf{D}_3$ and recall that Lemma 6.5 states that a point $\mathfrak{f} \in \mathfrak{L}\setminus\{s\}$ is not collinear with $\mathfrak{e}$ precisely if $H_\mathfrak{f}$ is a polar space of type $\mathsf{B}_3$.

Suppose first $x = s$. Then clearly $H$ contains $x$. If $H \cap D$ is a singular geometric hyperplane of $D$, then $H$ is a cone over a polar space of type $\mathsf{D}_2$ (with a line as vertex). If $H \cap D$ is nonsingular, then $H$ is a cone over a polar space of type $\mathsf{B}_2$ (with a point as vertex). In the former case, $H_\mathfrak{f}$ cannot be of type $\mathsf{B}_3$ for any $\mathfrak{f} \in \mathfrak{L}\setminus\{s\}$, as it contains 3-spaces. In the latter case, $H_\mathfrak{f}$ is of type $\mathsf{B}_3$ for all $\mathfrak{f} \in \mathfrak{L}\setminus\{s\}$, as otherwise it would have to be isomorphic to a cone (with a line containing $x$ as vertex) over a polar space of type $\mathsf{B}_2$, a contradiction.

Now suppose $x \neq s$. As $H$ does not contain $x$, it is a polar space of type $\mathsf{D}_3$. Then $T_\mathfrak{e} \cap p^\perp$, being a polar space of type $\mathsf{D}_4$, contains exactly two points collinear with all points of $H$. Obviously, $x$ is one of them. The other point is contained in $H_\mathfrak{f}$ for a unique point $\mathfrak{f} \in \mathfrak{L}\setminus\{s\}$, which is then the unique point of $\mathfrak{L}$ collinear with $\mathfrak{e}$. □

**Lemma 6.7** *Let $x, y$ be two opposite points of $\widehat{T}(p, q)$. Then $\widehat{E}(x, y) \cap \widehat{T}(p, q)$, endowed with the hyperbolic lines it contains, is a polar space of type $\mathsf{D}_4$ over $\mathbb{K}$ and a geometric hyperplane of $\widehat{E}(x, y)$.*

*Proof* Set $\beta(x) = U \subseteq \widehat{E}(p, q)$ and $\beta(y) = V \subseteq \widehat{E}(p, q)$. The points of $E(x, y) \cap \widehat{T}(p, q)$ correspond bijectively (under $\beta$) to maximal singular subspaces of $\widehat{E}(p, q)$ (viewed as a polar space) intersecting both $U$ and $V$ in a line. Let $\pi$ be a (hyperbolic) plane in $U$, and consider a point $a \in \pi$. Let $b$ be the unique point of $V$ symplectic to all points of $\pi$, and let $\rho$ be the unique (hyperbolic) plane of $V$ all of whose points are symplectic to $a$. Then for any (hyperbolic) line $L \subseteq U$ through $a$ inside $\pi$, the unique (hyperbolic) line $M \subseteq V$ all of whose points are symplectic with all points of $L$ belongs to $\rho$ and contains $b$. Hence, the maximal singular (hyperbolic) subspaces spanned by $L$ and $M$, as $L$ ranges over the set of hyperbolic lines through $a$ inside $\pi$, range over the set of maximal singular (hyperbolic) subspaces of $\widehat{E}(p, q)$ through the hyperbolic line $h(a, b)$ intersecting each of the two opposite maximal singular (hyperbolic) subspaces $U$ and $V$ in a line. By Theorem 5.33(ii), this set is thus exactly a(n imaginary) line of the half spin geometry of type $\mathsf{D}_5$ corresponding to $\widehat{E}(p, q)$, and so the set of points of $\widehat{T}(p, q)$ corresponding to it under $\beta$ is a hyperbolic line in $E(x, y)$. So every line of $U$ defines a point of $E(x, y) \cap \widehat{T}(p, q)$ and line pencils in $U$ correspond to hyperbolic lines in $E(x, y) \cap \widehat{T}(p, q)$, and these correspondences are bijective. It follows that $E(x, y) \cap \widehat{T}(p, q)$ is a polar subspace of type $\mathsf{D}_3$ fully embedded in $E(x, y)$, which is itself a polar space of type $\mathsf{B}_3$.

Since $\widehat{E}(x, y) \cap \widehat{T}(p, q)$ is a subspace of $\widehat{E}(x, y)$ by Lemma 5.17 and Proposition 5.32(ii), it is a polar space which contains a (point) residue of type $\mathsf{D}_3$ (namely, the cone with vertex $x$ and base $E(x, y) \cap \widehat{T}(p, q)$). Suppose that the polar space were degenerate, and let $a$ be a point of it collinear to all points of it. In particular, $x \perp\!\!\!\perp a \perp\!\!\!\perp y$, so $a \in E(x, y)$, implying that $E(x, y) \cap \widehat{T}(p, q)$ is also degenerate, a contradiction. Hence $\widehat{E}(x, y) \cap \widehat{T}(p, q)$ is a polar subspace of $\widehat{E}(x, y)$ of type $\mathsf{D}_4$. □

There is an immediate corollary.

**Corollary 6.8** *Let $x, y$ be two opposite points of $\widehat{T}(p, q)$. Then $x^{\perp\!\!\!\perp} \cap \widehat{T}(p, q) \cap \widehat{E}(x, y)$ has the structure of a cone with vertex $x$ over a polar space of type $\mathsf{D}_3$ over $\mathbb{K}$.* □

Let $\widehat{E}$ be the extended equator geometry associated with $\widehat{T}$. Recall that, for a point $p \in \widehat{E}$, the set $\widehat{T} \cap p^\perp$, endowed with the lines and the hyperbolic lines contained in it, is a *standard* $\mathsf{D}_4$ *in* $\widehat{T}$. A polar space of type $\mathsf{D}_4$ obtained by intersecting $\widehat{T}$ with $\widehat{E}(x, y)$, with $x, y \in \widehat{T}$ opposite, will be called a *hyperbolic* $\mathsf{D}_4$ *in* $\widehat{T}$. A geometric hyperplane therein isomorphic to a polar space of type $\mathsf{B}_3$ over $\mathbb{K}$ will be referred to as a *hyperbolic* $\mathsf{B}_3$ *in* $\widehat{T}$. A hyperplane therein arising as $x^{\perp\!\!\!\perp} \cap \widehat{T}(p,q) \cap \widehat{E}(x,y)$ is called a *hyperbolic* $\mathsf{D}_3$-*cone with vertex* $x$.

The next lemma is probably well-known to the specialists. It will enable us to identify the point set of the imaginary completion of $\widehat{E}(p,q)$ with the set of standard and hyperbolic $\mathsf{D}_4$s of $\widehat{T}(p,q)$.

**Lemma 6.9** *Let $\Theta$ be a half spin geometry of type $\mathsf{D}_5$ and let $\Theta^*$ be the corresponding polar space of type $\mathsf{D}_5$. Then the set of points of $\Theta$ incident with a fixed point of $\Theta^*$ induces a fully embedded polar space of type $\mathsf{D}_4$ in $\Theta$. Conversely, every fully embedded polar space of type $\mathsf{D}_4$ in $\Theta$ arises this way. In particular, there is a unique fully embedded polar space of type $\mathsf{D}_4$ in $\Theta$ containing two points at distance two in $\Theta$.*

*Proof* The first assertion is obvious. So let $\Omega$ be a fully embedded polar geometry of type $\mathsf{D}_4$ of $\Theta$. Consider two non-collinear points $x$, $y$ of $\Omega$. They correspond to two maximal singular subspaces $U_x$, $U_y$, respectively, of $\Theta^*$ intersecting in a point $a$ of $\Theta^*$. Every point $z \in x^\perp \cap y^\perp$ in $\Omega$ corresponds to a maximal singular subspace $U_z$ of $\Theta^*$ intersecting both $U_x$ and $U_y$ in planes. Since these two planes certainly meet inside $U_z$, they both must contain $a$, and so $a \in U_z$. Now every point of $\Omega$ is collinear with two non-collinear points of $x^\perp \cap y^\perp$, and the same argument then implies that the corresponding maximal singular subspace contains $a$. Hence we have shown that $\Omega$ is contained in the residue of $a$, as a fully embedded geometry. But it is now easy to see that it must coincide with that residue, since that residue is also a polar space of type $\mathsf{D}_4$.

The last assertion follows from the fact that $\Theta$ is a parapolar space (see 13.4.2, example 4 of [18]). □

*Remark 6.10* By the previous lemma, and using the same notation, a point of $\Theta^*$ corresponds to a subspace of $\Theta$ isomorphic to a geometry of type $\mathsf{D}_4$. Two points of $\Theta^*$ are collinear if, and only if, the corresponding geometries of type $\mathsf{D}_4$ intersect in a 3-space. Two non-collinear points of $\Theta^*$ correspond to disjoint geometries of type $\mathsf{D}_4$.

**Corollary 6.11** *The standard and hyperbolic $\mathsf{D}_4$s in $\widehat{T}(p,q)$ are the only fully embedded polar spaces of type $\mathsf{D}_4$ in $\widehat{T}(p,q)$ (the latter viewed as a half spin geometry of type $\mathsf{D}_5$). The standard $\mathsf{D}_4$s in $\widehat{T}(p,q)$ arise from points of $\widehat{E}(p,q)$ and the hyperbolic ones arise from imaginary points of $\widehat{E}(p,q)$, i.e. points of $\Theta(\widehat{T}(p,q)) \setminus \widehat{E}(p,q)$.*

*Proof* Let $\Omega$ be a fully embedded polar space of type $\mathsf{D}_4$ in $\widehat{T}(p,q)$. By Lemma 6.9, it arises as the point-residue of a point $x \in \Theta(\widehat{T}(p,q))$ (cf. Definition 5.34). Clearly, $\Omega$ contains two points $y, z$ at distance two, measured in the collinearity graph of $(\mathscr{P}, \mathscr{L})$. Hence, $\{y, z\}$ is either a special pair or an opposite pair in $\Gamma$. If $\{y, z\}$ is special, then $y \bowtie z \in \widehat{E}(p,q)$ by Proposition 5.32(*iii*) and so $\{y, z\}$ is contained in $\Omega_{y \bowtie z}$, the standard $\mathsf{D}_4$ defined by $y \bowtie z$. The last assertion of Lemma 6.9 then implies that $\Omega = \Omega_{y \bowtie z}$ and hence $x = y \bowtie z \in \widehat{E}(p,q)$. On the other hand, if $\{y, z\}$ is opposite, then again Lemma 6.9 implies that $\Omega = \widehat{E}(y, z) \cap \widehat{T}(p,q)$. Since $\Omega$ in this case does not contain lines of $\Gamma$, the point $x$ must be imaginary. □

**Lemma 6.12** *Let $Q$ be a polar space of type $\mathsf{D}_4$ fully embedded in $\widehat{T}(p,q)$ whose lines are the hyperbolic lines of $\widehat{T}(p,q)$ contained in it. Then, the hyperbolic solids $U$ of $Q$ all of whose points are collinear (in $\Gamma$) with a point of $\widehat{E}(p,q)$ are precisely the members of one system $\mathscr{M}$ of maximal singular hyperbolic subspaces of $Q$. Moreover, for each $U$ in $\mathscr{M}$, there is exactly one point $x_U$ in $\widehat{E}(p,q)$ collinear with each point of $U$. The map $U \mapsto x_U$ is injective.*

*Proof* Suppose first that $Q$ arises from the residue of a point $x$ of $\widehat{E}(p,q)$, i.e. $Q = x^\perp \cap \widehat{T}(p,q)$. Then all points of $Q$ are collinear with $x$. In this case, $Q$ clearly contains lines of $\Gamma$, contradicting our hypothesis. Hence $Q$ arises from the residue $R_z$ of an imaginary point $z$ in the imaginary completion of $\widehat{E}(p,q)$ to a polar space $\Theta^*$ of type D$_5$. So, the set of points of $Q$ are in bijective correspondence with one system of maximal singular subspaces of $R_z$, and vice versa (the principle of triality appears here). Hence, the set of lines of $\Theta^*$ through $z$ are in bijective correspondence to one system $\mathscr{M}$ of maximal singular subspaces of $Q$. Now we fix $U \in \mathscr{M}$ and the corresponding line $L_U$ through $z$. Consider the unique point $x_U$ on $L_U$ which belongs to $\widehat{E}(p,q)$—this point exists indeed since $\widehat{E}(p,q)$ is a geometric hyperplane of $\Theta^*$. Then, $U$ is a maximal singular subspace of the subgeometry of type D$_4$ arising from the residue of $x_U$ (because everything different from points of $\Theta^*$ incident with $L_U$ belongs to the residue of $x_U$). Hence all points of $U$ are collinear in $\Gamma$ to $x_U$. This establishes existence. Uniqueness follows immediately from Corollary 5.30. The injectivity of $U \mapsto x_U$ follows from the bijective correspondence between $\mathscr{M}$ and the set of lines of $\Theta^*$ through $z$.

Finally, we prove that no member of the other system $\mathscr{M}'$ of maximal singular subspaces of $Q$ is collinear with a point of $\widehat{E}(p,q)$. Suppose that $U' \in \mathscr{M}'$ is collinear with a point $x' \in \widehat{E}(p,q)$. Then there is a member $U$ of $\mathscr{M}$ intersecting $U'$ in a plane. Since $U$ and $U'$ are contained in $\widehat{E}(p',q')$, for opposite points $p', q'$ in $U \cup U'$, it follows by Proposition 5.32(i) that $x'$ is collinear with the unique point of $\widehat{E}(p,q)$ collinear with $U$, contradicting Lemma 5.17. □

**Lemma 6.13** *(i) Let $x, y$ be two opposite points of $\widehat{T}(p,q)$. Then, $\widehat{E}(p,q) \cap \widehat{T}(x,y)$ endowed with the hyperbolic lines it contains is a polar space of type D$_4$ over $\mathbb{K}$ and a geometric hyperplane of $\widehat{E}(p,q)$. This set is precisely the set of points of $\widehat{E}(p,q)$ collinear with the members of one system of maximal singular subspaces of the polar space $\widehat{E}(x,y) \cap \widehat{T}(p,q)$ of type D$_4$.*

*(ii) Let $Q$ be a geometric hyperplane of the polar space $\widehat{E}(p,q)$ isomorphic to a polar space of type D$_4$. Then, there exist exactly two tropic circle geometries $T_1$ and $T_2$ with $T_i \cap \widehat{E}(p,q) = Q$, $i = 1, 2$.*

*Proof* (i) By Lemma 6.7, $\widehat{E}(x,y) \cap \widehat{T}(p,q)$ endowed with the hyperbolic lines in it, is a polar space of type D$_4$ over $\mathbb{K}$ and a geometric hyperplane of $\widehat{E}(p,q)$. From Lemma 6.12, we know that exactly one of its systems of maximal singular subspaces is such that each of its members is collinear with a unique point of $\widehat{E}(p,q)$. Clearly, the set $Z$ of points thus obtained is contained in $\widehat{E}(p,q) \cap \widehat{T}(x,y)$. Since $\widehat{E}(x,y) \cap \widehat{T}(p,q)$ contains two disjoint maximal singular subspaces, Proposition 5.32(iv) implies that $Z$ contains a pair of opposite points (in $\widehat{E}(p,q) \cap \widehat{T}(x,y)$). It follows from Lemma 6.7 that $Z = \widehat{E}(p,q) \cap \widehat{T}(x,y)$. This proves (i).

(ii) Note that, by symmetry, members of precisely one family of maximal singular subspaces of $\widehat{E}(p,q) \cap \widehat{T}(x,y)$ correspond to points of $\widehat{E}(x,y) \cap \widehat{T}(p,q)$. Now look at the other family of maximal singular subspaces of $\widehat{E}(p,q) \cap \widehat{T}(x,y)$. The corresponding set of points of $\widehat{T}(p,q)$ contains two opposite points, say $x'$ and $y'$. As in the first paragraph of this proof, it follows that this set coincides with $\widehat{E}(x',y') \cap \widehat{T}(p,q)$. Using this argument, we can start from any geometric hyperplane $Q$ of $\widehat{E}(p,q)$ of type D$_4$ and conclude that the points of $\widehat{T}(p,q)$ corresponding to the maximal singular subspaces of $Q$ form two disjoint polar spaces of type D$_4$, which arise as the intersections of $\widehat{T}(p,q)$ with two extended equator geometries $\widehat{E}_1$ and $\widehat{E}_2$, respectively. Let $\widehat{T}_1$ and $\widehat{T}_2$ be their respective tropic circle geometries. By (i), $\widehat{T}_1$ and $\widehat{T}_2$ both intersect $\widehat{E}(p,q)$ in a polar space of type D$_4$ which clearly coincides with $Q$. Since for any tropic circle geometry $\widehat{T}^*$, with corresponding extended equator geometry $\widehat{E}^*$, intersecting $\widehat{E}(p, q)$ in $Q$, the members of one family of maximal singular subspaces of $Q$ have to be collinear with all points of $\widehat{E}^* \cap \widehat{T}(p, q)$, it follows that $\widehat{E}^* \cap \widehat{T}_1 \in \{\widehat{E}_1 \cap \widehat{T}(p, q), \widehat{E}_2 \cap \widehat{T}(p, q)\}$ and hence $\widehat{E}^* \in \{\widehat{E}_1, \widehat{E}_2\}$ (use Proposition 5.22), establishing that $\widehat{T}_1$ and $\widehat{T}_2$ are the only tropic circle geometries intersecting $\widehat{E}(p, q)$ in $Q$. □

Both Corollary 5.38 and Lemma 6.7 produce geometric hyperplanes of $\widehat{E}(p, q)$ of type $\mathsf{D}_4$. The connection between these two constructions is given by the following lemma.

**Lemma 6.14** *Let $r$ be any point not belonging to $\widehat{H}(p, q)$. With reference to Corollary 5.38, let $H_r$ be the subspace of $\widehat{E}(p, q)$ of type $\mathsf{D}_4$ consisting of the points of $\widehat{E}(p, q)$ special to $r$. Then, there exists a unique extended equator geometry $\widehat{E}_r$ containing $r$ and intersecting $\widehat{T}(p, q)$ in a hyperbolic $\mathsf{D}_4$. The tropic circle geometry $\widehat{T}_r$ corresponding to $\widehat{E}_r$ intersects $\widehat{E}(p, q)$ precisely in $H_r$. Also, the set of points of $\widehat{T}(p, q)$ which are symplectic to $r$ is contained in $\widehat{E}_r \cap \widehat{T}(p, q)$ (and hence constitutes a hyperbolic $\mathsf{B}_3$).*

*Proof* Let $\Omega_1$ and $\Omega_2$ be the two natural systems of maximal singular subspaces of $H_r$, and let $\omega_1$ and $\omega_2$ be the corresponding sets of points of $\widehat{T}(p, q)$ (so $\beta(x) \in \Omega_i$ for all $x \in \omega_i$, $i = 1, 2$). We claim that, if $U_i \in \Omega_i$, with $U_1 \cap U_2$ a hyperbolic plane, then exactly one of $\beta(U_1)$ and $\beta(U_2)$ is special to $r$, whereas the other is either symplectic to or opposite $r$.

Indeed, let $h$ be a hyperbolic line in $U_1 \cap U_2$. The elements of $\Omega_i$, $i = 1, 2$, containing $h$ form a (hyperbolic) line in the half spin $\mathsf{D}_5$ geometry corresponding to the dual polar space of type $\mathsf{B}_4$ associated with $\widehat{E}(p, q)$. Hence, the points of $\omega_i$ corresponding to the elements of $\Omega_i$ that contain $h$ form a hyperbolic line $h_i$ in $\widehat{T}(p, q)$. Then $h_1$ and $h_2$ are contained in $h^\perp$, which, in turn, is contained in $S(h)$. Hence, we see that $h^\perp$ is a 3-space when endowed with the lines and hyperbolic lines it contains, and, moreover, $h^\perp$ is entirely contained in $\widehat{T}(p, q)$. Let $x$ be any point of $h^\perp \setminus (h_1 \cup h_2)$. Then $\beta(x) \notin \Omega_1 \cup \Omega_2$, and hence it contains points of $\widehat{E}(p, q)$ that are opposite $r$. Hence $x$ cannot be symplectic to $r$. Note also that by Lemma 5.37 no point of $h^\perp$ is collinear with $r$. Now there are two possibilities.

- *$r$ is close to $S(h)$*. In this case, $r$ is collinear with all points of a line $L$ of $S(h)$, and by our previous remark, $L$ does not meet $h^\perp$. Hence the set of points of $h^\perp$ symplectic to $r$ is $h^\perp \cap L^\perp$, which is a hyperbolic line $g$ of $\Gamma$. But $g$ is disjoint from $h^\perp \setminus (h_1 \cup h_2)$ and must thus be contained in $h_1 \cup h_2$. Obviously, this implies that $g = h_1$ or $g = h_2$. Since $\beta(U_i) \in h_i$, $i = 1, 2$, $r$ is special to $\beta(U_i)$ if $g = h_i$ and opposite $\beta(U_i)$ if $g \neq h_i$. The claim now follows in this case.
- *$r$ is far from $S(h)$*. In this case, $r$ is symplectic to a unique point $s \in S(h)$. Since points of $S(h)$ not collinear with $s$ are opposite $r$, we deduce that $s \in h^\perp$. Hence, as before, $s \in h_1 \cup h_2$, say $s \in h_1$, and then all points of $h_2$ are special to $r$ (since $h_2 \subseteq h_1^\perp$), whereas all points of $h_1 \setminus \{s\}$ are opposite $r$. The claim follows in this case as well.

Now, since the bipartite graph on $\Omega_1 \cup \Omega_2$, where adjacency is intersecting in a 2-space, is connected, it follows that we can choose indices so that all elements of $\omega_1$ are symplectic to or opposite $r$, and all elements of $\omega_2$ are special to $r$. Moreover, from the two cases above, it follows that the set of points in $\omega_1$ symplectic to $r$ is a geometric hyperplane of $\omega_1$, viewed as a polar space of type $\mathsf{D}_4$. Hence there are at least two opposite points in $\omega_1$ which are symplectic to $r$. It follows that the unique extended equator geometry $\widehat{E}_r$ intersecting $\widehat{T}(p, q)$ in $\omega_1$ (see Lemma 6.13) contains $r$.

Now the other assertions follow easily from Lemma 6.13. □

We can now distinguish the points of $\widehat{E}(p, q)$, of $\widehat{T}(p, q)$, of $\widehat{H}(p, q) \setminus (\widehat{E}(p, q) \cup \widehat{T}(p, q))$ and the rest by their relations with the points of $\widehat{E}(p, q)$.

**Lemma 6.15** *Let $x$ be any point of $\Gamma$. Then*

(i) *$x \in \widehat{E}(p,q)$ if, and only if, no point of $\widehat{E}(p,q)$ is collinear or special to $x$. So, in this case, each point of $\widehat{E}(p,q)$ is either equal to, symplectic to or opposite $x$ and all possibilities occur, for all $x \in \widehat{E}(p,q)$.*
(ii) *$x \in \widehat{T}(p,q)$ if, and only if, no point of $\widehat{E}(p,q)$ is equal to, symplectic to, or opposite $x$. So, in this case, each point of $\widehat{E}(p,q)$ is either collinear or special to $x$ and both possibilities occur, for all $x \in \widehat{T}(p,q)$.*
(iii) *$x \in \widehat{H}(p,q) \setminus (\widehat{E}(p,q) \cup \widehat{T}(p,q))$ if, and only if, each point of $\widehat{E}(p,q)$ is either collinear with, symplectic to, special to, or opposite $x$, and all possibilities occur. In this case, there is a unique point $z \in \widehat{E}(p,q)$ collinear with $x$, there is a unique hyperbolic solid $U$ of $\widehat{E}(p,q)$ through $z$ all of whose points except $z$ are symplectic to $x$, all other points of $\widehat{E}(p,q)$ symplectic to $z$ are special to $x$ and all points of $\widehat{E}(p,q)$ opposite $z$ are also opposite $x$.*
(iv) *$x \notin \widehat{H}(p,q)$ if, and only if, each point of $\widehat{E}(p,q)$ is either special to or opposite $x$. In this case, the points of $\widehat{E}(p,q)$ that are special to $x$ form a hyperbolic D$_4$ and hence a geometric hyperplane of $\widehat{E}(p,q)$.*

*Proof* We prove all the "only if" statements and show that all possibilities do occur. Since the cases above are mutually exclusive and exhaustive, the "if" parts then also follow.

If $x \in \widehat{E}(p,q)$, then the result follows from Lemma 5.17. Suppose now $x \in \widehat{T}(p,q)$, and put $\beta(x) = U$. Then all points of $U$ are collinear with $x$. Now let $y \in \widehat{E}(p,q) \setminus U$. Consider a symplecton $S$ through $y$ and a point $u$ of $U$. Then $x$ does not belong to $S$ as otherwise by projecting $y$ onto the line $xu$ we find a second point of $\widehat{T}(p,q)$ on that line, a contradiction to Corollary 5.39. Hence $x$ is close to $S$ and since $y$ is not collinear with $u$, the point $y$ is special to $x$.

Now suppose $x \in \widehat{H}(p,q) \setminus (\widehat{E}(p,q) \cup \widehat{T}(p,q))$. Let $L$ be the unique line of $\Gamma$ through $x$ that intersects both $\widehat{E}(p,q)$ and $\widehat{T}(p,q)$ in respective points $y$ and $z$. Clearly, the point $y$ is collinear with $x$, and all other points of $\beta(z)$ are symplectic to $x$, by Lemma 5.7. Now let $u$ be a point of $\widehat{E}(p,q)$ symplectic to $y$, but not belonging to $\beta(z)$. By considering the symplecton through $u$ and $y$, Fact 5.5(1) implies that $x$ is special to $u$. Finally, let $v$ be a point of $\widehat{E}(p,q)$ opposite $y$. Since $v$ is opposite $y$ and special to $z$, each other point on the line $yz$, in particular $x$, is opposite $v$, as follows from Lemma 5.7.

Finally, if $x \notin \widehat{H}(p,q)$, then the result follows from Lemma 6.14. Now the complete assertion is clear. □

There is a nice consequence.

**Corollary 6.16** *Let $\mathfrak{e}$ be an extended equator geometry intersecting $\widehat{T}(p,q)$ in a hyperbolic D$_4$. Then, no point of $\mathfrak{e} \setminus \widehat{T}(p,q)$ belongs to $\widehat{H}(p,q)$.*

*Proof* Let $x$ be a point of $\mathfrak{e}$ not belonging to $\widehat{T}(p,q)$. Then, by Lemma 6.12, there is one system of maximal singular subspaces of $\mathfrak{e} \cap \widehat{T}(p,q)$ each member of which is collinear with a unique point of $\widehat{E}(p,q)$. By Lemma 6.13(i), these points are precisely the points of $\widehat{E}(p,q) \cap T_{\mathfrak{e}}$, which is a hyperbolic D$_4$ that we will denote by $Q$. Then, since every point of $Q$ is collinear with some point symplectic to $x$, Lemma 5.7 implies that $x$ cannot be opposite any point of $Q$. It follows that $x \notin \widehat{E}(p,q)$, because otherwise some point of $Q$ is opposite $x$. For the same reason, Lemma 6.15(iii) also rules out $x \in \widehat{H}(p,q) \setminus (\widehat{E}(p,q) \cup \widehat{T}(p,q))$. Of course $x \notin \widehat{T}(p,q)$ by assumption, so that only leaves $x \notin \widehat{H}(p,q)$. □

**Proposition 6.17** *Let $\widehat{E}(p,q)$ be an extended equator geometry. Put*

$$X = \widehat{E}(p,q) \cup \left\{ \widehat{E}(x,y) : x, y \in \widehat{T}(p,q), x \text{ opposite } y \right\}.$$

*Then X endowed with the members of $\mathscr{L}$ contained in X, is a polar space of type $\mathsf{D}_5$ over $\mathbb{K}$.*

*Proof* Let $\Theta^* = \Theta(\widehat{T}(p,q))$, as in Definition 5.34. To each point $x$ of $\widehat{E}(p,q)$, we can associate the standard $\mathsf{D}_4$ given by $x^\perp \cap \widehat{T}(p,q)$. Moreover, if $x, y \in \widehat{T}(p,q)$ are two opposite points, then by Lemma 6.7, $\widehat{E}(x,y) \cap \widehat{T}(p,q)$ is a hyperbolic $\mathsf{D}_4$ in $\widehat{T}(p,q)$, hence it corresponds to a point $u$ of $\Theta^* \setminus \widehat{E}(p,q)$. We write $\zeta(u) = \widehat{E}(x,y) \cap \widehat{T}(p,q)$. Clearly, $\widehat{E}(x,y) \cap \widehat{T}(p,q)$ determines $\widehat{E}(x,y)$ and we write $\gamma(\widehat{E}(x,y) \cap \widehat{T}(p,q)) = \widehat{E}(x,y)$. So, by Lemma 6.9, there is a natural bijection $\sigma$ from the point set of $\Theta^*$ to $X$ which is the identity on $\widehat{E}(p,q)$ and $\gamma \circ \zeta$ on $\Theta^* \setminus \widehat{E}(p,q)$. Since the structure of a polar space with given point set is uniquely determined by collinearity, it suffices to show that $\sigma$ preserves collinearity in both directions. Note that two points in $\Theta^*$ are collinear if, and only if, the corresponding standard/hyperbolic $\mathsf{D}_4$s intersect nontrivially, and then they intersect in a hyperbolic solid.

So let $\mathfrak{r}_1$ and $\mathfrak{r}_2$ be two points of $\Theta^*$. There are essentially three possibilities.

- Both $\sigma(\mathfrak{r}_1)$ and $\sigma(\mathfrak{r}_2)$ are points of $\widehat{E}(p,q)$. In this case, the assertion follows since $\sigma$ is the identity. (The fact that we do not seem to have to prove anything here is due to the fact that $\beta$ induces an isomorphism between the half spin geometry $\widehat{T}(p,q)$—endowed with ordinary and hyperbolic lines it—and $\widehat{E}(p,q)$, viewed as half spin geometry by considering its imaginary completion.)
- Both $\sigma(\mathfrak{r}_1) = \mathfrak{e}$ and $\sigma(\mathfrak{r}_2) = \mathfrak{e}'$ are new points. If $\mathfrak{r}_1$ and $\mathfrak{r}_2$ are collinear in $\Theta^*$, then we know that their images under $\zeta$ intersect in a hyperbolic solid. Then, clearly, also their images under $\sigma$ do.
  Conversely, if $\mathfrak{e}$ and $\mathfrak{e}'$ are collinear, we show that the hyperbolic solid $\mathfrak{e} \cap \mathfrak{e}'$ is contained in $\widehat{T}(p,q)$. By Corollary 6.16, a point of $\mathfrak{e} \cup \mathfrak{e}'$ outside $\widehat{T}(p,q)$ is not contained in $\widehat{H}(p,q)$ and hence, by Lemma 6.14, it is contained in a *unique* extended equator geometry intersecting $\widehat{T}(p,q)$ in a hyperbolic $\mathsf{D}_4$. Hence $\mathfrak{e} \cap \mathfrak{e}' \subseteq \widehat{T}(p,q)$ and the assertion follows.
- $\sigma(\mathfrak{r}_1) = \mathfrak{e}$ is a new point and $\sigma(\mathfrak{r}_2) = x$ is a point of $\widehat{E}(p,q)$. If $\mathfrak{e}$ and $x$ are collinear, then $x \in T_\mathfrak{e}$. So, $x$ is collinear in $\Gamma$ to the points of a projective 3-space of $\mathfrak{e}$, which is entirely contained in $\mathfrak{e} \cap \widehat{T}(p,q)$ by Lemma 5.37. Hence $\mathfrak{r}_1$ and $\mathfrak{r}_2$ are collinear in $\Theta^*$. Now assume that $\mathfrak{r}_1$ and $\mathfrak{r}_2$ are collinear in $\Theta^*$. This means that $\mathfrak{e}$ and $x^\perp$ intersect in a hyperbolic solid, which implies that $x \in T_\mathfrak{e}$ and the assertion is proved. □

**Definition 6.18** (*The quads*) We will denote the polar space of the previous proposition by $\Sigma(\widehat{E}(p,q))$. Note that, alternatively, one can define this as

$$\Sigma(\widehat{E}(p,q)) = \widehat{E}(p,q) \cup \left\{ \mathfrak{e} : \mathfrak{e} \cap \widehat{T}(p,q) \text{ is a hyperbolic } \mathsf{D}_4 \right\}.$$

We call the polar spaces of type $\mathsf{D}_5$ of Propositions 6.6 and 6.17 the *quads* of $(\mathscr{P}, \mathscr{L})$ and denote the family of quads by $\mathscr{Q}$.

After $\mathscr{P}$ and $\mathscr{L}$, the family of quads defines the third type of vertices of the associated building of type $\mathsf{E}_6$ we wish to construct. The points have type 1, the lines have type 3, the quad having type 6. Our next goal is to define the elements of the remaining three types and the incidence relation; and to show that the structure obtained is indeed a building of type $\mathsf{E}_6$ on which a symplectic polarity acts with corresponding fixed point building exactly the building of type $\mathsf{F}_4$ associated with $\Gamma$.

### 6.4 Maximal singular 4-spaces of $(\mathscr{P}, \mathscr{L})$

We need to define the elements of type 4 (these will be singular planes of $(\mathscr{P}, \mathscr{L})$), the elements of type 5 (these will be certain singular 4-spaces of $(\mathscr{P}, \mathscr{L})$) and the elements of type 2 (these will be singular 5-spaces of $(\mathscr{P}, \mathscr{L})$) in the geometry $(\mathscr{P}, \mathscr{L})$. We start with the set $\mathscr{U}$ of elements of type 5, defined as the set of intersections of two (distinct) quads which meet in at least two collinear points. In the next three lemmas we introduce three types of members of $\mathscr{U}$, prove that these are singular 4-spaces, and in Lemma 6.22 we show that these are all the elements of $\mathscr{U}$. Note that $\mathscr{U}$ will certainly not consist of *all* singular 4-spaces of $(\mathscr{P}, \mathscr{L})$. However, $\mathscr{U}$ is the family of all *maximal* singular 4-spaces. We will not show this, as this will not be needed (but in fact it is easy to do). We will later on (see Lemma 6.28) construct the 5-spaces, which will be maximal singular subspaces of dimension 5 of $(\mathscr{P}, \mathscr{L})$. The hyperplanes of the 5-spaces will be the 4'-spaces. The 4-spaces (i.e. the elements of $\mathscr{U}$) and the 4'-spaces are all the 4-dimensional singular subspaces of $(\mathscr{P}, \mathscr{L})$.

**Lemma 6.19** *Let $L$ be a line of $\Gamma$. Then, the set $L^\perp$, endowed with all ordinary and hyperbolic lines of $\Gamma$ contained in it, is a projective 4-space over $\mathbb{K}$, denoted by $U(L)$. Moreover, $U(L) = x_1^\perp \cap x_2^\perp = \Sigma(x_1) \cap \Sigma(x_2) \subseteq \Gamma$ for every pair $\{x_1, x_2\}$ of distinct points of $L$, and so $U(L) \in \mathscr{U}$.*

*Proof* Consider a point $p \in L$, then all points of $L^\perp$ belong to $\Sigma(p)$. Clearly, every pair of points of $L^\perp$ is collinear or symplectic, and the line or hyperbolic line joining them is completely contained in $L^\perp$, by Lemmas 5.6 and 5.7. Hence, $L^\perp$ is a singular subspace of $\Sigma(p)$ of dimension at most 4, by Proposition 6.6. Since $L^\perp \cap S$, for any symplecton $S$ containing $L$, is a projective space of dimension 3, and is properly contained in $L^\perp$, the dimension of $L^\perp$ must be 4. The last assertion of the lemma follows from Lemma 5.7, the definition of $\Sigma(p)$ and the fact that the intersection of $\Sigma(x_1)$ and $\Sigma(x_2)$ does not contain new points. □

**Lemma 6.20** *Let $h$ be a hyperbolic line of $\Gamma$. Then $U(h) = h^\perp \cup \{\mathfrak{e} \in \mathscr{E} : h \subseteq \mathfrak{e}\} \subseteq \mathscr{P}$, endowed with the lines of $\mathscr{L}$ contained in it, is a projective 4-space over $\mathbb{K}$. Moreover, $U(h) = \Sigma(x_1) \cap \Sigma(x_2)$ for any pair $\{x_1, x_2\}$ of distinct points of $h$, and so $U(h) \in \mathscr{U}$.*

*Proof* We first claim that $U(h)$ is a singular subspace of $(\mathscr{P}, \mathscr{L})$. Indeed, since any two points of $h^\perp$ lie in $S(h)$, they are collinear or symplectic and the line joining them is contained in $h^\perp$ as well. As any two new points $\mathfrak{e}$ and $\mathfrak{e}'$ of $U(h)$ have $h$ in common, Lemma 6.5 implies that $\mathfrak{e}$ and $\mathfrak{e}'$ are collinear in $(\mathscr{P}, \mathscr{L})$. Further, all new points of $\langle \mathfrak{e}, \mathfrak{e}' \rangle$ also contain $h$. The unique ordinary point of that new line must belong to $h^\perp$ and hence also belongs to $U(h)$. Finally, if we consider a point $x \in h^\perp$ and $\mathfrak{e} \in \mathscr{E}$ with $h \subseteq \mathfrak{e}$, then $x$ is collinear with at least two points of $\mathfrak{e}$ and hence belongs to $T_\mathfrak{e}$. Consequently, by definition, $x$ is collinear with $\mathfrak{e}$ in $(\mathscr{P}, \mathscr{L})$, and all new points of $\langle x, \mathfrak{e} \rangle$ contain $x^\perp \cap \mathfrak{e} \supseteq h$. The claim is proved.

But now we see that $U(h)$ is contained in $\Sigma(x)$, for each $x \in h$. Since $U(h)$ properly contains the 3-space $h^\perp$, Proposition 6.6 implies that its dimension is 4. It follows that, if $x_1$ and $x_2$ are two distinct points of $h$, then $U(h)$ belongs to $\Sigma(x_1) \cap \Sigma(x_2)$. Now, Lemma 5.6 and Proposition 5.9 readily imply that $U(h) = \Sigma(x_1) \cap \Sigma(x_2)$. □

**Lemma 6.21** *Let $V$ be a hyperbolic solid and let $x = \beta(V)$. Let $P_V$ be the set of points of $\Gamma$ collinear with at least two points of $V$. Then,*

$$U(V) = V \cup \{\mathfrak{e} \in \mathscr{E} : \mathfrak{e} \cap P_V \text{ is a hyperbolic } \mathsf{D}_3\text{-cone with vertex } x\}$$

is a projective 4-space in $(\mathcal{P}, \mathcal{L})$. Moreover, $U(V) = \Sigma(x) \cap \Sigma(\mathfrak{f})$ for any extended equator geometry $\mathfrak{f}$ containing $V$, and so $U(V) \in \mathcal{U}$.

*Proof* Since $x \in \mathfrak{e}$ for every $\mathfrak{e} \in \mathcal{E} \cap U(V)$, we see that $U(V) \subseteq \Sigma(x)$. We now show that $U(V)$ is a singular subspace of $\Sigma(x)$. Towards that aim, we consider two arbitrary points of $U(V)$, show that they are collinear in $(\mathcal{P}, \mathcal{L})$, and that all points of the line connecting them belong to $U(V)$. There are three possibilities.

- If both points belong to $V$, then the assertion follows from the definition of a hyperbolic solid.
- Suppose $v \in V$ and $\mathfrak{e} \in U(V) \cap \mathcal{E}$. We fix an arbitrary extended equator geometry $\mathfrak{f}$ containing $V$. The definition of $T_\mathfrak{f}$ implies that $P_V \subseteq T_\mathfrak{f}$. Hence, $\mathfrak{e} \cap T_\mathfrak{f}$ contains $\mathfrak{e} \cap P_V$, which in turn contains an opposite pair of points, as $\mathfrak{e} \in U(V)$. Lemma 6.7 implies that $\mathfrak{e} \cap T_\mathfrak{f}$ is a hyperbolic $\mathsf{D}_4$ in $T_\mathfrak{f}$. Hence there is a point $x' \in \mathfrak{e} \cap T_\mathfrak{f}$ opposite $x$. Let $V' = x'^\perp \cap \mathfrak{f}$. By the proof of Lemma 6.7, the points of $T_\mathfrak{f} \cap E(x, x')$ map through $\beta$ (relative to $\mathfrak{f}$) to all hyperbolic solids of $\mathfrak{f}$ intersecting both $V$ and $V'$ in hyperbolic lines. Now, as $v$ is on such a hyperbolic line, it is collinear with at least one point of $T_\mathfrak{f} \cap E(x, x')$. Since $v$ is also collinear with $x$, it belongs to $T_\mathfrak{e}$. Then, by definition, $\mathfrak{e}$ and $v$ are collinear as points of $(\mathcal{P}, \mathcal{L})$.
  Since $v \notin V'$, there is a unique hyperbolic solid $W$ of $\mathfrak{f}$ containing $v$ and intersecting $V'$ in a hyperbolic plane. Put $\beta(W) = w$. Then $x' \perp w \perp v \perp x$, and so, by Proposition 6.1, every new point $\mathfrak{e}'$ of the line $\langle \mathfrak{e}, v \rangle$ contains a point $x''$ of $x'w \setminus \{w\}$, more precisely, $\mathfrak{e}' = \widehat{E}(x, x'')$. By Proposition 5.32(i), $x'' \in T_\mathfrak{f}$, so by Corollary 6.8, $x^\perp \cap T_\mathfrak{f} \cap \widehat{E}(x, x'')$, which equals $P_V \cap \widehat{E}(x, x'')$, is a hyperbolic $\mathsf{D}_3$-cone with vertex $x$. Hence all points of the new line $\langle v, \mathfrak{e} \rangle$ belong to $U(V)$.
- Suppose $\mathfrak{e}, \mathfrak{e}' \in U(V) \cap \mathcal{E}$. Then, $\mathfrak{e} \cap P_V$ and $\mathfrak{e}' \cap P_V$ are two cones (with vertex $x$) over standard line Grassmannians of projective 3-spaces inside a cone (with vertex $x$) over the line Grassmannian of a projective 4-space, see Lemma 6.4. Hence their intersection is a cone (with vertex $x$) over the line Grassmannian of a plane (the intersection of the two 3-spaces inside the 4-space). It follows that $\mathfrak{e} \cap \mathfrak{e}'$ is a hyperbolic solid, contained in $P_V$.
  Let $\mathfrak{f}$ be as above. Then, $\mathfrak{e} \cap T_\mathfrak{f}$ and $\mathfrak{e}' \cap T_\mathfrak{f}$ are two hyperbolic $\mathsf{D}_4$s in $T_\mathfrak{f}$. By Lemma 6.9, they correspond to imaginary points $z, z'$ of $\Theta(T_\mathfrak{f})$. Since $x$ is contained in both $\mathfrak{e} \cap T_\mathfrak{f}$ and $\mathfrak{e}' \cap T_\mathfrak{f}$, the points $z$ and $z'$ are collinear in $\Theta(T_\mathfrak{f})$. Hence, the joining line contains a unique point $v$ of $\mathfrak{f}$, and clearly the corresponding standard $\mathsf{D}_4$ also contains $x$. Hence, $v \in V$. We now see that $v \in \langle \mathfrak{e}, \mathfrak{e}' \rangle$ and the assertion follows from the second case.

Hence, $U(V)$ is a singular subspace of $\Sigma(x)$, and also of $\Sigma(\mathfrak{f})$, with $\mathfrak{f}$ as above. Since it properly contains the 3-space $V$, it has dimension 4. It remains to show that $\Sigma(x) \cap \Sigma(\mathfrak{f}) = U(V)$. Clearly, the only points of $\Gamma$ contained in both $\Sigma(x)$ and $\Sigma(\mathfrak{f})$ are the points of $V$. Suppose a new point $\mathfrak{e} \in \mathcal{E}$ is contained in $\Sigma(x) \cap \Sigma(\mathfrak{f})$. Then, by definition of $\Sigma(\mathfrak{f})$ and Lemma 6.7, we have that $\mathfrak{f} \cap T_\mathfrak{e}$ is a hyperbolic $\mathsf{D}_4$, which implies by Lemma 6.13(i) that $\mathfrak{e} \cap T_\mathfrak{f}$ is a hyperbolic $\mathsf{D}_4$. Since also $x \in \mathfrak{e} \cap T_\mathfrak{f}$ and $P_V$ is the set of all points of $T_\mathfrak{f}$ collinear or symplectic to $x$, by Corollary 6.8, the intersection $\mathfrak{e} \cap P_V$ is a hyperbolic $\mathsf{D}_3$-cone with vertex $x$. Hence $\mathfrak{e} \in U(V)$.

In view of the natural bijection between the elements $L$ of $\mathcal{F}$ and the hyperbolic solids $V$ appearing as the intersection of any pair of collinear new points of $L$ (see Proposition 6.1 and Definition 6.2), we define $U(L)$ as $U(V)$. This, together with Lemmas 6.19 and 6.20, defines $U(\mathfrak{l})$ for every $\mathfrak{l} \in \mathcal{L}$. We note that, for all $\mathfrak{l}, \mathfrak{l}' \in \mathcal{L}$, we have $U(\mathfrak{l}) = U(\mathfrak{l}')$ if and only if $\mathfrak{l} = \mathfrak{l}'$ (this is easy to see; we refer to this property as *the injectivity of $U(\cdot)$*). We now show that $\mathscr{U} = \{U(L) \mid L \in \mathscr{L}\}$.

**Lemma 6.22** *Any two quads containing a common line of $(\mathscr{P}, \mathscr{L})$ intersect in a singular 4-dimensional subspace of $(\mathscr{P}, \mathscr{L})$. A singular 4-dimensional subspace of $(\mathscr{P}, \mathscr{L})$ is the intersection of two quads if, and only if, it is of the form $U(\mathfrak{l})$, with $\mathfrak{l} \in \mathscr{L}$. Moreover, at least one of these quads can be chosen to be of type $\Sigma(p)$, with $p$ a point of $\Gamma$.*

*Proof* By Lemmas 6.19, 6.20 and 6.21, it remains to show that the intersection of two quads sharing at least one line is of the form $U(\mathfrak{l})$, with $\mathfrak{l} \in \mathscr{L}$.

Let $\Sigma, \Sigma'$ be two distinct quads sharing at least two collinear points. Since they then share a line of $(\mathscr{P}, \mathscr{L})$, and since every member of $\mathscr{L}$ contains a point of $\Gamma$, we already know that $\Sigma$ and $\Sigma'$ share a point $s$ of $\Gamma$. There are now three possibilities.

- $\Sigma = \Sigma(x)$ and $\Sigma' = \Sigma(x')$, *for distinct points* $x, x'$ *of* $\Gamma$. If $x \perp x'$, then $\Sigma \cap \Sigma' = U(xx')$ by Lemma 6.19. If $x \perp\!\!\!\perp x'$, then $\Sigma \cap \Sigma' = U(h(x, x'))$ by Lemma 6.20. Now $x$ and $x'$ cannot be special, as otherwise $s$ would be equal to $x^\perp \cap x'^\perp = x \bowtie x'$. Also, by Lemma 5.17, no new point would be contained in $\Sigma \cap \Sigma'$. Finally, $x$ and $x'$ cannot be opposite either, since in that case $\Sigma$ and $\Sigma'$ would not share any point of $\Gamma$ as $x^\perp \cap x'^\perp = \emptyset$.
- $\Sigma = \Sigma(x)$, *for some point* $x$ *of* $\Gamma$, *and* $\Sigma' = \Sigma(\mathfrak{e})$, *for some* $\mathfrak{e} \in \mathscr{E}$. Let $\mathfrak{e} = \widehat{E}(p, q)$ for a pair of opposite points $p, q \in \mathfrak{e}$. Then $x \in \widehat{H}(p, q)$, as $s \in x^\perp \cap \mathfrak{e}$. If a new point $\mathfrak{f}$ would belong to $\Sigma \cap \Sigma'$, then $x \in \mathfrak{f} \cap \widehat{H}(p, q)$, contradicting Corollary 6.16. Hence $\Sigma \cap \Sigma' \subseteq \Gamma$ and, and the structure of $\Sigma'$ implies that any line in $\Sigma \cap \Sigma'$ is a hyperbolic line $h$ which is contained in $\widehat{E}(p, q)$. The structure of $\Sigma$ implies that $x \notin h$ and $x$ is collinear to all points of $h$. Consequently, $x \in \widehat{T}(p, q)$ and $\Sigma \cap \Sigma' = U(x^\perp \cap \mathfrak{e})$, as follows from Lemma 6.21.
- $\Sigma = \Sigma(\mathfrak{e})$ *and* $\Sigma' = \Sigma(\mathfrak{e}')$, *for some distinct* $\mathfrak{e}, \mathfrak{e}' \in \mathscr{E}$. Since $\Sigma$ and $\Sigma'$ must have a point of $\Gamma$ in common, there are only two possibilities, in view of Lemma 6.5. The first one is that $\mathfrak{e}$ and $\mathfrak{e}'$ are collinear in $(\mathscr{P}, \mathscr{L})$. Put $V = \mathfrak{e} \cap \mathfrak{e}'$ and let $x = \beta(V)$, so $x$ is the unique point of $\Gamma$ on the new line $\langle \mathfrak{e}, \mathfrak{e}' \rangle$. Then, by Lemma 6.21, $U(\langle \mathfrak{e}, \mathfrak{e}' \rangle) = U(V) = \Sigma \cap \Sigma(x) = \Sigma' \cap \Sigma(x)$. Hence, $U(\langle \mathfrak{e}, \mathfrak{e}' \rangle) \subseteq \Sigma \cap \Sigma'$. We now show that $U(\langle \mathfrak{e}, \mathfrak{e}' \rangle) = \Sigma \cap \Sigma'$. Clearly, the points of $\Gamma$ in $\Sigma \cap \Sigma'$ are precisely those of $V$ and hence they also belong to $U(\langle \mathfrak{e}, \mathfrak{e}' \rangle)$. Now suppose $\mathfrak{f}$ is a new point of $\Sigma \cap \Sigma'$. In order to show that $\mathfrak{f}$ belongs to $U(\langle \mathfrak{e}, \mathfrak{e}' \rangle)$, it suffices to show that $\mathfrak{f} \in \Sigma(x)$, i.e. $x \in \mathfrak{f}$. Since $\mathfrak{f} \in \Sigma \cap \Sigma'$, both $\mathfrak{e} \cap T_\mathfrak{f}$ and $\mathfrak{e}' \cap T_\mathfrak{f}$ are hyperbolic $\mathsf{D}_4$s in $T_\mathfrak{f}$. We claim that $\mathfrak{e} \cap \mathfrak{e}' \cap T_\mathfrak{f} = V$. Indeed, suppose that a point $v \in V$ does not belong to $T_\mathfrak{f}$. Then Corollary 6.16 implies that we can apply Lemma 6.14 to $v$ and $T_\mathfrak{f}$ and conclude that $\mathfrak{e} = \mathfrak{e}'$, a contradiction. The claim follows. It now follows from Lemma 6.9 that the imaginary points $z$ and $z'$ of $\Theta(T_\mathfrak{f})$ that correspond to $\mathfrak{e} \cap T_\mathfrak{f}$ and $\mathfrak{e}' \cap T_\mathfrak{f}$, respectively, are collinear. Hence the unique point of $\mathfrak{f}$ on the imaginary line joining $z$ and $z'$ is collinear with $V$. By Lemma 5.30, this point coincides with $x$ and so $x \in \mathfrak{f}$.
  The second possibility is that $\mathfrak{e} \cap \mathfrak{e}'$ is a single point $x$. Clearly, $x$ is the unique point of $\Gamma$ in $\Sigma \cap \Sigma'$. Hence it suffices to show that no new point of $\Sigma \cap \Sigma'$ is collinear with $x$. Suppose, for a contradiction, that there is a new point $\mathfrak{f}$ of $\Sigma \cap \Sigma'$ collinear with $x$. Then, $T_\mathfrak{f}$ contains $x$ and intersects both $\mathfrak{e}$ and $\mathfrak{e}'$ in hyperbolic $\mathsf{D}_4$s. But it follows from Lemma 6.9 that two hyperbolic $\mathsf{D}_4$s in $T_\mathfrak{f}$ are either disjoint or share a 3-space, a contradiction to $(\mathfrak{e} \cap T_\mathfrak{f}) \cap (\mathfrak{e}' \cap T_\mathfrak{f}) = \{x\}$.

The lemma is proved. □

**Corollary 6.23** *For every $\mathfrak{l} \in \mathscr{L}$, and for every pair of distinct points $\mathfrak{p}, \mathfrak{q} \in \mathfrak{l}$, we have $U(\mathfrak{l}) = \Sigma(\mathfrak{p}) \cap \Sigma(\mathfrak{q})$.*

*Proof* If $L$ is contained in $\Gamma$, then this follows from Lemmas 6.19 and 6.20. If $L$ is not contained in $\Gamma$, then Lemma 6.21 implies that $U(L) = \Sigma(x) \cap \Sigma(\mathfrak{f})$ for any new point $\mathfrak{f}$ of $L$ and $x$ the unique point of $L$ in $\Gamma$. Now let $\mathfrak{f}$ and $\mathfrak{f}'$ be two distinct new points of $L$. Then $U(L) \subseteq \Sigma(\mathfrak{f}) \cap \Sigma(\mathfrak{f}')$ and equality follows from Lemma 6.22. □

**Lemma 6.24** *Let $p \in \mathscr{P}$ and $L \in \mathscr{L}$. Then $U(L) \subseteq \Sigma(p)$ if, and only if, $p \in L$.*

*Proof* The "if"-part follows from Corollary 6.23. We now show the "only if"-part. So suppose $U(L) \subseteq \Sigma(p)$. We consider the different cases for $L$ separately.

- If $L$ is a line of $\Gamma$, then $U(L)$ cannot be contained in an extended equator geometry by Lemma 5.17, hence $p$ is a point of $\Gamma$, which is collinear with every point of $U(L)$. Since $\Gamma$ does not contain 3-spaces whose lines are ordinary lines of $\Gamma$, $p$ must be contained in every plane of $\Gamma$ in $U(L)$, hence $p \in L$.
- If $L$ is a hyperbolic line of $\Gamma$, then, again, clearly $U(L)$ cannot be contained in a extended equator geometry since $V_L = U(L) \cap \Gamma$ contains lines of $\Gamma$. Hence $p$ is a point of $\Gamma$, which must belong to $V_L^\perp = L$.
- If $L$ is a new line and $p$ is a point of $\Gamma$, then $p$ is the unique point of $\Gamma$ on $L$ by Lemma 5.29. If $p$ is a new point, then, as an equator geometry, it contains all of $\Gamma$ of $U(L)$, and hence it belongs to $L$ by the definition of new lines and Proposition 6.1.

This completes the proof of the lemma. □

**Lemma 6.25** *Two distinct members of $\mathscr{U}$ intersect in a projective subspace of $(\mathscr{P}, \mathscr{L})$ of dimension at most 2.*

*Proof* Suppose $U_1, U_2 \in \mathscr{U}$ are such that their intersection is a singular 3-space. Let $U_i = U(L_i), i = 1, 2$, with $L_i \in \mathscr{L}$. Let $p_i \in L_i$ be a point of $\Gamma$ and let $\mathfrak{q}_i \in L_i$ be a second point on $L_i$, possibly in $\Gamma$, possibly in $\mathscr{E}$. By Corollary 6.23, we have that $U_1 = \Sigma(p_1) \cap \Sigma(\mathfrak{q}_1)$ and $U_2 = \Sigma(p_2) \cap \Sigma(\mathfrak{q}_2)$. Since $U_1 \cap U_2$ is a singular 3-space, at least one of $U_1 \cap \Sigma(p_2)$, $U_1 \cap \Sigma(\mathfrak{q}_2)$ equals $U_1 \cap U_2$ and hence $U_1 \cap U_2 = \Sigma(p_1) \cap \Sigma(\mathfrak{q}_1) \cap \Sigma(x)$ for some $x \in \{p_2, \mathfrak{q}_2\}$.

We claim that $U_1 \cap U_2$ is the intersection of three quads related to three points that are pairwise collinear in $(\mathscr{P}, \mathscr{L})$, at most one of which is related to a new point. Indeed, suppose $\mathfrak{q}_1, \mathfrak{q}_2 \in \mathscr{E}$ and that we cannot choose $x$ equal to $p_2$. Then $\Sigma(p_1) \cap \Sigma(\mathfrak{q}_1) \cap \Sigma(p_2) = U_1$. Moreover, by Lemma 6.21, $U_1$ contains a hyperbolic solid $V$. Hence both $p_1$ and $p_2$ are collinear with all points of $V$, so $p_1 = p_2$ by Corollary 5.30. As $\Sigma(\mathfrak{q}_1) \cap \Sigma(\mathfrak{q}_2)$ contains $U_1 \cap U_2$, Lemma 6.22 implies that $\Sigma(\mathfrak{q}_1) \cap \Sigma(\mathfrak{q}_2) = U(\mathfrak{l})$, for some $\mathfrak{l} \in \mathscr{L}$. Then, Lemma 6.24 implies that both $\mathfrak{q}_1$ and $\mathfrak{q}_2$ belong to $\mathfrak{l}$, hence, in particular, they are collinear. If $p_3$ is the ordinary point of $\Gamma$ on the line $\langle \mathfrak{q}_1, \mathfrak{q}_2 \rangle$, then, by Corollary 6.23, we have $\Sigma(p_3) \cap \Sigma(\mathfrak{q}_1) = \Sigma(\mathfrak{q}_1) \cap \Sigma(\mathfrak{q}_2)$. Hence $U_1 \cap U_2 = \Sigma(p_1) \cap \Sigma(p_3) \cap \Sigma(\mathfrak{q}_1)$, which proves the claim. As above, it follows that $p_1, p_3, \mathfrak{q}_1$ are pairwise collinear in $(\mathscr{P}, \mathscr{L})$.

Hence there are two possibilities. We forget the above notation in the rest of the proof.

- $U_1 \cap U_2 = \Sigma(p_1) \cap \Sigma(p_2) \cap \Sigma(p_3)$, with $p_1, p_2, p_3$ points of $\Gamma$. By the above, $p_i$ and $p_j$ for $i \neq j$ are either collinear are symplectic. Suppose first that $p_1 \perp p_2$. Every 3-space inside the 4-space $U(p_1 p_2)$ contains a point of $p_1 p_2$. Hence we may assume that $p_3 \perp p_1$. If $p_2$ is not collinear with $p_3$ (in $\Gamma$), then $p_3^\perp \cap p_1^\perp \cap p_2^\perp$, which is 3-dimensional as it is precisely $U_1 \cap U_2$, contains a hyperbolic plane. Lemma 5.21 and Theorem 5.33 imply that $p_3 \in p_1p_2$, a contradiction, as otherwise $U_1 = U_2$. Hence $p_3 \perp p_2$. So $p_3 \in p_1^\perp \cap p_2^\perp = \Sigma(p_1) \cap \Sigma(p_2)$, and then it is easy to see, since planes are the maximal singular subspaces of $\Gamma$, that $U_1 \cap U_2$, which is contained in $p_1^\perp \cap p_2^\perp \cap p_3^\perp$, has dimension at most 2, a contradiction.

Hence, we may assume that $p_1 \perp\!\!\!\perp p_2 \perp\!\!\!\perp p_3 \perp\!\!\!\perp p_1$. So $p_1, p_2, p_3$ are contained in a hyperbolic plane, and as in the previous paragraph, we deduce that $p_1^\perp \cap p_2^\perp \cap p_3^\perp$ is a line $L$. Since $L$ is a subspace of $U_1 \cap U_2$ of codimension at most 1, the dimension of $U_1 \cap U_2$ is at most 2, a contradiction.

- $U_1 \cap U_2 = \Sigma(p_1) \cap \Sigma(p_2) \cap \Sigma(\mathfrak{e})$, with $p_1, p_2$ points of $\Gamma$ and $\mathfrak{e} \in \mathscr{E}$. As above, either $p_1 \perp p_2$ or $p_1 \perp\!\!\!\perp p_2$. If $p_1 \perp p_2$, then $U_1 \cap U_2$ does not contain new points and hence, since it also contained in $\Sigma(\mathfrak{e})$, it is a hyperbolic solid in $p_1^\perp \cap p_2^\perp$. However, every singular 3-space $S$ of $(\mathscr{P}, \mathscr{L})$ in $p_1^\perp \cap p_2^\perp$ contains a point of $p_1p_2$, which is then collinear with all other elements of $S$, so $S$ cannot be hyperbolic, a contradiction. If $p_1 \perp\!\!\!\perp p_2$, then $p_1^\perp \cap p_2^\perp$ does not contain a hyperbolic plane, contradicting the fact that $p_1^\perp \cap p_2^\perp \cap \mathfrak{e}$ is a subspace of $U_1 \cap U_2$ of codimension at most 1. □

**Lemma 6.26** *Let $W$ be a hyperbolic solid in $\Gamma$ and let $x = \beta(W)$. Let $V^+$ be the set of points on the lines of $\Gamma$ joining points of $W$ to $x$ and let $V^-$ be the set of elements of $x^\perp$ collinear with at least two points of $W$. Then,*

(i) *$V^+$ and $V^-$ are maximal singular subspaces of $\Sigma(x)$. The lines of $\Gamma$ in $V^+$ and $V^-$ are precisely the lines in $V^+$ and $V^-$ through $x$. Hence both $V^+$ and $V^-$ are cones over hyperbolic solids with vertex $x$.*
(ii) *Let $W'$ be a hyperbolic solid in $V^- \setminus \{x\}$. Then, $x = \beta(W')$ and $V^+$ coincides with the set of elements of $x^\perp$ collinear with at least two points of $W'$.*
(iii) *$V^-$ is independent of the choice of the hyperbolic solid $W$ in $V^+$.*

*Proof* (i) First, note that $W$ is a singular subspace of the quad $\Sigma(x)$. Since $x$ is collinear with all points of $W$, it follows that $x$ and $W$ generate a singular subspace of $\Sigma(x)$ which must necessarily be of dimension 4 and clearly coincides with $V^+$.

Let $v \in V^- \setminus \{x\}$. We first claim that $v^\perp \cap V^+$ is a singular 3-dimensional subspace of $(\mathscr{P}, \mathscr{L})$ of $V^+$ containing $x$. Indeed, $v$ is collinear with at least two points of $\mathfrak{e}$, where $\mathfrak{e}$ is an arbitrary extended equator geometry containing $W$. Hence $v$ belongs to $T_\mathfrak{e}$. Since $v \perp x$, Proposition 5.32(i) implies that $\beta(v) \cap \beta(x)$ is a hyperbolic plane $\pi$ in $\mathfrak{e}$. By Lemma 5.6, $v$ is collinear with all points of the 3-space $Z$ of $V^+$ generated by $x$ and $\pi$. If $Z \subsetneq v^\perp \cap V^+$, then Lemmas 5.6 and 5.7 imply that $W \subseteq v^\perp$, contradicting Corollary 5.30. The claim is proved.

In the notation of Lemma 6.4, we have $V^- \subseteq P_W$. Since for every $u \in V^- \setminus \{x\}$, the intersection $u^\perp \cap W$ is a hyperbolic plane by the previous claim, $V^- \subseteq P_W^\Gamma$. The inverse inclusion follows analogously. Hence $V^-$ is a cone over a hyperbolic solid and (i) is proved.

(ii) Now it follows from Lemma 6.4(ii) that $V^-$, endowed with the ordinary and hyperbolic lines contained in it, is also a maximal singular subspace of $\Sigma(x)$. The same lemma implies that $V^-$ is a cone over a hyperbolic solid, say $W'$, and we have $x = \beta(W')$. By the previous paragraph, collinearity induces a duality from $W'$ and $W$, because a point of $W'$ is collinear with a hyperbolic plane in $W$. Hence $V^+$ is the set of points of $\Gamma$ collinear with $x$ and collinear with at least two points of $W'$.

(iii) This follows immediately from the fact that, by Lemma 5.7 $V^-$ can be defined as the set of points collinear with all points of at least two (ordinary) lines of $V^+$ (necessarily containing $x$). □

**Definition 6.27** The set $V^+$ of the previous lemma will be called a *hyperbolic cone (with vertex x)* and $V^-$, which is a hyperbolic cone too, is called the *twin* of $V^+$. The previous lemma also implies that the twin of $V^-$ is $V^+$.

## 6.5 Singular 5-spaces of $(\mathscr{P}, \mathscr{L})$

We will now describe two types of singular 5-spaces of $(\mathscr{P}, \mathscr{L})$. In fact, together with the subspaces $U(L)$, $L \in \mathscr{L}$, defined earlier, these will be all the maximal singular subspaces of $(\mathscr{P}, \mathscr{L})$, but there is no need to show this as it will follow once we have proved that $(\mathscr{P}, \mathscr{L})$ defines a building of type $\mathsf{E}_6$.

**Lemma 6.28** *(i) Every symplecton of $\Gamma$, endowed with the ordinary and hyperbolic lines contained in it, is a projective 5-space.*
*(ii) Let $V^+$ be a hyperbolic cone in $\Gamma$. Define*

$$M(V^+) = V^{\cup} \left\{ \mathfrak{f} \in \mathscr{E} : V^+ \subseteq T_{\mathfrak{f}} \right\}.$$

*Then, $M(V^+)$ endowed with the members of $\mathscr{L}$ contained in it is a projective 5-space.*

*Proof* Assertion $(i)$ is clear. We show $(ii)$.

Let $V^-$ be the twin of $V^+$. Let $x$ be the common vertex of $V^+$ and $V^-$.

We first claim that, if a tropic circle geometry $T$, say $T = T_{\mathfrak{f}}$ for $\mathfrak{f} \in \mathscr{E}$, contains some singular 3-space $Z$ of $V^+$ through $x$, then it contains $V^+$ entirely, and $\mathfrak{f} \cap V^-$ is a hyperbolic solid. Indeed, let $L$ be any line of $\Gamma$ in $V^-$ through $x$. Then, by the proof of Lemma 6.26, $L$ is collinear with all points of a 3-space $Z_L$ of $V^+$ containing $x$ and there is a pair of symplectic points $u_1, u_2$ in $Z \cap Z_L$. The set of points of $\Gamma$ collinear (in $\Gamma$) with $\{x, u_1, u_2\}$ is a plane $\pi$ (when endowed with all lines and hyperbolic lines in it) contained in $u_1 \Diamond u_2$ and $x$ lies on each line of $\Gamma$ contained in $\pi$. As each point of $\pi$ is collinear with $u_1$ and $u_2$, we have that $\pi$ is contained in $V^-$, and so $L$ is one of the lines through $x$ inside $\pi$. But, with respect to $\mathfrak{f}$ and $T$, the intersection $\beta(x) \cap \beta(u_1) \cap \beta(u_2)$ is easily seen to be a hyperbolic line $h$ of $\mathfrak{f}$, which must then be contained in $\pi$. Consequently, the point $L \cap h$ belongs to $\mathfrak{f}$. So we have shown that $\mathfrak{f} \cap V^-$ is a hyperbolic solid. But then every point of $V^+$ is collinear with at least two symplectic points of $\mathfrak{f}$, and hence belongs to $T$. The claim is proved.

Let $v$ be any point of $V^- \setminus \{x\}$. Let $\Pi \ni x$ be the 3-space of $V^+$ all of whose points are collinear with $v$. Then $\Pi$ is a singular 3-space of the quad $\Sigma(v)$. Since $\Sigma(v)$ is of type $\mathsf{D}_5$, $\Pi$ is contained in exactly two maximal singular subspaces of $\Sigma(v)$, one of which is the subspace $U(vx)$. The latter only contains points of $\Gamma$. Let $M_v$ be the other maximal singular subspace through $\Pi$. Since $v \notin \Pi$, no point of $M_v \setminus \Pi$ is collinear in $\Sigma(v)$ with $v$ and hence each point of $M_v \setminus \Pi$ is a new point. It follows that $M_v \setminus \Pi$ is the set of extended equator geometries containing $v$, such that $\Pi$ is contained in the associated tropic circle geometry. By our above claim, the latter is equivalent to requiring that $V^+$ is contained in the tropic circle geometry. Hence $M_v \subseteq M(V^+)$. The second assertion of our above claim now implies that $M(V^+)$ is the union of $V^+$ and all $M_w$, for $w$ ranging over $vx \setminus \{x\}$.

The set of points $v \in V^-$ for which a given new point of $M(V^+)$ belongs to $\Sigma(v)$ is a hyperplane of $V^-$, viewed as a projective 4-space. Also, the set of points $v \in V^-$ for which a given point $u \in \Gamma$ of $M(V^+)$ belongs to $\Sigma(v)$ is the hyperplane $V^- \cap u^\perp$. It follows that every triple of points of $M(V^+)$ is contained in at least one quad $\Sigma(v)$, $v \in V^-$, where they lie in a maximal subspace, and hence generate in $(\mathscr{P}, \mathscr{L})$ a plane. That plane is entirely contained in $M(V^+)$, as is easily checked. Hence $M(V^+)$ is a linear space where any three points generate a projective plane. By the celebrated theorem of Veblen and Young ([28], see Theorem 2.3 in [5] for a modern version), $M(V^+)$ is a projective space. Since every new line contains a unique ordinary point of $\Gamma$, $V^+$ is a hyperplane of $M(V^+)$ and so the dimension of $M(V^+)$ is equal to five. □

We denote the family of 5-spaces obtained in Lemma 6.28 by $\mathscr{M}$. It is clear that $M(V) \neq M(V')$ for distinct hyperbolic cones $V$, $V'$ (since $\Gamma \cap M(V) = V \neq V' = M(V') \cap \Gamma$); we refer to this property as *the injectivity of* $M(\cdot)$. Finally, we denote the family of projective planes contained as singular subspace in at least one quad by $\mathscr{T}$. Our next goal is to show that the 6-tuple $\mathfrak{E} = (\mathscr{P}, \mathscr{L}, \mathscr{T}, \mathscr{M}, \mathscr{U}, \mathscr{Q})$, with a suitable incidence relation, defines a geometry of type E$_6$, hence a building by [2]. We now define the incidence relation.

**Definition 6.29** (*Incidence relation in* $\mathfrak{E}$) The incidence between two elements of different types, where one or both of these elements are members of $\mathscr{P} \cup \mathscr{L} \cup \mathscr{T}$, or both belong to $\mathscr{U} \cup \mathscr{Q}$, is given by symmetrized set-theoretic (strict) containment. A 5-space (member of $\mathscr{M}$) is incident with a 4-space (member of $\mathscr{U}$) if they intersect in a singular 3-space. A 5-space (member of $\mathscr{M}$) is incident with a quad (member of $\mathscr{Q}$) if they intersect in a singular 4-space (which is in fact, with earlier terminology, a $4'$-space). Two elements of the same type are never incident. We will denote this incidence relation with $*$.

### 6.6 The symplectic polarity

Our eventual goal is to show that the pair $(\mathfrak{E}, *)$ is a geometry of type E$_6$. This proof will be facilitated by the use of the symplectic polarity $\theta$ that eventually will define $\Gamma$. We now define $\theta$.

**Definition 6.30** (*The Symplectic Polarity*) For $x \in \mathscr{P}$, we set $\theta(x) = \Sigma(x)$, and for $\Sigma \in \mathscr{Q}$, we define $\theta(\Sigma)$ as the point $x \in \mathscr{P}$ for which $\Sigma(x) = \Sigma$. Similarly, for $L \in \mathscr{L}$, we set $\theta(L) = U(L)$, and for $U \in \mathscr{U}$, we define $\theta(U)$ as the line $L \in \mathscr{L}$ for which $U(L) = U$. This is well defined by the injectivity of $U(\cdot)$. For $M \in \mathscr{M}$, we define $\theta(M) = M$ if $M$ is a symplecton of $\Gamma$. If $M$ can be written as $M(V^+)$, as in the second statement of Lemma 6.28, for appropriate $V^+$, then we define $\theta(M) = M(V^-)$, with $V^-$ the twin of $V^+$. This is well defined by the injectivity of $M(\cdot)$ and the fact that, by Lemma 6.26, twins are unique.

Note that in the previous definition we did not define $\theta$ on the planes. This will be done later, after Lemma 6.34 below. The next two lemmas will identify the two classes of maximal singular subspaces of the quads. The first lemma is an analog of Lemma 6.25.

**Lemma 6.31** *Two distinct members of $\mathscr{M}$ never have a subspace of dimension* 3 *in common.*

*Proof* By Fact 5.4, we may assume that at least one of the two members of $\mathscr{M}$ is of the form $M(V^+)$, with $V^+$ a hyperbolic cone with vertex $x$. Let the second member first be a symplecton $S$. Since both $S$ and $M(V^+)$ are singular subspaces of $(\mathscr{P}, \mathscr{L})$, their intersection is a projective subspace. Suppose that it contains a 3-space. Then, $S \cap V^+$ contains a plane $\pi$. Since no plane in $V^+$ consists solely of lines of $\Gamma$, the plane $\pi$ contains a hyperbolic line $h$. Hence $S = S(h)$. It follows that $x \in S \cap V^+$. Since $S$ does not contain planes without ordinary lines of $\Gamma$, and since every plane in $V^+$ not through $x$ is hyperbolic in $\Gamma$, we see that $S \cap V^+ = \pi$. But now $S \cap M(V^+) = S \cap V^+$ since $S$ does not contain any new point, a contradiction.

Now let the second member be $M(W^+)$, with $W^+$ a hyperbolic cone with vertex $y$. If $M(V^+) \cap M(W^+)$ has dimension at least 3, then $V^+ \cap M(W^+)$ has dimension at least 2. Since $V^+$ only contains points of $\Gamma$, we have $V^+ \cap M(W^+) = V^+ \cap W^+$. Hence we may assume for a contradiction that $V^+ \cap W^+$ contains a plane $\pi$, and that, if $V^+ \cap W^+ = \pi$, then $M(V^+) \cap M(W^+)$ contains at least one new point. Assume first that $x \in \pi$, and note that this is equivalent with $x = y$. There are two possibilities. The first one is that $V^+ \cap W^+ = \pi$. In this case some new point $\mathfrak{f}$ is contained in $M(V^+) \cap M(W^+)$. Then $T_\mathfrak{f}$ contains $V^+ \cup W^+$. Lemma 6.4(ii) yields $V^+ = W^+$, since both are cones with common vertex over a 3-space and only ordinary lines through the vertex. The second possibility is that $V^+ \cap W^+$ is a 3-space $Z \ni x$. Let $\mathcal{N}_x$ be as in Lemma 5.36 and let $\mathsf{D}_4(\mathcal{N}_x)$ be the corresponding geometry of type $\mathsf{D}_4$. Then $V^+$ and $W^+$ correspond to 3-spaces of $\mathsf{D}_4(\mathcal{N}_x)$ containing no line of the dual polar subspace of type $\mathsf{B}_3$ arising from the residue of $x$ in $\Gamma$. This means that, viewing $\mathsf{D}_4(\mathcal{N}_x)$ as a half spin geometry of type $\mathsf{D}_4$ of some quadric $Q$ of type $\mathsf{D}_4$, and the subspace of type $\mathsf{B}_3$ arising from the residue of $x$ in $\Gamma$ as the intersection with a subquadric $Q'$ of $Q$ of type $\mathsf{B}_3$, both $V^+$ and $W^+$ arise from the set of appropriate 3-spaces through points $v^+$ and $w^+$, respectively, of $Q \setminus Q'$ (the other possibility, namely that $V^+$ or $W^+$ would arise from the set of appropriate 3-spaces intersecting a given 3-space in planes is not feasible since in this case $V^+$ or $W^+$ would contain lines of the dual polar space related to $Q'$, a contradiction). But the set of appropriate 3-spaces through both $v^+$ and $w^+$ is either empty, or corresponds to a line in the half spin geometry, hence to $V^+ \cap W^+ = \emptyset$, or $V^+ \cap W^+$ a plane, both are contradictions to our assumption.

Assume now that $x \neq y$. Then $V^+ \cap W^+ = \pi$ is a hyperbolic plane. We claim that $x \perp y$. Indeed, let $a, b, c \in \pi$ be not on a common hyperbolic line. Then $c$ is close to $S(h(a,b))$ and both $x$ and $y$ must belong to the unique line of $S(h(a,b))$ consisting of points collinear with $c$. The claim follows. Again, there is some new point $\mathfrak{f}$ contained in $M(V^+) \cap M(W^+)$. Then $T_\mathfrak{f}$ contains $V^+ \cup W^+$. Lemma 6.4(ii) yields $xy \in V^+$, a contradiction. □

We denote the subset of elements of $\mathscr{U}$ incident with the quad $\Sigma$ by $\mathscr{U}(\Sigma)$ and we denote by $\mathscr{M}(\Sigma)$ the set of 4-spaces of $\Sigma$ obtained by intersecting $\Sigma$ with the members of $\mathscr{M}$ that are incident with $\Sigma$. The next lemma states that $\mathscr{U}(\Sigma)$ and $\mathscr{M}(\Sigma)$ are the two natural families of maximal singular subspaces of $\Sigma$.

**Lemma 6.32** *Let $\Sigma$ be an arbitrary quad. Then, every maximal singular subspace of $\Sigma$ is a member of exactly one of $\mathscr{U}(\Sigma), \mathscr{M}(\Sigma)$. Moreover, any two distinct elements of $\mathscr{U}(\Sigma)$ ($\mathscr{M}(\Sigma)$, respectively) intersect in either a point or a plane. Hence, $\mathscr{U}(\Sigma)$ and $\mathscr{M}(\Sigma)$ are the two natural systems of maximal singular subspaces of the hyperbolic quadric $\Sigma$; i.e. any 3-space of $\Sigma$ is contained in a unique member of $\mathscr{U}(\Sigma)$ and in a unique member of $\mathscr{M}(\Sigma)$.*

*Proof* We begin by showing that every maximal singular subspace $W$ of $\Sigma$ either belongs to $\mathscr{U}(\Sigma)$ or to $\mathscr{M}(\Sigma)$. Let $W_\Gamma$ be $W \cap \Gamma$. There are two cases, depending on the type of $\Sigma$. First we assume that $\Sigma = \Sigma(p)$, for some point $p$ of $\Gamma$. Then there are again two cases, depending on the dimension of $W_\Gamma$. Assume first that $W_\Gamma = W$. Clearly this is equivalent to $p \in W$. Using Lemma 5.36, we see that there are three possibilities for $W$. To identify these, we view $\mathsf{D}_4(\mathcal{N}_p)$ as a half spin $\mathsf{D}_4$-geometry, and we fix a polar space $Q$ of type $\mathsf{B}_3$ in the corresponding polar space $Q^+$ of type $\mathsf{D}_4$.

– *$W$ is a cone with vertex $p$ over a hyperbolic solid $V$.* (This case corresponds to the set of all maximal singular subspaces of $Q^+$ of the given half spin type through a fixed point of $Q^+ \setminus Q$.) Here, clearly $W = M(W) \cap \Sigma \in \mathscr{M}(\Sigma)$ (with the notation of Lemma 6.28).
– *$W$ is the intersection of $p^\perp$ with a symplecton $S$ through $p$.* (This case corresponds to the set of all maximal singular subspaces of $Q^+$ of the given half spin type through a fixed point of $Q$.) Here, clearly $W = S \cap \Sigma \in \mathscr{M}(\Sigma)$.
– *$W$ is a cone with some line $L$ through $p$ as vertex line and base a hyperbolic plane.* (This case corresponds to the set of all maximal singular subspaces of $Q^+$ of the given half spin type intersecting a fixed maximal singular subspace of $Q+$ of the other half spin type in a plane.) Here, clearly $W = U(L) \in \mathcal{U}(\Sigma)$.

Now assume that $p \notin W$. Then $W_\Gamma$ is a hyperplane section not containing $p$ of the maximal singular subspace of $\Sigma$ generated by $p$ and $W_\Gamma$. The three possibilities above give rise to the following three respective possibilities for $W_\Gamma$.

- $W_\Gamma$ *is a hyperbolic solid.* If $\mathfrak{f}$ is an arbitrary extended equator geometry containing $W_\Gamma$, then clearly $p = \beta(W_\Gamma)$. Hence, by Lemma 6.21, $U(W_\Gamma)$ is contained in $\Sigma$. Since $U(W_\Gamma)$ is the only maximal singular subspace of $\Sigma$ containing $W_\Gamma$ and not containing $p$, we conclude $W = U(W_\Gamma) = U(\langle p, \mathfrak{f}\rangle) \in \mathcal{U}(\Sigma)$.
- $W_\Gamma$ *can be written as* $p^\perp \cap q^\perp$, with $q$ symplectic to $p$. Put $h = h(p, q)$. Here, clearly all points of $U(h)$ belong to $\Sigma$. Since, as before, $U(h)$ is the only maximal singular subspace of $\Sigma$ containing $W_\Gamma$ and not containing $p$, we conclude $W = U(h) \in \mathcal{U}(\Sigma)$.
- $W_\Gamma$ *is a cone over a hyperbolic plane.* In $\mathsf{D}_4(\mathcal{N}_p)$, the planes are the intersections of two 3-spaces of different types; translated to the half spin setting, a plane is the set of 3-spaces of $Q^+$ of the given half spin type through a point $x$ of $Q^+$ and intersecting in an ordinary plane a 3-space of $Q^+$ of the other type incident with $x$. If the plane is hyperbolic, then $x \in Q^+ \setminus Q$, and so we see that $x$ defines a unique hyperbolic solid $V^+$ of $\mathsf{D}_4(\mathcal{N}_p)$ containing $W_\Gamma$ (we view $V^+$ as a cone over a hyperbolic solid, say with vertex $q$). Let $V^-$ be the twin of $V^+$, then clearly $p \in V^-$ (by the definition of twin). Let $W'$ be a hyperbolic solid in $V^-$ containing $p$, then, since each point of $V^+$ is collinear with at least two points of $W'$ (because $V^+$ is the twin of $V^-$) we see that $V^+ \subseteq T_\mathfrak{e}$, for every $\mathfrak{e} \in \mathcal{E}$ containing $W'$. Hence such $\mathfrak{e}$ belongs to $M(V^+)$, belongs to $\Sigma(p)$ and is collinear with all points of $W_\Gamma$. It follows that $M(V^+)$ intersects $\Sigma$ in a 4-space containing $W_\Gamma$ and not containing $p$; consequently $M(V^+) \cap \Sigma = W$ and $M(V^+) \in \mathcal{M}$.

The second case is $\Sigma = \Sigma(\mathfrak{e})$, with $\mathfrak{e} \in \mathcal{E}$. Since $\Sigma(\mathfrak{e}) \cap \Gamma = \mathfrak{e}$, the set $W_\Gamma \subseteq \mathfrak{e}$ is a hyperbolic solid of $\mathfrak{e}$. Put $x = \beta(W_\Gamma)$. We view $T_\mathfrak{e}$ as a half spin geometry of type $\mathsf{D}_5$ with corresponding hyperbolic quadric $Q^+$; the quadric of type $\mathsf{B}_4$ corresponding to $\mathfrak{e}$ is denoted by $Q$. Then $x \in T_\mathfrak{e}$ represents a maximal singular subspace of $Q^+$ of the given half spin type. That is already one maximal singular subspace $W_1$ of $\Sigma(\mathfrak{e})$ containing $W_\Gamma$. The second maximal singular subspace $W_2$ of $\Sigma(\mathfrak{e})$ containing $W_\Gamma$ is a maximal singular subspace of $Q^+$ of the type distinct from the given half spin type.

Suppose first that $W = W_1$. Besides the points of $\beta(x) = W_\Gamma$, the subspace $W$ further consists of the extended equator geometries $\mathfrak{f}$ containing $x$ and such that $T_\mathfrak{f}$ contains $W_\Gamma$. Since all these new points also belong to $\Sigma(x)$, it follows that $W \subseteq \Sigma(x) \cap \Sigma(\mathfrak{e}) = U(\langle x, \mathfrak{e}\rangle)$. Since $W$ is a 4-space, it follows that $W = U(\langle x, \mathfrak{e}\rangle) \in \mathcal{U}(\Sigma)$.

Suppose now that $W = W_2$. Then, $W$ is represented by a full pencil with centre $x$. A point of $Q^+$ belongs to it if, and only if, its residue (hence the set of maximal singular subspaces of the given half spin type incident with it) induces a half spin $\mathsf{D}_4$ in that full pencil. Such half spin $\mathsf{D}_4$s containing $x$ are given by the points of the full pencil collinear in $\Gamma$ with a point of $W_\Gamma$; those not containing $x$ are given by the hyperbolic $\mathsf{D}_4$s contained in the full pencil. If $\mathfrak{f} \in \mathcal{E}$ intersects the full pencil in a hyperbolic $\mathsf{D}_4$, then clearly $x \in T_\mathfrak{f}$. Hence, if we denote the cone with vertex $x$ and base $W_\Gamma$ by $V^+$, then it follows that $W$ is contained in $M(V^+) \in \mathcal{M}$.

So we have shown that every maximal singular subspace of $\Sigma$ either belongs to $\mathcal{U}$ (hence to $\mathcal{U}(\Sigma)$), or is contained in a unique member of $\mathcal{M}$ (and hence belongs to $\mathcal{M}(\Sigma)$). Now consider the graph $\Omega$ with vertices the maximal singular subspaces of $\Sigma$, two of those being adjacent if they meet in a 3-space. Then by Lemmas 6.25, 6.31, and the above, the sets $\mathcal{U}(\Sigma)$ and $\mathcal{M}(\Sigma)$ form a bipartition of $\Omega$. Since $\Omega$ is a connected bipartite graph with partitions the maximal singular subspaces of the two types, we see that $\mathcal{U}(\Sigma)$ is one of the natural types of maximal singular subspaces of $\Sigma$, and $\mathcal{M}(\Sigma)$ the other one. Hence, any two distinct elements of $\mathcal{U}(\Sigma)$ ($\mathcal{M}(\Sigma)$, respectively) intersect in either a point or a plane. This completes the proof of the lemma. □

This has the following consequence.

**Corollary 6.33** *Every $M \in \mathcal{M}$ sharing a 3-space with a given quad $\Sigma$ intersects $\Sigma$ in a 4-dimensional projective subspace of $\Sigma$.*

*Proof* Let $M \in \mathcal{M}$ share a 3-space $W$ with some quad $\Sigma$. By Lemma 6.32, there exists $Y \in \mathcal{M}(\Sigma)$ with $W \subseteq Y$. By the definition of $\mathcal{M}(\Sigma)$, there exists $M' \in \mathcal{M}$ with $Y \subseteq M'$. But then $W \subseteq M \cap M'$ implies, by Lemma 6.31, that $M = M'$. So $Y = M \cap \Sigma$ is a 4-dimensional subspace. □

We can now extend $\theta$ to $\mathcal{T}$.

**Lemma 6.34** *(i) For every point $p \in \mathcal{P}$ and every quad $\Sigma \in \mathcal{Q}$, we have $p * \Sigma$ if, and only if, $\theta(p) * \theta(\Sigma)$.*
*(ii) Let $\pi \in \mathcal{T}$. Then the intersection of all quads $\Sigma(x)$, where $x$ runs over the points of $\pi$, is a plane $\pi'$. Also, the intersection of all quads $\Sigma(x')$, where $x'$ runs over the points of $\pi'$, is precisely $\pi$ again.*

*Proof* Assertion (i) follows from the symmetry in the following three easy assertions. Let $p, q$ be points of $\Gamma$ and let $\mathfrak{e}, \mathfrak{f}$ be extended equator geometries of $\Gamma$.

(a) $p \in \Sigma(q)$ if, and only if, $p \perp q$;
(b) $p \in \Sigma(\mathfrak{e})$, if, and only if, $p \in \mathfrak{e}$ if, and only if, $\mathfrak{e} \in \Sigma(p)$;
(c) $\mathcal{E} \ni \mathfrak{e} \in \Sigma(\mathfrak{f})$ if, and only if, $T_\mathfrak{e} \cap \mathfrak{f}$ is a hyperbolic $\mathsf{D}_4$ (this is just part of Definition 6.18), which happens, by Lemma 6.13(i), if, and only if, $T_\mathfrak{f} \cap \mathfrak{e}$ is a hyperbolic $\mathsf{D}_4$.

We now prove (ii). By definition, $\pi$ is contained in at least one quad, say $\Sigma_1$. By Lemma 6.32, we can select two members $U_2, U_3$ of $\mathcal{U}$, contained in $\Sigma_1$, such that $\pi = U_2 \cap U_3$. By the definition of $\mathcal{U}$, we can select a quad $\Sigma_i \neq \Sigma_1$ such that $U_i = \Sigma_1 \cap \Sigma_i$, for $i = 2, 3$. So $\pi = \Sigma_1 \cap \Sigma_2 \cap \Sigma_3$. Now let $p \in \pi$ be arbitrary. Then $\Sigma(p)$ contains the point $\theta(\Sigma_i)$, $i = 1, 2, 3$, by (i). By Lemma 6.24, the points $\theta(\Sigma_1), \theta(\Sigma_2)$ and $\theta(\Sigma_3)$ are not contained in a member of $\mathcal{L}$, but since they are pairwise collinear, they generate a plane $\pi'$ inside $\Sigma(p)$. Again by Lemma 6.24, the intersection of all $\Sigma(p)$, with $p \in \pi$, is not a 4-space, and so it is at most a plane by Lemma 6.25. Since it contains $\pi'$, we see that $\pi'$ is the intersection of all $\Sigma(x)$, with $x \in \pi$. Since $\theta(\Sigma_i)$, $i = 1, 2, 3$, belongs to $\pi'$, and $\Sigma_i$ contains $\pi$, it follows that $\pi$ is the intersection of all $\Sigma(x')$, with $x' \in \pi'$. □

The plane $\pi'$ of the previous lemma is by definition the image of $\pi$ under $\theta$ and hence will be denoted by $\theta(\pi)$.

**Lemma 6.35** *The map $\theta$ preserves incidence and non-incidence and has order 2.*

*Proof* For two distinct elements $X$ and $Y$ of $\mathcal{P} \cup \mathcal{L} \cup \mathcal{T} \cup \mathcal{U} \cup \mathcal{Q} \cup \mathcal{M}$, we have to show that $X * Y$ if, and only if, $\theta(X) * \theta(Y)$. If $X \in \mathcal{P}$ and $Y \in \mathcal{L} \cup \mathcal{T} \cup \mathcal{Q}$, this follows from Lemmas 6.24 and 6.34. Now suppose $X \in \mathcal{P}$ and $Y \in \mathcal{U}$. Let $Y = U(L)$, with $L \in \mathcal{L}$. Then $X \in U(L)$ if, and only if, $X \in \Sigma(a) \cap \Sigma(b)$, for $a, b \in L$, $a \neq b$, if, and only if, $\Sigma(X)$ contains $a$ and $b$ if, and only if, $\Sigma(X)$ contains $L$ if, and only if, $\Sigma(X) = \theta(X)$ is incident with $L = \theta(Y)$.

Next suppose $X \in \mathscr{P}$ and $Y \in \mathscr{M}$. If $Y$ is a symplecton of $\Gamma$ and $X * Y$, then $X$ is a point of $\Gamma$ and $\Sigma(X)$ has the 4-space $X^\perp \cap Y$ in common with $Y = \theta(Y)$, hence $\theta(X) * \theta(Y)$. If $X$ is a point of $\Gamma$ and $X \notin Y$, then $\Sigma(X)$ has at most a line in common with $Y$ by Fact 5.5, and so $\theta(X)$ is not incident with $\theta(Y)$. Since for a new point $X \in \mathscr{E}$ the quad $\Sigma(X)$ does not contain ordinary lines of $\Gamma$, we see that $\theta(X)$ cannot be incident with $Y = \theta(Y)$. Now let $Y = M(V^+)$, with $V^+$ a hyperbolic cone with vertex $x$ over a hyperbolic solid. Let $V^-$ be the twin of $V^+$. Then $\theta(Y) = M(V^-)$ by definition. Suppose $X \in V^+$. If $X = x$, then $x^\perp$ contains $V^-$ and hence $\Sigma(x) \cap M(V^-) = V^-$, which means $\Sigma(x) * M(V^-)$. Now suppose $X \neq x$. Then $X^\perp \cap V^-$ is a 3-subspace, as is shown in the second paragraph of the proof of Lemma 6.26. Since by Lemma 6.32 every 3-subspace of a quad is contained in a member of $\mathscr{M}$, and members of $\mathscr{M}$ do not share any 3-space by Lemma 6.31, we see that $M(V^-) * \Sigma(X)$. Now let $X$ be a new point with $V^+ \subseteq T_X$. In the fourth paragraph of the proof of Lemma 6.28 it is shown that $X$, as an extended equator geometry, intersects $V^-$ in a hyperbolic 3-space, hence $\Sigma(X) * M(V^-)$ again. Now suppose $\theta(X) * M(V^+)$, for $X$ an ordinary point of $\Gamma$. Then $X^\perp$ intersects $V^+$ in at least a 3-space, and so $X$ belongs to $V^-$ by the definition of twin. If $X \in \mathscr{E}$ and $\Sigma(X) * M(V^+)$, then $X$ shares a hyperbolic solid with $V^+$, and by the definition of $V^-$, the latter is contained in $T_X$. So $X * M(V^-)$ again. Hence we have shown that $X * M(V^+)$ if, and only if, $\theta(X) * \theta(M(V^+))$, as $\theta(M(V^+)) = M(V^-)$.

So, if $X \in \mathscr{P}$, incidence and non-incidence is preserved. As $\theta(X) \in \mathscr{Q}$, the same holds for all $X \in \mathscr{Q}$. There remain the cases when $X, Y \in \mathscr{L} \cup \mathscr{T} \cup \mathscr{U} \cup \mathscr{M}$. First suppose $X \in \mathscr{L}$ and $Y \in \mathscr{T}$. Let $x, y \in X$ be two different points. Then $X * Y$ if, and only if, $X \subset Y$ if, and only if, $x \in Y$ and $y \in Y$ if, and only if, $\theta(x) * \theta(Y)$ and $\theta(y) * \theta(Y)$ if, and only if, $\theta(Y) \subseteq \Sigma(x) \cap \Sigma(y)$ if, and only if, $\theta(Y) \subseteq U(X)$ if, and only if, $\theta(Y) * \theta(X)$. Reading this from right to left also proves the assertion for $X \in \mathscr{U}$ and $Y \in \mathscr{T}$. Now suppose $X \in \mathscr{L}$ and $Y \in \mathscr{U}$, and let $x, y \in X$ be distinct again. Then, by definition, we can write $Y = U(L)$. So $X \subset Y$ if, and only if, $x, y \in Y$ if, and only if, $\Sigma(x) \cap \Sigma(y) \supseteq L$ if, and only if, $U(X) \supseteq L$ if, and only if, $\theta(X) * \theta(Y)$ (since $\theta(X) = U(X)$ and $\theta(Y) = L$).

So we may suppose $X \in \mathscr{M}$. First let $Y \in \mathscr{L}$. Let again $x, y \in Y$ be distinct points. Then we have that $X * Y$ if, and only if, $\Sigma(x) * \theta(X)$ and $\Sigma(y) * \theta(X)$ if, and only if, $\Sigma(x) \cap \Sigma(y) \cap \theta(X)$ is 3-dimensional (here we use Corollary 6.33 and the fact that no element of $\mathscr{U}$ is contained in an element of $\mathscr{M}$) if, and only if, $U(Y) \cap \theta(X)$ is 3-dimensional, and this is by definition equivalent to $\theta(Y) * \theta(X)$. Reading this the other way around takes care of the case $X \in \mathscr{M}$ and $Y \in \mathscr{U}$. So the remaining case is $X \in \mathscr{M}$ and $Y \in \mathscr{T}$. Now we choose three non-collinear points $x, y, z \in Y$. Then, if $Y \subseteq X$, then also $x, y, z \in X$ and hence $\Sigma(x) * \theta(X)$, $\Sigma(y) * \theta(X)$ and $\Sigma(z) * \theta(X)$. This implies $\Sigma(x) \cap \theta(X)$ is 4-dimensional, and likewise for $\Sigma(y) \cap \theta(X)$ and $\Sigma(z) \cap \theta(X)$. Hence $\Sigma(x) \cap \Sigma(y) \cap \Sigma(z) \cap \theta(X)$ is at least 2-dimensional. This is equivalent with $\theta(Y) \cap \theta(X)$ is 2-dimensional and so $\theta(Y) \subseteq \theta(X)$. By symmetry, the other direction also holds.

The fact that $\theta$ has order 2 is easy (by definition, it has order 2 on the points, and so $\theta^2$ must be the identity). □

We can now prove the main results.

### 6.7 Proofs of Theorems 1 and 2

The following result yields Theorem 1.

**Theorem 6.36** *(i) The geometry $(\mathfrak{E}, *)$ has type* E$_6$, *where, with the Bourbaki labeling specified in Sect. 2 (see Fig. 1), the sets $\mathscr{P}, \mathscr{M}, \mathscr{L}, \mathscr{T}, \mathscr{U}, \mathscr{Q}$ are the elements of types 1, 2, 3, 4, 5, 6, respectively.*

*(ii) The map $\theta$ is a symplectic polarity of $(\mathfrak{E}, *)$ with fixed point structure $\Gamma$. More exactly, the set of absolute points and absolute lines are precisely the set of points and lines of $\Gamma$, and the fixed planes and fixed 5-spaces are the planes and symplecta, respectively, of $\Gamma$.*

*Proof* It follows from the definition and Lemma 6.32 that the elements incident with a quad form the oriflamme complex of a polar space of type $\mathsf{D}_5$ with point set the elements of $\mathscr{P}$ in the quad, line set the elements of $\mathscr{L}$ in the quad, plane set the elements of $\mathscr{T}$ in the quad, maximal singular subspaces of one type the elements of $\mathscr{U}$ in the quad, and the maximal singular subspaces of the other type the elements of $\mathscr{M}$ incident with the quad. Applying the map $\theta$, the elements of $\mathscr{L} \cup \mathscr{T} \cup \mathscr{M} \cup \mathscr{U} \cup \mathscr{Q}$ incident with a fixed point also form the oriflamme complex of a polar space of type $\mathsf{D}_5$. This determines all rank 2 residues of $(\mathfrak{E}, *)$ except for the residue of cotype $\{1, 6\}$. But it is easy to see that, if a point is incident with a plane, and that plane is incident with a quad, then the point is incident with the quad. Hence the diagram of $(\mathfrak{E}, *)$ has type $\mathsf{E}_6$. This proves *(i)*.

By [2], we now know that $(\mathfrak{E}, *)$ is a building of type $\mathsf{E}_6$. Since $\theta$ has order 2 it is bijective, and since it preserves incidence, it is an automorphism of the building. Since it does not preserve the types of the elements, it is a polarity. If $p$ is a point of $\Gamma$, then clearly $p * \theta(p) = \Sigma(p)$. If $\mathfrak{e} \in \mathscr{E}$, then we claim that no point of $(\mathfrak{E}, *)$ collinear with $\mathfrak{e}$ is incident with $\theta(\mathfrak{e})$. Indeed, if $\mathfrak{f} \in \mathscr{E}$ is collinear with $\mathfrak{e}$, and $\mathfrak{f}$ is incident with $\theta(\mathfrak{e})$, then, as extended equator geometries, they share a hyperbolic solid. But since $T_\mathfrak{f} \cap \mathfrak{e}$ is a geometric hyperplane of $\mathfrak{e}$ by Lemma 6.7 and Definitions 6.18 and 6.29, it intersects $\mathfrak{e} \cap \mathfrak{f}$ nontrivially, a contradiction. Also, an ordinary point of $\mathfrak{e}$ never lies in $T_\mathfrak{e}$, hence the claim follows. By Proposition 4.14, $\theta$ is a symplectic polarity. Clearly, the set of absolute points is the point set of $\Gamma$, the set of absolute lines is the line set of $\Gamma$, and the set of fixed planes is the plane set of $\Gamma$. The latter is true since every point of an absolute plane must be absolute (otherwise the plane cannot coincide with its image), and taking into account the definition of $\theta$ on points, it follows that all points of the plane must be collinear in $\Gamma$ with each other. Hence the plane is an "ordinary" plane of $\Gamma$ (and vice versa). Also, the fixed 5-spaces are the symplecta of $\Gamma$, as is immediate from the definition of $\theta$ on $\mathscr{M}$.

This completes the proof of the theorem. □

This now implies the following uniqueness result and Theorem 2.

**Theorem 6.37** *A building $\Delta$ of type $\mathsf{E}_6$ admits, up to conjugacy, a unique symplectic polarity.*

*Proof* Let $\mathbb{K}$ be the field underlying $\Delta$, i.e., each residue plane of $\Delta$ is isomorphic to the Pappian projective plane defined over $\mathbb{K}$ (and by the classification of spherical buildings of rank at least 3, see [22], this defines $\Delta$ in a unique way). Let $\Gamma$ be the split building of type $\mathsf{F}_4$ over $\mathbb{K}$, i.e., the planes of $\Gamma$ are the Pappian projective planes defined over $\mathbb{K}$ and the rank 2 residue of type $\mathsf{C}_2$ corresponds to a symplectic generalized quadrangle over $\mathbb{K}$. Then Theorem 6.36 implies that $\Gamma$ can be seen as the fixed point structure of a symplectic polarity of a building of type $\mathsf{E}_6$ over $\mathbb{K}$, hence isomorphic to $\Delta$. This shows existence of the symplectic polarity in $\Delta$.

Let now $\theta$ and $\theta'$ be two symplectic polarities in $\Delta$. Then their fixed point structures, say $\Gamma$ and $\Gamma'$, respectively, are both split buildings of type $\mathsf{F}_4$ over $\mathbb{K}$. Hence, by Theorem 10.2 of [22], $\Gamma$ and $\Gamma'$ are isomorphic to each other, say, under the isomorphism $\psi : \Gamma \to \Gamma'$. Also, by Lemma 4.18, the absolute points of $\theta$ and $\theta'$ form geometric hyperplanes $H$ and $H'$, respectively, of the natural point-line geometry associated with $\Delta$.

Let $Q$ be a quad of $\Delta$. Then, by Lemma 4.20, either it is absolute under $\theta$ and intersects $H$ in a singular geometric hyperplane of $Q$, or it is not absolute and intersects $H$ in a parabolic quadric of type B$_4$. Suppose first $Q$ is absolute and put $x = Q^\theta$. Then, by Lemma 4.19, $Q \cap H$ consists of the union of lines of $\Gamma$ through $x$. Suppose now that $Q$ is not absolute, and again put $x = Q^\theta$. We claim that $Q \cap H$ does not contain lines of $\Gamma$. Indeed, let, for a contradiction, $L$ be a line of $\Gamma$ in $Q \cap H$. Pick $z \in L$ arbitrarily. Since $L$ is absolute, $z^\theta$ contains $L$. But $z^\theta$ also contains $x$. Since $z^\theta$ is a polar space, there is some point on $L$ collinear with $x$. Hence, by Fact 4.3, $x$ neighbours $Q$. But this contradicts the second statement of Proposition 4.14. The claim is proved. Consequently, by Proposition 4.15, all the lines of $Q$ in $H$ are hyperbolic lines of $\Gamma$. Let $p, q$ be two points of $Q \cap H$ not collinear in $\Delta$. Then, by the foregoing, $p$ and $q$ are either symplectic, special or opposite. By Corollary 4.4, if $p$ and $q$ are special, then $p \bowtie q$ also belongs to $Q \cap H$ and hence also the line of $\Gamma$ joining $p$ and $p \bowtie q$, a contradiction. If $p$ and $q$ are symplectic, then, again by Corollary 4.4, the set $p^\perp \cap q^\perp$ ($\perp$ denotes collinearity in $\Gamma$) belongs to $Q$, a contradiction as this set contains lines of $\Gamma$. Hence $p$ and $q$ are opposite. Now, using Corollary 4.4, it is easy to see that $\widehat{E}(p, q) \subseteq \Gamma$ is contained in $Q \cap H$. But since $\widehat{E}(p, q)$ is also convex with respect to the hyperbolic lines of $\Gamma$, we easily deduce that $\widehat{E}(p, q) = Q \cap H$. Conversely, every extended equator geometry $\widehat{E}(a, b)$ arises as the intersection of a non-absolute quad with $H$. Indeed, the points $a$ and $b$ are not collinear in $\Delta$ and the unique quad $Q$ containing $a$ and $b$ cannot be absolute as otherwise $Q^\theta$ is collinear with both $a$ and $b$, contradicting Lemma 4.20. Hence there is a natural bijective correspondence between the family of quads of $\Delta$ and the union of the family of extended equator geometries of $\Gamma$ and the family of point-perps in $\Gamma$.

The same thing holds for $\theta'$. Hence $\psi$ induces a natural permutation $\xi$ of the family of quads of $\Delta$ by first intersecting a quad $Q$ with $H$, then applying $\psi$, and then taking the unique quad of $\Delta$ whose intersection with $H'$ coincides with $(H \cap Q)^\psi$. Let $Q$ and $Q^*$ be two quads. Suppose they share a 4-space. Then clearly the intersection $Q \cap Q^* \cap H$ contains a 3-space, and so does $(Q \cap Q^* \cap H)^\psi$. Hence $Q^\xi \cap (Q^*)^\xi \cap H'$ contains a 3-space and so $Q^\xi$ and $(Q^*)^\xi$ share a 4-space. A similar argument, or just arguing with the inverse of $\psi$, implies that, if, for two quads $Q, Q'$ of $\Delta$, the quads $Q^\xi$ and $(Q^*)^\xi$ intersect in a 4-space, then $Q$ and $Q'$ intersect in a 4-space.

This now means that $\xi$ induces an isomorphism, which we can again denote by $\xi$, from $\Delta$ to $\Delta'$, since "collinearity" of quads is preserved in both directions. Clearly, $\xi$ extends $\psi$. Hence $\theta$ and $\theta'$ are conjugate (more exactly, $\theta = \xi^{-1}\theta'\xi$). □

## 7 Proof of Theorem 3

Let $\Gamma$ be a split building of type F$_4$, viewed as a symplectic metasymplectic parapolar space, defined over the field $\mathbb{K}$. Suppose that it is point-line-embedded in the natural point-line geometry associated with a building $\Delta$ of type E$_6$. We identify the points and lines of $\Gamma$ with points and lines of $\Delta$. If two points $x, y \in \Gamma$ are collinear in $\Gamma$, then we say that $x$ and $y$ are $\Gamma$-*collinear*, and write $x \perp_\Gamma y$. Likewise, if two points $x, y \in \Delta$ are collinear in $\Delta$, then we say that they are $\Delta$-*collinear*, and write $x \perp_\Delta y$. The set of points equal to or $\Gamma$-collinear ($\Delta$-collinear, respectively) to a point $x$ of $\Gamma$ (of $\Delta$, respectively) is denoted by $x^{\perp_\Gamma}$ ($x^{\perp_\Delta}$, respectively).

### 7.1 The type in $\Delta$ of the symplecta of $\Gamma$

We now assign to every type of elements of $\Gamma$ a unique type of elements of $\Delta$. This will be trivial for the points and lines of $\Gamma$, and easy for the planes. It requires more work for the symplecta. We will show that every symplecton of $\Gamma$ is contained in a unique 5-space of $\Delta$.

We start with an easy lemma, showing that the planes of $\Gamma$ correspond to (full) planes of $\Delta$.

**Lemma 7.1** *Let $\pi$ be a plane of $\Gamma$. Then the points and lines of $\pi$ are all the points and lines of a unique plane of $\Delta$, i.e., there is a unique element of type 4 of $\Delta$ whose points and lines are precisely those of $\pi$.*

*Proof* This follows immediately from the fullness of the embedding. ∎

An immediate consequence of this lemma is the following.

**Corollary 7.2** *The defining field of $\Delta$ is $\mathbb{K}$.*

*Proof* By Lemma 7.1, the planes of $\Gamma$ and $\Delta$ are isomorphic, hence they are defined over the same field, which is $\mathbb{K}$ by our choice of $\Gamma$. ∎

We now determine which elements of $\Delta$ the symplecta of $\Gamma$ correspond to. If a symplecton $S$ of $\Gamma$ is contained in a singular subspace of $\Delta$, then, since such a subspace has dimension at most 5, $S$ is embedded in a 5-space in the usual way, i.e., as the absolute geometry of a symplectic polarity. In this case, $S$ corresponds to that unique 5-space of $\Delta$ and the point set of $S$ coincides with the point set of that 5-space. We next consider the case when some symplecton is not embedded in a 5-space of $\Delta$.

**Lemma 7.3** *If the point set of a symplecton $S$ of $\Gamma$ does not coincide with the point set of a 5-space of $\Delta$, then $S$ is contained in a unique quad of $\Delta$. Moreover, in this case, the characteristic of the field $\mathbb{K}$ is equal to 2.*

*Proof* Since $S$ is not contained in a singular subspace of $\Delta$, there are two points $x, y$ of $S$ that are not $\Delta$-collinear. Then they are not $\Gamma$-collinear either, and so $T := x^{\perp_\Gamma} \cap y^{\perp_\Gamma}$ is a symplectic quadrangle (i.e., a symplectic polar space of rank 2). By Fact 4.1, there is a unique quad $Q$ containing $x, y$. By Corollary 4.4, $Q$ also contains $T$. Hence $T \subseteq Q_{x,y} := x^{\perp_\Delta} \cap y^{\perp_\Delta} \subseteq Q$, and $Q_{x,y}$ is a polar space of type $\mathsf{D}_4$.

Suppose for a contradiction that some point $p \in S$ is not contained in $Q$. Note that, by the previous paragraph, $p$ is not $\Gamma$-collinear with either $x$ or $y$. There are two possibilities. The first one is that $(p^{\perp_\Gamma} \cap x^{\perp_\Gamma}) \setminus T \neq \emptyset$. Let $q \in (p^{\perp_\Gamma} \cap x^{\perp_\Gamma}) \setminus T$. Then $q$ belongs to $(x^{\perp_\Gamma} \cap S) \setminus T$ and is $\Gamma$-collinear with $p$. Also, the line $pq$ intersects $y^{\perp_\Gamma}$ in a point $z$ distinct from $q$. Then $Q$ contains $q$ and $z$ and, since $Q$ is a subspace, it contains all points of the line $qz$. This contradicts our assumption that $p$ is not contained in $Q$. The second possibility is that $p^{\perp_\Gamma} \cap x^{\perp_\Gamma} \subseteq T$. Then, since every line of $S$ through $x$ contains a point $\Gamma$-collinear with $p$ and one $\Gamma$-collinear with $y$, we have $p^{\perp_\Gamma} \cap x^{\perp_\Gamma} = T$ and so $T \subseteq p^{\perp_\Gamma}$. But then, by the foregoing, every point of every line of $S$ through $p$ is contained in $Q$, and hence so is $p$, a contradiction. Hence $S \subseteq Q$.

Now for $T$ there are two possibilities. First assume that $T$ is contained in a singular subspace $W$ of $Q_{x,y}$. Since $T$ is a symplectic generalized quadrangle over $\mathbb{K}$, and $T$ is fully embedded in $W$, we see that $W$ has dimension 3 and $T$ is embedded in the standard way as the absolute geometry of a symplectic polarity. Since $S$ is a symplectic polar space of rank 3, it contains a point $z \notin \{x, y\}$ which is $\Gamma$-collinear with all points of $T$. If we denote by $\perp_S$ the collinearity in $S$ (inherited from $\Gamma$), then $x^{\perp_S}$, $y^{\perp_S}$ and $z^{\perp_S}$ are three singular 4-dimensional subspaces of $Q$ containing $W$, a contradiction since $Q$ is a polar space of type $\mathsf{D}_5$.

Hence there are two points $a, b \in T$ which are not $\Delta$-collinear. Then the lines $ax$ and $by$ are opposite in both polar spaces $S$ and $Q$. In $Q$, the subspace generated by $ax$ and $by$ is a grid. Since $S$ is also a polar space and since $\Gamma$-collinearity implies $\Delta$-collinearity, it also contains that grid. This implies that the characteristic of $\mathbb{K}$ is equal to 2. Indeed, this either follows from a direct calculation, or we can geometrically argue as follows. A grid in $S$ is contained in a 3-space and hence it is the intersection of the perp of two non-collinear points of $S$, implying that it is contained in a symplectic quadrangle $S'$ over $\mathbb{K}$. So $S'$ contains three mutually non-intersecting lines which are concurrent to three other mutually non-intersecting lines. Dually, this yields three mutually non-collinear points in a parabolic polar space $\Omega$ of type $\mathsf{B}_2$ over $\mathbb{K}$ collinear with three other mutually non-collinear points. Noting that, if the characteristic of $\mathbb{K}$ is not equal to 2, then $\Omega$ arises from a polarity $\rho$ in $\mathsf{PG}(4, \mathbb{K})$, we find two planes (generated by the respective non-collinear points) contained in each other's image under $\rho$, a contradiction to the dimensions.

The lemma is proved. □

Our next aim is to show that, in the case when $\mathsf{char}\mathbb{K} = 2$ also, the symplecta correspond to 5-spaces of $\Delta$. Hence we aim for the following proposition.

**Proposition 7.4** *The point set of every symplecton of $\Gamma$ coincides with the point set of some 5-space of $\Gamma$.*

The way to accomplish this in the case when $\mathsf{char}\mathbb{K} = 2$, is to use the following strategy. The rank 2 residues of $\Gamma$ isomorphic to projective planes and consisting of symplecta and planes of $\Gamma$ can be interpreted in $\Delta$ as representations of projective planes in projective spaces where the points of the projective plane are lines of the projective space, and the lines of the projective plane are reguli and partial planar line pencils (see below for precise definitions). In the next four lemmas, we prepare to prove nonexistence of those representations that would arise if some symplecton of $\Gamma$ would not correspond with a 5-space of $\Delta$.

**Lemma 7.5** *Let $p$ be a point of the projective plane $\mathsf{PG}(2, \mathbb{K})$, with $\mathbb{K}$ a field of characteristic 2. Let $\mathscr{C}$ be a family of conics with nucleus $p$ satisfying each of the following properties.*

*(i) All conics of $\mathscr{C}$ have the same set of tangents.*
*(ii) Every pair of conics in $\mathscr{C}$ intersects in a unique point.*
*(iii) Every pair of points on distinct tangents lies on a unique member of $\mathscr{C}$.*

*Then $\mathbb{K}$ is a perfect field.*

*Proof* Let $p$ have coordinates $(0, 0, 1)$ and let one of the conics $C_0$ in $\mathscr{C}$ have equation $XY = Z^2$. An arbitrary irreducible conic with nucleus $p$ has equation $XY = aX^2 + bY^2 + cZ^2$, $a, b, c \in \mathbb{K}, c \neq 0$. This conic belongs to $\mathscr{C}$ only if it meets every tangent to $C_0$ (see $(i)$). Such a tangent is a line with equation $X = 0$, or $Y = 0$, or $X = sY$, with $s$ a nonzero square in $\mathbb{K}$. Hence, in this case, $bc \in \mathbb{K}^2$, $ac \in \mathbb{K}^2$ and $sY^2 = as^2Y^2 + bY^2 + cZ^2$ has a unique solution with $Y = 1$. Thus, $cs + acs^2 + cb$ is a square, and as $bc$ and $ac$ are squares, this happens if, and only if, $cs$, and hence $c$, is a square. Hence $a, b, c$ are squares and we can replace $a, b, c$ by $a^2, b^2, c^2$ and obtain as equation of the conic $XY = (aX + bY + cZ)^2$. Moreover, in this case we automatically have that every line with equation $X = nY$, with $n$ a non square of $\mathbb{K}$, does not intersect the conic. The unique point of intersection $(x, y, z)$ of this conic with $C_0$ (see $(ii)$) satisfies $xy = z^2$ and $ax + by = (c + 1)z$. If $c + 1 \neq 0$, then $xy = s_1^2 x^2 + s_2^2 y^2$, with $s_1, s_2 \in \mathbb{K}^2$ and, without loss, $s_1 \neq 0$. However, then $(x + s_1^{-2} y, y, z + s_1^{-1} y)$ is also a point of the intersection, contradicting $(ii)$. Hence $c = 1$. But then $ax + by = 0$ and, since $a \neq 0$ without loss, $by^2 = az^2$, implying $ab \in \mathbb{K}^2$.

Suppose now for a contradiction that $\mathbb{K}$ is not perfect.

We select two points $(1, 0, x)$ and $(0, 1, y)$, with $x \in \mathbb{K}^2$ and $y \in \mathbb{K}\setminus\mathbb{K}^2$, on the distinct tangent lines with equations $Y = 0$ and $X = 0$, respectively. The only conic with equation $XY = (aX + bY + Z)^2$ containing these points satisfies $a = x$ and $b = y$, and so has to belong to $\mathscr{C}$ (see $(iii)$), contradicting $\mathbb{K}^2 \ni ab = xy \notin \mathbb{K}^2$. □

A *regulus of lines* in a 3-dimensional projective space over the field $\mathbb{K}$ is the set of generators of one type of a hyperbolic quadric. The *complementary regulus* is the set of generators of the other type. An elementary (and easy to prove) property is that every plane containing an element of a regulus also contains a unique element of the complementary regulus. We will frequently use this property without notice.

**Lemma 7.6** *Let $\Pi$ be a projective plane and $\mathsf{PG}(4, \mathbb{K})$ a projective 4-space. Suppose that each point of $\Pi$ is identified with a line of $\mathsf{PG}(4, \mathbb{K})$ and that the set of points of each line of $\Pi$ corresponds bijectively to the set of lines of a regulus of lines in some 3-space of $\mathsf{PG}(4, \mathbb{K})$. Then all such 3-spaces coincide, and hence all lines of $\mathsf{PG}(4, \mathbb{K})$ corresponding to points of $\Pi$ are contained in some common 3-space $\mathsf{PG}(3, \mathbb{K})$.*

*Proof* For ease of notation we will call the lines of $\mathsf{PG}(4, \mathbb{K})$ which are identified with points of $\Pi$ briefly $\Pi$-*lines*, the reguli corresponding to the lines of $\Pi$ are briefly called $\Pi$-*reguli*.

Suppose at least two $\Pi$-reguli are contained in a common 3-space. Since any $\Pi$-regulus is contained in the 3-space generated by any pair of its $\Pi$-lines, and every $\Pi$-line is contained in a $\Pi$-regulus sharing a $\Pi$-line with each of the two given $\Pi$-reguli in the 3-space, it is clear that all $\Pi$-reguli are contained in that common 3-space. Hence we may assume for a contradiction that distinct $\Pi$-reguli generate distinct 3-spaces of $\mathsf{PG}(4, \mathbb{K})$. Let $\Sigma$ be one such 3-space, and let $R$ be the corresponding regulus. No $\Pi$-line not belonging to $R$ is contained in $\Sigma$, hence every such $\Pi$-line $E$ intersects $\Sigma$ in a unique point $p_E$, and this point is not on any line of $R$. Fix such $E$ and $p_E$. Consider an arbitrary line $L$ of $R$. Let $K$ be the unique line of the complementary regulus of the $\Pi$-regulus $R^*$ containing $L$ and $E$ passing through $p_E$. Then $K$ intersects $L$, say in the point $x$, and hence $K$ is contained in $\Sigma$. Clearly, $K$ does not intersect any member of $R\setminus\{L\}$, and hence $K$ is a tangent line to the quadric $Q$ defined by $R$. Now, since $\Pi$ is a projective plane, every $\Pi$-line distinct from $E$ is contained in a $\Pi$-regulus containing $E$ and a line of $R$. Hence we conclude that the set $I$ of intersection points with $\Sigma$ of the $\Pi$-lines not in $R$ coincides with the set of points on the tangent lines of $Q$ through $p_E$ except for the points of $Q$ itself. If $\mathsf{char}\mathbb{K} \neq 2$, then the hyperbolic quadric $Q$ defined by $R$ is the absolute geometry of an orthogonal polarity and so this set $I$ is a cone with vertex $p_E$ where one conic (corresponding to the intersection of that cone with $Q$) is removed; since $p_E$ was arbitrary, and since that cone has a unique vertex, this is a contradiction. If $\mathsf{char}\mathbb{K} = 2$, then $Q$ is a subset of the set of absolute points of some symplectic polarity and hence $I$ is the set of points of a plane $\pi$ with one conic removed; this conic has nucleus $p_E$ since every line through $p_E$ in $\pi$ must have a unique point not contained in $I$. Again, since $p_E$ was arbitrary, every other point is also the nucleus of the conic, a contradiction. □

The next lemma shows that the hypothesis in Lemma 7.6 is untenable.

**Lemma 7.7** *Let $\Pi$ be a projective plane and $\mathsf{PG}(3, \mathbb{K})$ a projective 3-space. Suppose each point of $\Pi$ is associated with a line of $\mathsf{PG}(3, \mathbb{K})$. Then it is impossible that the set of points of each line of $\Pi$ corresponds bijectively to the lines of a regulus of lines in $\mathsf{PG}(3, \mathbb{K})$.*

*Proof* This time $\mathbb{K}$ cannot be finite as there are not enough mutually disjoint lines in $\mathsf{PG}(3, \mathbb{K})$ (if $|\mathbb{K}| = q$, then $\mathsf{PG}(3, \mathbb{K})$ has at most $\frac{q^3+q^2+q+1}{q+1} = q^2+1$ mutually disjoint lines, whereas we need $q^2 + q + 1$ such lines if $\Pi$ is of order $q$; note that $\Pi$ necessarily has order $q$ as each regulus contains $q + 1$ lines). We use the same terminology as in the previous proof.

Let $\alpha$ be a plane of $\mathsf{PG}(3, \mathbb{K})$ through some $\Pi$-line $L$. We intersect every other $\Pi$-line with $\alpha$ and obtain a set of points $\mathscr{A}$ with the following properties.

(1) Each $\Pi$-regulus through $L$ is mapped onto a line of $\alpha$ distinct from $L$ (i.e., the transversal to the $\Pi$-regulus in the plane $\alpha$); since these lines cannot have intersection points off $L$, all these lines intersect in a fixed point $p \in L$. We call such a line a $* - \Pi$-line.
(2) The $\Pi$-reguli not through $L$ correspond to conics of $\alpha$ completely contained in $\mathscr{A}$ and each $* - \Pi$-line intersects each such conic in a unique point; hence $p$ is the nucleus of each such conic. We call such a conic a $* - \Pi$-conic.
(3) Each pair of $* - \Pi$-conics intersect in a unique point.
(4) Through each pair of points on distinct $* - \Pi$-lines passes a unique $* - \Pi$-conic.

Since nuclei exist by (2), the field $\mathbb{K}$ has characteristic 2, and as $L$ is a line through the nucleus external to every $* - \Pi$-conic, $\mathbb{K}$ is not perfect. This contradicts Lemma 7.5. □

**Lemma 7.8** *Let $\Pi$ be a linear space and $\mathsf{PG}(4, \mathbb{K})$ a projective 4-space. Suppose each point of $\Pi$ is identified with a line of $\mathsf{PG}(4, \mathbb{K})$ and suppose that the set of points of any line of $\Pi$ corresponds bijectively either to the set of lines of a regulus of lines in some 3-space of $\mathsf{PG}(4, \mathbb{K})$, or to the set of lines of a partial planar line pencil (i.e., a set of concurrent coplanar lines of size at least 3). If the elements of $\Pi$ span $\mathsf{PG}(4, \mathbb{K})$, then $\Pi$ is not a projective plane.*

*Proof* Suppose, by way of contradiction, that $\Pi$ is a projective plane. We call a partial planar line pencil corresponding to a line of $\Pi$ briefly a $\Pi$-pencil, and use the terminology of the previous proofs concerning $\Pi$-lines and $\Pi$-reguli. The *vertex* of a $\Pi$-pencil is the common point of its members.

Note that by Lemma 7.6, there is at least one $\Pi$-pencil. Also, if no $\Pi$-reguli exist, then the $\Pi$-lines are either contained in a plane of $\mathsf{PG}(4, \mathbb{K})$ (if there exist two $\Pi$-pencils with distinct vertices), or in a solid of $\mathsf{PG}(4, \mathbb{K})$ (if all $\Pi$-pencils have the same vertex), a contradiction.

So there exists a $\Pi$-regulus $R$ and a $\Pi$-pencil $P$, sharing a line $L$. Let $x$ be the vertex of $P$. Clearly, if $P$ is contained in the 3-space $\Sigma$ generated by $R$, then all $\Pi$-lines are contained in $\Sigma$. Hence, by hypothesis, we may assume that $P$ intersects $\Sigma$ in $L$. It follows that every member of $P\setminus\{L\}$ is skew to every member of $R\setminus\{L\}$. Consider two distinct members $M_1, M_2 \in P\setminus\{L\}$ and two distinct members $K_1, K_2 \in R\setminus\{L\}$. Let $L'$ be the unique line of the complementary regulus to $R$ through $x$. Then both $K_1$ and $K_2$ intersect $L'$ and so the 3-spaces generated by $K_1$ and $M_1$, and by $K_2$ and $M_2$ intersect in a plane containing the line $L'$. It follows that $L'$ intersects the $\Pi$-line $K$ corresponding to the point of $\Pi$ associated with the line common to the two $\Pi$-reguli determined by $K_1, M_1$ and $K_2, M_2$. Hence $K$ intersects a member $M$ of $R$ on the line $L'$, but not in $x$. Consequently, the line of $\Pi$ through the points of $\Pi$ associated with $K$ and $M$ corresponds to a $\Pi$-pencil, which has no element in common with $P$. This contradicts the fact that $\Pi$ is a projective plane. □

In the next lemma we exclude the case where symplecta correspond to both quads and 5-spaces.

**Lemma 7.9** *Suppose that some symplecton $S$ of $\Gamma$ is contained in a quad $Q$ of $\Delta$. Then, each symplecton of $\Gamma$ is contained in a quad of $\Delta$.*

*Proof* Assume for a contradiction that there is a symplecton $S'$ contained in a 5-space $W$ (and note that $S \neq S'$ as $S'$ cannot be contained in $W \cap Q$). Since, by Fact 5.5, the graph on the symplecta of $\Gamma$ with adjacency "meeting in a plane" is connected, we may assume that $S$ and $S'$ share a plane $\pi$. This implies that $Q$ and $W$ are incident in $\Delta$ (this is the dual of the fact that, if a point is incident with a plane, and a 5-space is incident with that plane, then that point and 5-space are mutually incident) and hence intersect in a $4'$-space $U$. Let $x \in U$ be the image of $U$ with respect to the symplectic polarity $\rho$ induced on $W$ by $S'$. Since $\pi \subseteq U$, it follows that $x = U^\rho \subseteq \pi^\rho = \pi$. Let $L$ be a line in $\pi$ not containing $x$ and let $p \in L$. The residue in $\Gamma$ of $\{p, L\}$ is a projective plane $\Pi'$, with point set the set of planes of $\Gamma$ through $L$ and line set the symplecta of $\Gamma$ through $L$. Now, $\Pi'$ is contained in the residue $\mathscr{R}$ of the flag $\{p, L\}$ in $\Delta$. That residue is a projective 4-space $\mathsf{PG}(4, \mathbb{K})$ with point set the 5-spaces through $L$, line set the planes through $L$, plane set the set of 4-spaces through $L$ and set of hyperplanes the quads through $L$.

The planes of $\Gamma$ in $S'$ through $L$ are contained in a common 3-space $V$ which is not contained in $Q$, by our choice of $L \not\ni x$. Since $V$ can be regarded as a flag consisting of the 5-space $W$ and a 4-space $U'$, we see that the set of planes of $\Gamma$ through $L$ in $S'$ corresponds to a set $P$ of the lines of $\mathsf{PG}(4, \mathbb{K})$ through a point (corresponding to $W$) inside a plane (corresponding to $U'$). Hence the lines of $\Pi'$ that correspond to 5-spaces form complete planar line pencils of $\mathsf{PG}(4, \mathbb{K})$ (i.e., the set of *all* lines through a given point inside a given plane). We now turn to $Q$.

In $Q$, the symplecton $S$ is embedded as an orthogonal polar space. The latter contains subspaces isomorphic to parabolic quadrics of type $\mathsf{B}_3$. Since, by Lemma 7.3, the characteristic of $\mathbb{K}$ is equal to 2, these are obtained automatically if we consider the subbuilding $\Gamma'$ of $\Gamma$ defined from $\Gamma$ by viewing $\Gamma$ as a building of mixed type $\mathsf{F}_4(\mathbb{K}, \mathbb{K})$ and then $\Gamma'$ is a subbuilding of type $\mathsf{F}_4(\mathbb{K}, \mathbb{K}^2)$. The latter is exactly the dual of the former (see Theorem 10.2 of [22]). Hence a symplecton in the latter is isomorphic to a parabolic quadric of type $\mathsf{B}_3$. Now, the residue of $\{p, L\}$ inside $Q$ consists of the set of planes through $L$ contained in a fixed parabolic subquadric of Witt index 3; hence in the 4-space $\mathsf{PG}(4, \mathbb{K})$, we have a set of lines in a hyperplane, so in a 3-space $V_Q$, defined by $Q$, as is apparent from the diagram. The structure of the residue of $\{p, L\}$ in $Q$ is visible in $V_Q$ through the Klein correspondence: The points of that residue correspond to lines of $V_Q$, and the points of the residue of $\{p, L, S\}$ in $\Gamma$ form a parabolic quadric of type $\mathsf{B}_1$, i.e., a conic, and consequently through the Klein correspondence this is a regulus $R$ in $V_Q$.

Now we restrict the projective plane $\Pi'$ to $\Gamma'$ and obtain a projective plane $\Pi$ represented in $\mathsf{PG}(4, \mathbb{K})$ as follows: its points are certain lines of $\mathsf{PG}(4, \mathbb{K})$ and its lines are certain reguli and partial planar line pencils. Since $V$ is not contained in $Q$, we see that the elements of $\Pi$ span $\mathsf{PG}(4, \mathbb{K})$. But then Lemma 7.8 implies that this leads to a contradiction. $\square$

We finally also rule out the case where all symplecta of $\Gamma$ are contained in quads of $\Delta$.

**Lemma 7.10** *No symplecton of $\Gamma$ is contained in a quad of $\Delta$.*

*Proof* Suppose that some symplecton is contained in a quad. Then, by Lemma 7.9, each symplecton corresponds to a quad. As in the third paragraph of the proof of Lemma 7.9, this gives rise to a projective plane $\Pi$ represented in $\mathsf{PG}(4, \mathbb{K})$ as follows. The points of $\Pi$ are lines of $\mathsf{PG}(4, \mathbb{K})$, and the lines of $\Pi$ are reguli. But then Lemmas 7.6 and 7.7 yield a contradiction. $\square$

We conclude that every $\Gamma$-symplecton is contained in a 5-space of $\Delta$. This concludes the proof of Proposition 7.4. As a consequence, we have the following corollary.

**Corollary 7.11** *Hyperbolic lines of $\Gamma$ are full lines of $\Delta$.*

*Proof* Each hyperbolic line of $\Gamma$ is contained in a symplecton. □

### 7.2 The relation of $\Gamma$ with the quads of $\Delta$

In this subsection we show that quads of $\Delta$ have only two possible positions with respect to the embedded geometry $\Gamma$, and the precise position depends on the geometric hyperplane of the quad induced by the point set of $\Gamma$. The analysis of the quads will be fundamental and crucial for the construction of the corresponding symplectic polarity.

We start with a lemma analysing the collinearity in $\Delta$ between points of $\Gamma$.

**Lemma 7.12** *Let $x$, $y$ be two distinct points of $\Gamma$. Then $x$ and $y$ are collinear in $\Delta$ if, and only if, they are collinear or symplectic in $\Gamma$.*

*Proof* In view of Corollary 7.11, we only have to prove that if $x$ and $y$ are neither collinear nor symplectic in $\Gamma$, then they are not collinear in $\Delta$.

First suppose that $x$ and $y$ are special. By Corollary 5.2, there exists a symplecton $S$ containing $x$ and $x \bowtie y$. Then $y$ is close to $S$, and so the set of points of $S$ collinear or symplectic to $y$ is a 3-space $W$. If $x$ were also collinear with $y$ in $\Delta$, then we would have found a 4-space (generated by $x$ and $W$) inside the 5-space $S$ collinear with $y$, contradicting Fact 4.10.

Now suppose that $x$ and $y$ are opposite points in $\Gamma$. We can select a symplecton $S$ through $x$ far from $y$. Let $u$ be the unique point of $S$ symplectic to $y$. Then all points of $S$ equal or $\Gamma$-collinear with $u$ form a 4-space $W \subseteq S$ and they are not $\Delta$-collinear with $y$ (except for $u$) by the first part of the proof. If $x$ and $y$ were $\Delta$-collinear, then by Fact 4.10, $y$ would be $\Delta$-collinear with all points of a 3-space of $S$, and hence with all points of a plane of $W$, a contradiction. □

**Lemma 7.13** *Let $x$ be an arbitrary point of $\Gamma$. Then the line pencil through $x$ in $\Gamma$ coincides with the line pencil through $x$ in $\Delta$ of a certain quad $Q$ through $x$, i.e. $x^{\perp_\Gamma} = x^{\perp_\Delta} \cap Q$. Moreover, there is only one quad $Q$ having this property.*

*Proof* Let $L_1, L_2$ be any two lines of $\Gamma$ through $x$ not contained in a symplecton. If $p_1 \in L_1 \setminus \{x\}$ and $p_2 \in L_2 \setminus \{x\}$, then $\{p_1, p_2\}$ is a special pair. By Lemma 7.12 and Fact 4.1, there is a unique quad $Q$ containing $p_1$ and $p_2$. By Fact 4.3, $Q$ contains $x$ and hence it contains $L_1$ and $L_2$. Let $L$ be any line of $\Gamma$ through $x$ such that $L$ and $L_i$, $i = 1, 2$, are contained in a symplecton of $\Gamma$. Then clearly $L_1$ and $L$ are contained in a plane of $\Delta$, just like $L_2$ and $L$. Using Fact 4.3 again, we see that $L$ has to be contained in $Q$.

Now note that Lemma 5.36 states that $\mathcal{N}_x$, furnished with all its subsets of lines intersecting an ordinary or a hyperbolic line contained in $x^{\perp_\Gamma} \setminus \{x\}$, has the structure of a polar space $P$ of type D$_4$. The above property now translates to the following. The lines $L_1, L_2$ correspond to non-collinear points $\ell_1, \ell_2$ of $P$, and $L$ to a point $\ell$ of $P$ collinear with both these points. The points of $P$ collinear with both $\ell_1$ and $\ell_2$ form a polar space of rank 3, and its set of points collinear with any other (fixed) point $z$ of $P$ forms a geometric hyperplane of it (or coincides with it). By Lemma 5.10, $z$ is collinear with two non-collinear points that are collinear with both $\ell_1, \ell_2$. Translated back to $\Delta$, this means that every line of $\Gamma$ through $x$ is indeed contained in $Q$.

Now, since planes of $\Gamma$ are full planes of $\Delta$ by Lemma 7.1, and since hyperbolic lines of $\Gamma$ are full lines of $\Delta$ by Corollary 7.11, one sees that the embedding of the polar space of type $\mathsf{D}_4$ obtained by considering the lines through $x$ in $\Gamma$, in the polar space of type $\mathsf{D}_4$ obtained by considering the lines of $\Delta$ through $x$ in $Q$ (which is a polar space of type $\mathsf{D}_5$), is full, and hence it is a bijection.

Also, by Fact 4.2 and Lemma 5.36, the quad containing $x^{\perp_\Gamma}$, i.e. $Q$, is unique. This completes the proof of the lemma. □

From now on we denote the unique quad $Q$, only depending on $x \in \Gamma$, of the previous lemma by $Q_x$.

For two points $x, y \in \Delta$ that are not $\Delta$-collinear, we denote by $Q(x, y)$ the unique quad containing $x, y$. Note that, if $x$ and $y$ are special in $\Gamma$, then $Q(x, y) = Q_{x \bowtie y}$.

**Lemma 7.14** *If $x$ and $y$ are two opposite points of $\Gamma$, then the quad $Q(x, y)$ intersects $\Gamma$ precisely in $\widehat{E}(x, y)$.*

*Proof* Every point of $\widehat{E}(x, y)$ belongs to $Q(x, y)$ by the definition of $\widehat{E}(x, y)$, Lemma 7.12 and Corollary 4.4. Now let $z$ belong to $\Gamma$ and to $Q(x, y)$. Since $\widehat{E}(p, q)$ is a polar space relative to its hyperbolic lines we can find two opposite (in $\Gamma$) points $p, q \in \widehat{E}(x, y)$ both $\Delta$-collinear with $z$, and so $z \in \widehat{E}(p, q) = \widehat{E}(x, y)$ by Proposition 5.22. □

An immediate consequence of the previous lemma is the following.

**Corollary 7.15** *Each extended equator geometry of $\Gamma$ is contained in a unique quad.*

**Lemma 7.16** *If $x$ is a point of $\Gamma$, then the quad $Q_x$ intersects $\Gamma$ precisely in $x^{\perp_\Gamma}$. In particular, for two points $u, v$ of $\Gamma$ we have $u \in Q_v$ if, and only if, $v \in Q_u$.*

*Proof* In view of Lemma 7.13, we only need to show that no point of $Q_x$ not $\Delta$-collinear with $x$ belongs to $\Gamma$. Suppose, by way of contradiction, that a point $u \in Q_x$ not $\Delta$-collinear with $x$ belongs to $\Gamma$. Since $Q_x$ is a polar space, every line of $\Gamma$ through $x$ contains a unique point $\Delta$-collinear with $u$. Hence, as $x$ is not symplectic to $u$, the pair $\{x, u\}$ is special by the last assertion of Lemma 5.7. So there is exactly one point $\Gamma$-collinear with both $x$ and $u$. Let $v$ be a point $\Gamma$-collinear with $x$ and symplectic to $u$. As $\mathcal{N}_x$ has the structure of a $\mathsf{D}_4$ (see Lemma 5.36), we can take a second point $w \in x^{\perp_\Gamma} \cap u^{\perp\!\!\!\perp}$ with $\{v, w\}$ a special pair (by selecting the line $xw$ opposite $xv$ and $w \neq x \bowtie u$). The symplecton $S(u, v)$ is far from $w$ since $\{v, w\}$ is a special pair and $x = w \bowtie v$ does not belong to $S(u, v)$ (as otherwise $x$ would be $\Delta$-collinear to $u$). But then $v$ should be $\Gamma$-collinear with $u$ by Fact 5.5, a contradiction.

The last assertion follows from the fact that $u \in Q_v$ if, and only if, $u \in v^{\perp_\Gamma}$. The latter is equivalent to $v \in u^{\perp_\Gamma}$ and finally to $v \in Q_u$. □

By the *dual* of $\Delta$, we mean the point-line geometry obtained from $\Delta$ as a building by switching the roles of points and quads, and of lines and 4-spaces (not contained in a 5-space, so elements of type 5, see Sect. 4.1). Hence the points of the dual of $\Delta$, as a point-line geometry, are the quads of $\Delta$, the lines are the 4-spaces not contained in a 5-space, and incidence is symmetrized containment. The principal of duality in buildings of type $\mathsf{E}_6$, as shown by the existence of a symplectic polarity, see Theorem 6.37 (see also [20]), implies that the dual of $\Delta$ is isomorphic to $\Delta$ itself.

**Lemma 7.17** *The map $x \mapsto Q_x$ defines a full embedding of $\Gamma$ in the dual of $\Delta$.*

*Proof* The injectivity of the map $x \mapsto Q_x$ follows from Lemma 7.16 and the fact that, for each point $x \in \Gamma$ and for each point $y \in x^{\perp_\Gamma}$, $y \neq x$, there exists a point $z \in x^{\perp_\Gamma}$ not collinear to $y$ (as the residue of $\Gamma$ in $x$ is not a linear space, but a dual polar space).

Let $L$ be a line of $\Gamma$. For $x, y \in L$, by Lemma 7.13, the quads $Q_x$ and $Q_y$ share the points of $\Gamma$ which are $\Gamma$-collinear with $L$. Furthermore, by Fact 4.2 and Lemma 6.19, $Q_x \cap Q_y$ coincides with the 4-space $W = L^{\perp_\Gamma}$. So, $W \subseteq z^{\perp_\Gamma} \subseteq Q_z$, for each $z \in L$. This shows that we have an embedding. Now we show fullness.

Let $Q$ be any quad containing $W$, $Q \neq Q_x$. We show that $Q = Q_t$ for some $t \in L$, thus completing the proof of the lemma. Let $u \in \Gamma$ be a point of $Q_x$ that is $\Gamma$-collinear with $x$ and symplectic to $y$ in $\Gamma$. Then clearly each plane of $\Gamma$ through $L$ inside the symplecton $S(y, u)$ contains a line of points $\Gamma$-collinear with $u$; so $Q_u$ intersects $W$ in at least two points, hence, again by Fact 4.2, $Q_u \cap Q$ is a 4-space $V$. We see that $V \neq W$, as otherwise all points of $W$ would be $\Gamma$-collinear with $u$; but $y \notin u^{\perp_\Gamma}$, a contradiction. So $V$ and $W$ are distinct 4-spaces inside the quad $Q$ having at least two points in common. Since they both occur as the intersection of quads, they belong to the same system of maximal singular subspaces of $Q$ (they have both type 5 in the building of type E₆, with Bourbaki labelling) and hence they intersect in a plane and so there is some point $v \in V \setminus W$ collinear with $u$.

Then $uv$ is a line of $\Gamma$ by Lemma 7.13. Since $v \notin W$, there is some point $w \in W$ not $\Delta$-collinear with $v$. But since $w \in W$, the point $w$ belongs to $\Gamma$. By Lemma 7.12, $w$ is either special to or opposite $v$. But since $v \in Q$ and $L \subseteq Q$, at least one point of $L$ is $\Delta$-collinear with $v$. It follows that $w$ is not opposite $v$. Hence $w$ is special to $v$. But then, if $t = v \bowtie w$, since $Q_t$ is the unique quad through $w$ and $v$, we have $Q_t = Q$. Lemma 7.16 implies that $t$ is $\Delta$-collinear with all points of $W$, hence $t \in W$. But $t$ is also $\Gamma$-collinear with all points of $W$. We claim that this implies that $t \in L$. Indeed, if not, then consider a plane $\alpha$ of $\Gamma$ through $L$ not containing $t$. Then $t$ is collinear with all points of $\alpha$, contradicting Fact 5.5 by including $\alpha$ in a symplecton. The claim and the lemma are proved. □

We will call the above embedding of $\Gamma$ in the dual of $\Delta$ briefly *the dual embedding*. This gives us yet another embedding of $\Gamma$ in a building of type E₆. Note that, as a consequence, Lemmas 7.12, 7.13, 7.14 and 7.16, and Corollary 7.15, also hold for this embedding, the only difference being that we have denoted the points $x$ of $\Gamma$ as $Q_x$ instead of just $x$ when seen in $\Delta$.

**Lemma 7.18** *Let $x$ and $y$ be opposite points of $\Gamma$. Then $Q_x \cap Q_y = \{s\}$, with $s$ a point of $\Delta \setminus \Gamma$, and $s \in Q_u$ for each point $u \in Q(x, y) \cap \Gamma$, i.e., $\{s\} = \bigcap_{u \in \widehat{E}(x,y)} Q_u$.*

*Proof* By Lemmas 7.12 and 7.17, $Q_x$ and $Q_y$ are non-collinear points in the dual embedding. Denote the unique quad through them in the dual embedding by $Q^s$, where, dually, $s$ is the unique point in $Q_x \cap Q_y$. Moreover, $Q^s$ contains all points $Q_t$ with $Q_t \in \widehat{E}(Q_x, Q_y)$ by Lemma 7.14, i.e. all points $Q_t$ with $t \in \widehat{E}(x, y)$, as this is a notion in $\Gamma$ independent of any embedding (in particular $\widehat{E}(Q_x, Q_y) = \{Q_z : z \in \widehat{E}(x, y)\}$). Note that the points in the quad $Q^s$ correspond to all quads through $s$, so in particular, $s \in Q_u$ for $u \in Q(x, y) \cap \Gamma = \widehat{E}(x, y)$, the latter equality also by Lemma 7.14. The assertion follows. □

**Definition 7.19** (*Tangent and Secant Quads*) If $x \in \Gamma$, then $Q_x$ is called a *tangent quad*. If $x, y \in \Gamma$ are opposite, then $Q(x, y)$ is called a *secant quad*.

**Lemma 7.20** *Every quad of $\Delta$ is either tangent or secant.*

*Proof* We first show the assertion for quads containing at least one point of $\Gamma$. So let $Q$ be a quad and $x \in Q$ a point of $\Gamma$. By Lemma 7.13 applied to the dual embedding, cf. Lemma 7.17, every quad meeting $Q_x$ in a 4-space through $x$ is of the form $Q_y$, where $y \in \Gamma$ is a point $\Gamma$-collinear with $x$. Hence, by Fact 4.2, we may assume that $Q$ intersects $Q_x$ in just $\{x\}$. We show that in this case $Q$ is a secant quad.

Towards this, we first show how to select a special pair $\{y_1, y_2\} \subseteq x^{\perp_\Gamma}$ such that $Q_{y_1} \cap Q$ is a 4-space $W_i$ through $x$. Dualizing the latter condition (hence considering the dual embedding) yields points $y_1$ and $y_2$ that are both $\Delta$-collinear with a point $q$ (corresponding to $Q$) in $Q_x$ (as $x \in Q$) with $q$ and $x$ not $\Delta$-collinear (as $Q_x \cap Q$ is not a 4-space). Hence we select $y_1$ and $y_2$ as points of $\Gamma$ both collinear to $x$, special to each other and both $\Delta$-collinear with $q$. Since $y_i \in Q_x$, we have $y_i \notin Q$, $i = 1, 2$, and so we see that, inside $Q_{y_i}$, $y_i$ is $\Gamma$-collinear with all points of a 3-space $V_i \subseteq W_i$ of $\Delta$. From $Q_{y_1} \cap Q_{y_2} = \{x\}$ and $y_1 \perp_\Gamma x \perp_\Gamma y_2$ we infer that $V_1 \cap V_2 = \{x\}$. All points of $V_1 \cup V_2$ belong to $\Gamma$, by Lemma 7.13. The subspaces $V_1$ and $V_2$ cannot be contained in a common subspace, as that subspace would have dimension at least 6, a contradiction. So we can select points $u_i \in V_i$, $i = 1, 2$, with $u_1$ not $\Delta$-collinear with $u_2$. But $u_1, u_2 \in Q \cap \Gamma$, and so $Q = Q(u_1, u_2)$. Note that this is necessarily a secant quad, as otherwise $u_1$ and $u_2$ would be special and $u_1 \bowtie u_2$ would be $\Gamma$-collinear with all points of $V_1 \cup V_2$, in particular, with $x$ and hence contained in $Q_x$. But a point of $Q_x$ is only $\Gamma$-collinear with at most a 3-space of $Q$, a contradiction.

Now we show that every quad contains at least one point of $\Gamma$. Suppose for a contradiction that the quad $Q$ does not contain any point of $\Gamma$. Since the graph on the quads with adjacency defined as 'intersecting in a 4-space' is connected (it has diameter 2, which is seen from Lemma 7.17 and Fact 4.1), we may assume that $Q$ is adjacent with a quad of the form $Q_x$ or $Q(x, y)$. In both cases, the intersection with $\Gamma$ is a geometric hyperplane, and the intersection with $Q$ is a 4-space. Hence the result. □

### 7.3 Construction of the polarity of $\Delta$ fixing $\Gamma$; end of the proof of Theorem 3

By the dual of Lemma 7.20, every point $p$ of $\Delta$ either belongs to $\Gamma$ or is the only point incident with the quads $Q_x$ for each $x$ in a certain extended equator geometry $\widehat{E}$, i.e. $\{p\} = \bigcap_{x \in \widehat{E}} Q_x$. Note that, if $\widehat{E} = \widehat{E}(x, y)$, then $p$ is determined as $Q_x \cap Q_y$ (cf. Lemma 7.18). We say that $p$ is the point *associated with* $\widehat{E}$ and vice versa. The forgoing arguments and Lemma 7.18 show that this association is a bijection between the complement of $\Gamma$ in $\Delta$ and the set $\mathscr{E}$ of all extended equator geometries.

**Definition 7.21** (*Polarity*) We now define a map $\theta$ from the point set of $\Delta$ to the set of quads of $\Delta$ and from the set of quads of $\Delta$ to the set of points of $\Delta$ as follows. We define the image under $\theta$ of a point $x$ of $\Gamma$ as the quad $Q_x$ and vice versa; the image under $\theta$ of the point associated (cf. the previous paragraph) with the extended equator geometry $\widehat{E}(x, y)$ is defined as the quad $Q(x, y)$ and vice versa.

The next proposition completes the proof of Theorem 3. We need the following lemma for its proof.

**Lemma 7.22** *Let $Q$ be a quadric of type $\mathsf{D}_5$ and let $E$ be a subquadric of type $\mathsf{B}_4$ obtained in the standard way by intersecting an embedding of $Q$ in $\mathsf{PG}(9, \mathbb{K})$ with a hyperplane. Then there is a (unique) natural pairing $\sigma$ of the points of $Q \setminus E$ such that, if $x \in Q \setminus E$, and $V$, $W$ are two 4-spaces of $Q$ of the same type through $x$ intersecting in just $x$, then $x^\sigma = V' \cap W'$, where $V'$ and $W'$ are the unique 4-spaces in $Q$ of the opposite type through $V \cap E$ and $W \cap E$, respectively.*

*Proof* If $Q$ has equation $X_{-1}X_1 + X_{-2}X_2 + X_{-3}X_3 + X_{-4}X_4 + X_{-5}X_5 = 0$, and $E$ is obtained by intersecting $Q$ with the hyperplane with equation $X_{-5} = X_5$, then the pairing $\sigma$ is given by mapping $(x_{-5}, x_{-4}, x_{-3}, x_{-2}, x_{-1}, x_1, x_2, x_3, x_4, x_5)$ to $(x_5, x_{-4}, x_{-3}, \ldots, x_3, x_4, x_{-5})$. It is easy to show and calculate that this pairing satisfies the conditions, and that it is unique. □

**Proposition 7.23** *The map $\theta$ is a polarity of $\Delta$ whose absolute structure is $\Gamma$.*

*Proof* We first show that $\theta$ preserves incidence among the points and the quads. Let $x$ be a point contained in a quad $Q$ of $\Delta$. There are four cases.

(i) $x \in \Gamma$ and $Q$ is tangent, say $Q = Q_y$, $y \in \Gamma$. In this case, $x$ and $y$ are $\Gamma$-collinear as $x \in Q_y \cap \Gamma = y^{\perp_\Gamma}$ by Lemma 7.16, and hence the same lemma implies $Q_y^\theta = y \in Q_x = x^\theta$.

(ii) $x \in \Gamma$ and $Q$ is secant, say $Q = Q(y, z)$, with $y, z$ opposite points of $\Gamma$. By Lemma 7.14, $x \in Q \cap \Gamma = \widehat{E}(y, z)$. So, $Q^\theta \in x^\theta = Q_x$ by the definition of $Q^\theta$.

(iii) $x \in \Delta \setminus \Gamma$ and $Q$ is tangent. The dual of the previous case holds.

(iv) $x \in \Delta \setminus \Gamma$ and let $Q = Q(y, z)$ be a secant quad, with $y, z$ opposite points of $\Gamma$. Let $t$ be $Q^\theta$. Then $\{t\} = Q_y \cap Q_z = \bigcap_{u \in \widehat{E}(y,z)} Q_u$. For each $u \in \widehat{E}(y, z)$, we have $t \in Q_u$ and as $t \notin \Gamma$ by Lemma 7.18, it follows from Lemma 7.13 that $t$ and $u$ are not $\Delta$-collinear. Hence, $t$ is opposite $Q$, as otherwise Fact 4.3 would imply that $t$ is $\Delta$-collinear with a point of $\widehat{E}(y, z)$, which is not the case.

We may choose $y$ and $z$ in $Q$ such that they are both $\Delta$-collinear with $x$. Consider the pairing $\sigma$ of Lemma 7.22. Then also $x^\sigma$ is $\Delta$-collinear with both $y, z$. By Fact 4.3, $Q_y \cap (x^\sigma)^{\perp_\Delta}$ is a 4′-space $V_y$ of $Q_y$. The point $t$ belongs to $Q_y$ (as $y \in Q$), but $t \notin V_y$ as $t$ is not $\Delta$-collinear with $x^\sigma$, being opposite $Q$. So $t^{\perp_\Delta} \cap V_y$ is a 3-space $W_y$ contained in $Q(t, x^\sigma)$ by Corollary 4.4. Hence $Q(t, x^\sigma) \cap Q_y$ is the 4-space $U_y$ generated by $t$ and $W_y$. The space $V_y$ and the point $x^\sigma$ are contained in a unique 5-space $T_y$. Likewise, $Q(t, x^\sigma) \cap Q_z$ is a 4-space $U_z$ containing $t$, and $W_z, V_z, T_z$ are defined analogously. Consider an arbitrary point $u \in W_y$. Since $u \in Q_y \cap y^{\perp_\Delta}$, Lemma 7.13 implies $u \in \Gamma$. Since the $Q_w$, $w \in \Gamma$, form a dual embedding, there is a unique point $w$ on the line $uy$ such that $Q_w$ intersects $Q$ in a 4-space $Y$ (this follows from the fact that, in the dual of $\Delta$, $y$ is a polar space with point set the original quads through $y$ and line set the original 4-spaces containing $y$, with natural incidence). The subspace $Y$ does not contain $x^\sigma$ as $x^\sigma$ is $\Delta$-collinear with $w$ and, by Lemma 7.13, this would force $x^\sigma$ to belong to $\Gamma$.

Now we claim that $w = u$. Indeed, suppose not. As $w \in V_y \setminus W_y$, the points $w$ and $t$ are not $\Delta$-collinear. The subspace $Y$ intersects $\Gamma$ in a 3-space $Y'$ consisting of the points of $Y$ collinear with $w$ by Lemma 7.16. Let $a$ be a point in $Y'$ distinct from $y$. Since $a \in Q \cap Q_w$, we know that $Q_a$ contains $Q^\theta = t$ and $w$. However, this implies that $Q_a$ coincides with $Q(w, t)$, which is equal to $Q_y$, a contradiction. Hence $w = u$.

As $u$ is $\Delta$-collinear to $Y'$ and $x^\sigma$, it is $\Delta$-collinear to the 4′-space $Y''$ of $Q$ generated by $x^\sigma$ and $Y'$. Likewise, we can select a point $v \in W_z$, and choose $u, v$ not $\Delta$-collinear (this is possible as otherwise $W_z$ and $W_y$ would generate a 6-space in $\Delta$). There is a corresponding 4-space $Z = Q \cap Q_v$ in $Q$ intersecting $\Gamma$ in a 3-space $Z'$ which, together with $x^\sigma$, generates a 4′-space $Z''$. Now suppose $Y' \cap Z'$ is not empty. Then $Y''$ and $Z''$, which by definition also contain $x^\sigma$, are two 4′-spaces intersecting in at least a line, hence $Y'' \cap Z''$ is at least a plane; so $Y' \cap Z'$ is at least a (hyperbolic) line, implying that $u$ and $v$ are $\Delta$-collinear, a contradiction. Hence $Y'$ and $Z'$ are disjoint and it follows from Lemma 7.22 that $Z$ and $Y$ intersect in $x$. Hence $x \in Q_u \cap Q_v$. From (iii) it follows that $x^\theta$ contains $u$ and $v$ and therefore, as those points are not collinear, $x^\theta = Q(u, v)$. But then $t \in x^\theta$ as $Q(u, v) = Q(t, x^\sigma)$. So $Q^\theta \in x^\theta$, as required.

It is trivial to check that $\theta$ has order 2. Now $\Delta$-collinear points are mapped onto quads that share a 4-space, since, if two points are $\Delta$-collinear they are contained in at least two quads, and so the images contain at least two points. It follows that $\theta$ is a collineation from $\Delta$ onto its dual, hence a polarity. Clearly $\Gamma$ is its absolute structure. □

*Proof of Theorem 3.* Let $\Gamma$ and $\Delta$ be as in the theorem. Then by Corollary 7.2 they are defined over the same field $\mathbb{K}$. By Proposition 7.23, $\Gamma$ arises from a polarity of $\Delta$, which is symplectic by definition (since $\Gamma$ has symplectic residues). □

**Acknowledgements** The work of the second author was carried when he was at Indian Statistical Institute, Bangalore centre, and during his visits to Department of Mathematics, Ghent University. He thanks both the institutions for extending their kind hospitality and excellent working conditions.

# Index of Symbols

| Symbol | Description |
|---|---|
| $x^\perp$ | The "perp" of the point $x$: all points equal or collinear to $x$ |
| $\Delta$ | Building of type $\mathsf{E}_6$ or it natural point-line geometry |
| $\theta$ | A symplectic polarity in $\Delta$ |
| $\Gamma$ | Building of type $\mathsf{F}_4$ or the corresponding symplectic metasymplectic parapolar space |
| $x \perp y$ | The point $x$ is collinear to the point $y$ |
| $x \perp\!\!\!\perp y$ | The point $x$ is symplectic to the point $y$ |
| $x \Diamond y$ | The unique symplecton through the symplectic points $x$ and $y$ |
| $x \bowtie y$ | The unique point collinear to both $x$ and $y$ when $\{x, y\}$ is a special pair |
| $x^{\perp\!\!\!\perp}$ | All points equal or symplectic to the point $x$ |
| $h(x, y)$ | The hyperbolic line containing the symplectic pair $\{x, y\}$ of points |
| $S(h)$ | The unique symplecton containing the hyperbolic line $h$ |
| $\mathscr{S}_p$ | The family of symplecta containing the point $p$ |
| $E(p, q)$ | The equator geometry of the pair $\{p, q\}$ of opposite points |
| $\widehat{E} = \widehat{E}(p, q)$ | The extended equator geometry of the pair $\{p, q\}$ of opposite points |
| $\widehat{T} = \widehat{T}(p, q)$ | The tropic circle geometry of the pair $\{p, q\}$ of opposite points |
| $\beta(x)$ | The unique hyperbolic solid in $\widehat{E}(p, q)$ collinear to $x \in \widehat{T}(p, q)$ |
| $\beta(U)$ | The unique point collinear to the hyperbolic solid $U$ |
| $\Theta(\widehat{T}(p, q))$ | The imaginary completion of $\widehat{T}(p, q)$ to a half spin $\mathsf{D}_5$ |
| $\widehat{H}(p, q)$ | The set of point collinear or equal to at least one point of $\widehat{E}(p, q)$ |
| $\mathscr{N}_x$ | The set of lines of $\Gamma$ through the point $x$ |
| $\mathsf{D}_4(\mathscr{N}_x)$ | The point-line geometry of type $\mathsf{D}_4$ defined on $\mathscr{N}_x$ |
| $\mathscr{P}$ | The point set of the point-line $\mathsf{E}_6$-geometry defined from $\Gamma$ |
| $\mathscr{L}$ | The line set of the point-line $\mathsf{E}_6$-geometry defined from $\Gamma$ |
| $\mathscr{E}$ | The family of new points of $(\mathscr{P}, \mathscr{L})$, i.e., the family of extended equator geometries of $\Gamma$ |
| $\mathscr{F}$ | The family of new lines of $(\mathscr{P}, \mathscr{L})$, i.e., those containing members of $\mathscr{E}$ |
| $T_\mathfrak{e}$ | The tropic circle geometry of the extended equator geometry $\mathfrak{e}$ |
| $\Sigma(p)$ | The quad of $(\mathscr{P}, \mathscr{L})$ corresponding to the point $x$ |
| $\Sigma(\widehat{E}(p, q))$ | The quad of $(\mathscr{P}, \mathscr{L})$ corresponding to the new point $\widehat{E}(p, q)$ |
| $\mathscr{Q}$ | The family of quads of $(\mathscr{P}, \mathscr{L})$ |
| $\mathscr{U}$ | The family of maximal singular 4-spaces of $(\mathscr{P}, \mathscr{L})$ |
| $U(L)$ | The projective 4-space associated to the line $L$ of $(\mathscr{P}, \mathscr{L})$ |
| $V^+, V^-$ | Twin hyperbolic cones |
| $\mathscr{M}$ | The family of singular 5-spaces of $(\mathscr{P}, \mathscr{L})$ |
| $\mathscr{T}$ | The family of singular planes of $(\mathscr{P}, \mathscr{L})$ |
| $\mathfrak{E}$ | The geometry of type $\mathsf{E}_6$ defined from $\Gamma$ |
| $*$ | The incidence relation of $\mathfrak{E}$ |
| $\mathscr{U}(\Sigma)$ | The subset of elements of $\mathscr{U}$ incident with the quad $\Sigma$ |
| $\mathscr{M}(\Sigma)$ | The set of 4-spaces of the quad $\Sigma$ obtained by intersecting $\Sigma$ with the members of $\mathscr{M}$ that are incident with $\Sigma$ |
| $x^{\perp_\Gamma}, x^{\perp_\Delta}$ | The perp of $x$ in $\Gamma$ and $\Delta$, respectively |
| $Q_x$ | The unique quad in $\Delta$ containing all lines of $\Gamma$ through $x$ |
| $Q(x,y)$ | The unique quad of $\Delta$ containing the non-collinear points $x$ and $y$ |

## Index of Notions